\documentclass{article}

\usepackage{amsmath}
\usepackage{amssymb}
\usepackage{amsthm}
\usepackage{amscd}

\newtheorem{theorem}{Theorem}[section]
\newtheorem{lemma}[theorem]{Lemma}
\newtheorem{conjecture}[theorem]{Conjecture}
\newtheorem{proposition}[theorem]{Proposition}

\theoremstyle{definition}

\newtheorem{example}[theorem]{Example}

\theoremstyle{remark}
\newtheorem{remark}[theorem]{Remark}

\numberwithin{equation}{section}

\makeatletter
\def\eqnarray{%
  \stepcounter{equation}%
  \let\@currentlabel=\theequation
  \global\@eqnswtrue
  \global\@eqcnt\z@
  \tabskip\@centering
  \let\\=\@eqncr
  $$\halign to \displaywidth\bgroup\@eqnsel\hskip\@centering
  $\displaystyle\tabskip\z@{##}$&\global\@eqcnt\@ne
  \hfil$\displaystyle{{}##{}}$\hfil
  &\global\@eqcnt\tw@$\displaystyle\tabskip\z@{##}$\hfil
  \tabskip\@centering&\llap{##}\tabskip\z@\cr}
\makeatother

\newfont{\germ}{eufm10}
\newfont{\slsmall}{cmsl8}

\def\Aff{\mbox{\sl Aff}}
\def\bt{\tilde{b}}
\def\cd{\cdots}
\def\et#1{\tilde{e}_{#1}}
\def\ft#1{\tilde{f}_{#1}}
\def\geh{\goth{g}}

\def\goth#1{\mbox{\germ #1}}
\def\La{\Lambda}
\def\la{\lambda}
\def\Lab{\overline{\La}}
\def\ol#1{\overline{#1}}
\def\ot{\otimes}

\def\qed{~\rule{1mm}{2.5mm}}

\def\Uqp{U'_q(\geh)}
\def\veps{\varepsilon}
\def\vphi{\varphi}
\def\wt{\mbox{\sl wt}\,}
\def\wts{\mbox{\slsmall wt}\,}
\def\C{{\mathbb C}}
\def\Z{{\mathbb Z}}

\def\geh{{\mathfrak g}}
\def\bb{\boldsymbol{b}}
\def\bp{\boldsymbol{p}}

\def\Aff{\mbox{\sl Aff}}
\def\et#1{\tilde{e}_{#1}}
\def\ft#1{\tilde{f}_{#1}}
\def\id{\mbox{\sl id}\,}

\def\lev{\mbox{\sl lev}\,}
\def\ch{\mbox{\sl ch}\,}
\def\P{{\mathcal P}}
\def\simarrow{\overset{\sim}{\rightarrow}}

\title{Paths, Crystals and 
Fermionic Formulae}

\author{
G. Hatayama\thanks{
Institute of Physics, University of Tokyo, Komaba, Tokyo 153-8902, Japan},
A. Kuniba,$\hspace{-1.2mm}^*$
M. Okado\thanks{Department of Informatics and Mathematical Science,
Graduate School of Engineering Science,
Osaka University,
Toyonaka, Osaka 560-8531,
Japan},
T. Takagi\thanks{
Department of  Applied Physics, National Defense Academy,
Yokosuka 239-8686, Japan} 
 and Z. Tsuboi$\hspace{-0.1mm}^*$}
\date{}
\begin{document}

\maketitle
\begin{center}
{\it Dedicated to Professor Barry McCoy on the occasion
of his sixtieth birthday}
\end{center}

\begin{abstract}
We introduce a fermionic formula associated 
with any quantum affine algebra $U_q(X^{(r)}_N)$.
Guided by the interplay between corner transfer matrix and Bethe ansatz
in solvable lattice models, we study several aspects related to 
representation theory, most crucially, the crystal basis theory.
They include  one dimensional sums over both finite and semi-infinite paths,
spinon character formulae, 
Lepowski-Primc type conjectural formula for vacuum string functions,
dilogarithm identities, 
$Q$-systems and their solution by characters of various 
classical subalgebras and so forth.
The results expand \cite{HKOTY1}
including the twisted cases and more details on 
inhomogeneous paths consisting of non-perfect crystals.
As a most intriguing example, certain inhomogeneous one dimensional sums 
conjecturally give rise to  branching functions of an integrable 
$G^{(1)}_2$-module related to the embedding 
$G^{(1)}_2 \hookrightarrow B^{(1)}_3 \hookrightarrow D^{(1)}_4$.
\end{abstract}


\section{Introduction}
Fermionic formulae are certain $q$-polynomials or $q$-series 
expressed as a specific sum of products of $q$-binomial coefficients.
In this paper, which is an extended version of \cite{HKOTY1}, 
we introduce the fermionic formulae associated with any 
quantum affine algebra and discuss several aspects related to 
integrable systems and representation theory.
For twisted cases, they are new even at $q=1$.
Before going into the details, 
let us venture a brief historical account of the 
results scattering over immense literatures.

The first fermionic formula ($q=1$) appeared in eq.(45) of  Bethe's paper 
\cite{Be} in 1931. 
It was given as the number of solutions to the Bethe 
equation for the Heisenberg spin chain based on the string hypothesis.
Subsequent developments took place mainly {}from the 80's.
Here we divide them into the three main streams as arranged in 
Sections \ref{subsec:sb}--\ref{subsec:ctm} below.

\subsection{Quasi-particle picture}\label{subsec:sb}

Probably it was the Stony Brook group 
who first used the physics terminology 
{\em fermionic} for affine Lie algebra characters. 
Based on the Bethe ansatz analysis of integrable 3-state Potts chain \cite{KM},
Kedem and McCoy discovered a novel fermionic expression of 
the branching functions for the coset pair of the affine Lie algebras 
$(A^{(1)}_3)_1\oplus (A^{(1)}_3)_1/(A^{(1)}_3)_2$.
There was no minus sign in the sum as opposed to the well known 
Feigin-Fuchs-Rocha-Caridi bosonic formula for Virasoro characters.
This outstanding feature was due to the fermionic 
nature of the lowlying excitations characterized by 
the roots of the Bethe type equation.
Based on an elaborate analysis of the fermionic selection rule on the roots, 
they actually went beyond the characters of chiral blocks, and succeeded 
in constructing the modular invariant partition functions, i.e.,
bilinear combinations thereof, as a fermionic $q$-series.

Focusing on the constituent chiral blocks, 
Kedem {\em et al.} subsequently proposed the fermionic formulae 
for a variety of conformal field theory (CFT) models  
including the non-unitary Virasoro minimal ones
\cite{KKMM1,KKMM2}. 
Through these studies,
the quasi-particle picture was put forward as a general concept of
lowlying excitations leading to fermionic formulae \cite{DKKMM},
the Bailey pair technique was explored \cite{BMSW} and 
the Rogers-Ramanujan type identities were interpreted as 
an equivalence of bosonic and fermionic descriptions 
\cite{KMM}.
We refer to \cite{M} for the original survey.

The progress by the Stony Brook group inspired 
further developments {}from the viewpoint of 
Yangian symmetry and the  ``spinon basis" 
in integrable $A^{(1)}_1$ modules \cite{BPS,BLS1,BLS2}.
These ideas led to a proof of fermionic formulae by means of 
the chiral vertex operators and brought a new description of the  CFT 
Hilbert spaces.
We remark that most of the formulae in these works 
are level truncated (or restricted)
in the sense explained in Section \ref{subsec:ctm}.
See also \cite{ANOT,NY2,NY3} for more aspects concerning 
the spinon interpretations.

In the early stage of these developments, 
main interests were directed to  $q$-series rather than 
$q$-polynomials. 
The idea of ``finitizing" the former to the latter came with \cite{Me} 
as an ultraviolet cutoff in the quasi-particle picture.
It was followed by many papers, most typically 
\cite{Ber,Wa,Sc,BMS}, where the fermionic formulae were established  
for a polynomial analogue of $A^{(1)}_1$ coset Virasoro characters
including the fractional level cases. 
See also \cite{FW1,FW2}.
An essential question in finitization is the meaning of the 
consequent $q$-polynomials.
We postpone the discussion on this point to Section \ref{subsec:ctm}.

\subsection{Combinatorics of Bethe ansatz}
\label{subsec:kr}

The advances achieved in the 90's sketched in Section \ref{subsec:sb} 
may be viewed as a renaissance in the history of fermionic formulae 
linked up with  CFT.
Actually, fermionic formulae 
had been obtained in the middle 80's in a different context 
by Kerov, Kirillov and Reshetikhin \cite{KKR,KR1} 
for the Kostka polynomial.

The Kostka polynomial is a $q$-analogue of the 
multiplicity of the irreducible $\goth{sl}_{\,n}$-module $V_{\la}$
in the $m$-fold tensor product $V_{(\mu_1)}\ot V_{(\mu_2)}\ot\cd
\ot V_{(\mu_m)}$. Here for a partition $\la=(\la_1,\cd,\la_n)$
($\la_1\ge\cd\ge\la_n\ge0$) $V_{\la}$ stands for the irreducible 
$\goth{sl}_{\,n}$-module with highest weight $\sum_{i=1}^{n-1}
(\la_i-\la_{i+1})\La_i$ with $\La_i$ being the fundamental weights
of $\goth{sl}_{\,n}$. In particular, $V_{(\mu_i)}$ is the symmetric
tensor representation of degree $\mu_i$. 
The Kostka polynomial plays an important role in algebraic combinatorics, and
is generalized to the cases where the components of the tensor product are
not necessarily symmetric tensor representations.
See e.g. \cite{LLT,LT,S,SW}.

As the Stony Brook group, 
the approach by Kerov {\em et al.} was also based on 
the Bethe ansatz\footnote{Their Bethe equation is different {}from 
the one in \cite{KM}.}, 
but it had several distinctive features.
Postulating the completeness of string hypothesis,
they derived the $q=1$ fermionic formula as the multiplicity of 
irreducible $A_n$-modules in the tensor product of 
finite-dimensional Yangian $Y(A_n)$-modules.
Their formula admitted a natural $q$-analogue matching  
the Kostka polynomial.
In this context, they treated $q$-polynomials 
without embarking on $q$-series and finitization.
Later, Kirillov and Reshetikhin pushed these ideas further to 
extend the $q=1$ fermionic formula {}from 
$A_n$ to arbitrary classical simple Lie algebra $X_n$ \cite{KR2}.
There was no level truncation in those formulae 
(``classically restricted" in our terminology),
and passage to the infinite series and 
connection to $A^{(1)}_n$ branching functions 
came afterward in 1995 by Kirillov as his Conjecture 4 in  \cite{Ki4}.

The studies \cite{Ki1,Ki2,KR2} on the completeness of the string hypothesis
are relevant to polynomial fermionic formula 
at $q=1$ with no level truncation.
In this sense they are most natural continuation of 
Bethe's original treatise.
The invention of 
rigged configurations, bijection to 
semistandard Young tableaux \cite{KKR,KR1} and 
the Kirillov-Reshetikhin conjecture 
\cite{KR2} (see also Conjecture \ref{conj:fd-module})
was one of the earliest applications of the Bethe ansatz 
to combinatorics and representation theory.

When inferring the $q$-analogue at the level of 
fermionic formulae, the Kostka polynomial 
was the crucial guide for the $A_n$ case.
To seek the clue for the other cases is a 
natural question which was also raised 
in the end of Section \ref{subsec:sb}.
More broadly, one may ask;  
what is the meaning, definition or characterization of the 
$q$-polynomial that should be represented by a 
fermionic formula?

\subsection{Corner transfer matrix and crystals}\label{subsec:ctm}

The answer to the question had been prepared 
in Baxter's corner transfer matrix (CTM) method \cite{B},
which is another fundamental tool in solvable lattice models.
The method introduces the quantity called 
{\em one dimensional sums} (1dsums).
They are $q$-polynomials growing with lattice size, which eventually 
tend to affine Lie algebra characters in the infinite lattice limit.
Such a phenomenon was first worked out for the celebrated 
ABF model \cite{ABF,FB} and subsequently for other 
restricted solid-on-solid (RSOS) \cite{DJKMO1,JMO} 
as well as vertex models \cite{DJKMO2}.
In these works in the 80's, identification of 1dsums with characters 
was done mainly through explicit formulae (mostly bosonic one), therefore 
it was hard to go much beyond tractable examples.

The advent of Kashiwara's crystal basis theory \cite{Ka1} and its applications 
to vertex \cite{KMN1} and RSOS models \cite{JMMO,DJO} 
swept the difficulty away.
The theory provided an intrinsic definition of 
1dsums as well as a conceptual proof  
that their infinite lattice limit
must be affine Lie algebra characters.
The foregoing curiosity on 1dsums was resolved by 
the fundamental theorem \cite{KMN1}:
crystals of integrable highest weight modules are isomorphic to  
semi-infinite tensor products of finite crystals.
This propaganda is simply called {\em path realization}.

The crystal basis theory is a theory of quantum groups at $q=0$ 
having a source in CTM with numerous applications.
Here we just mention two examples relevant to this paper.
The first one is the  proof of the 
Kirillov conjecture on branching functions 
by Nakayashiki and Yamada \cite{NY1}.
It was done by introducing the 
{\em classically restricted} 1dsums as an intrinsic 
characterization of the Kostka polynomial in terms of crystals.
They are intermediate ones  
between the level restricted ones for RSOS and 
the unrestricted ones for vertex models, and will be essential to our 
discussion on fermionic formulae without level truncation.
The second example is the study of Demazure crystals.
Recall that Demazure modules are certain finite-dimensional subspaces of 
integrable highest weight modules labeled by Weyl group elements.
Demazure modules at $q=0$ are called  Demazure crystals \cite{Ka2}.
It was shown \cite{KMOTU} that under the path realization, 
the restriction of 1dsums to a certain Demazure crystal 
agrees with the finitization in the CTM language.
See Table 1 in \cite{KMOTU} for a summary of distinctive features of the
three kinds of paths;  unrestricted, classically restricted and 
(level) restricted ones.
Our fermionic formulae $M_\infty$ and $M_l$ in this paper 
correspond to the latter two, respectively.

\subsection{Physical combinatorics}\label{subsec:pc}

Throughout the story so far, it should be emphasized that 
the two essential ingredients in solvable lattice models, 
Bethe ansatz and CTM, play complementary roles.
The CTM and crystal theory provide an intrinsic definition 
of 1dsums $X$ and characterization as 
branching coefficients of integrable highest weight modules 
($q$-series case) or Demazure modules (finitized $q$-polynomial case).
On the other hand, the Bethe ansatz offers a specific formula $M$ 
which is fermionic.
Therefore a synthesis of the two ideas reads
\[
X = M. 
\]
To prove and understand this equality in the maximal generality 
is a fundamental problem in ``Physical Combinatorics"\footnote{
Named by Jean-Yves Thibon.} in the beginning of the 21st century.

In \cite{HKOTY1} we have launched a fermionic formula 
$M$ for general nontwisted affine Lie algebras such that 
the above identity conjecturally holds.
It is a $q$-analogue of \cite{KR2}, which covers a 
large list of known 
fermionic formulae under suitable specializations.
In this paper we extend our proposal including the twisted cases where 
none was known even at $q=1$.
Our generalization is not only in the direction of algebras.
The earlier studies on 1dsums \cite{KMN1,KMOTU} were done
exclusively for perfect crystals and homogeneous paths
(cf. Section \ref{subsec:contents}). 
Thus it remains a challenge to characterize the
infinite lattice limit of 1dsums over inhomogeneous paths 
consisting of not necessarily perfect crystals.
This is also a main subject of the paper.
We propose a curious answer compatible with 
the thermodynamic Bethe ansatz as well as 
the crystal theoretic scheme in \cite{HKKOT}.
Roughly, inhomogeneity affects the path realization into that for    
tensor products of several highest weight modules.
Non-perfectness leads to a description in terms of embeddings of 
affine Lie algebras.
See \cite{AM,DMN,HKMW} for a physical background 
of considering inhomogeneous paths.

Before going into the contents of the paper, 
we note further aspects of the subject including some latest ones; 
applications to quasi-hypergeometric functions \cite{AI}, 
exclusion statistics \cite{BM,BS2},
decomposition of the Kirillov-Reshetikhin modules $W^{(k)}_s$ 
(see Section \ref{subsec:family})
into classical irreducible modules \cite{C1},
the Feigin-Stoyanovsky theory \cite{FS},
fermionic formulae for the spaces of coinvariants associated to 
$\widehat{sl}_2$ \cite{FKLMM},
modular representations of Hecke algebras \cite{FLOTW},
approaches by various bijections and tableaux \cite{FOW,KSS,SW,S,SS1}, 
a graphical computational algorithm of fermionic formulae \cite{Kl},
dilogarithm identities \cite{Ki4},
fermionic formulae for weight multiplicities \cite{KN},
connection to geometry of quiver varieties \cite{Lu,N},
spinons at $q=0$ \cite{NY3,NY4},
fermionic formulae {}from fixed boundary RSOS models \cite{OPW},
decomposition of level 1 modules by level 0 actions \cite{Ta,KKN,BS1},
bosonic formulae {}from crystals \cite{SS2,KMOTU}, etc.

\subsection{Contents of the paper}\label{subsec:contents}

Let us explain the contents of the paper in detail 
along Sections \ref{sec:cfc}-\ref{sec:qsys} and Appendices 
\ref{app:mlist}-\ref{app:example}.

In Section \ref{sec:cfc}, we fix basic notations 
on affine Lie algebra $\geh=X^{(r)}_N$ following \cite{Kac}, and 
crystals of $U_q(\geh)$ (or $U'_q(\geh)$)\footnote{
Although denoted by the same letter, 
the deformation parameter $q$ is unrelated with 
the argument of fermionic formulae having the 
meaning $q= e^{-\delta}$ in terms of the null root $\delta$.} 
following \cite{KMN1,KMN2}.
We consider two categories of representations and crystals; one for the 
integrable highest weight modules and the other for 
finite-dimensional modules.
Crystals in the latter category are called finite crystals, for which 
the notion of perfect crystals is  explained.
We introduce irreducible 
finite-dimensional $U'_q(\geh)$-modules $W^{(k)}_s$ and associated 
conjectural family of finite crystal basis $B^{k,s}$.
In particular, all the $B^{k,s}$'s are conjectured to be perfect
in the twisted cases $r > 1$.

In Section \ref{sec:paths} we define a path
as an element of the semi-infinite tensor product
$\cdots \ot B \ot B$ of a finite crystal $B$ obeying 
a certain boundary condition.
Here  $B$ has the form 
$B=B^{k_1,s_1}\ot B^{k_2,s_2}\ot\cd\ot B^{k_d,s_d}$ and the resulting 
paths are called inhomogeneous in the sense that  $(k_i,s_i)$ can be 
$i$-dependent in general.
Homogeneous paths ($d=1$) consisting of perfect $B= B^{k_1,s_1}$
were originally studied in \cite{KMN1}.
We state Conjectures 
\ref{conj:tensor-prod-th} and \ref{conj:tensor-prod-th-np},
which identify the set of paths with the tensor product
of crystals of integrable highest weight modules.
The former conjecture 
concerns the situation where all $B^{k_i,s_i}$'s are perfect.
It can actually be proved for a large list of crystals 
as mentioned in Remark \ref{rem:One}.
The latter conjecture deals with non-perfect cases, and is described, 
rather curiously, in terms of the 
embedding of the affine Lie algebras 
$B^{(1)}_n \hookrightarrow D^{(1)}_{n+1}, 
C^{(1)}_n \hookrightarrow A^{(1)}_{2n-1}, 
F^{(1)}_4 \hookrightarrow E^{(1)}_6$ and 
$G^{(1)}_2 \hookrightarrow B^{(1)}_3 \hookrightarrow D^{(1)}_4$.
See (\ref{eq:incl}).
It originates in the thermodynamic Bethe ansatz 
calculation of RSOS central charges \cite{Ku} and is supported by 
computer experiments.
In Section \ref{subsec:sums} 
we introduce one dimensional sums (1dsums) 
$X_l(B,\la,q)$ in (\ref{eq:defX}), 
where $l \in \Z_{\ge 0}$ or $l = \infty$ 
is the level of restriction and $\la$ 
is the weight of the relevant paths.
They are crystal theoretic formulation of a trace of corner transfer matrices 
going back to \cite{ABF,B}.
The 1dsums over inhomogeneous paths firstly appeared in 
\cite{NY1} for $A^{(1)}_n$ case.
It makes essential use of the combinatorial $R$ matrices, 
i.e., the isomorphism and energy functions 
on tensor products of finite crystals.
Our definition basically follows theirs except  
a subtle ``boundary energy" (2nd term in (\ref{eq:subtle})) which 
becomes crucial in general if  $\geh \neq A^{(1)}_n$.
Our fundamental observation is Conjecture \ref{conj:X=M} 
equating the 1dsums on finite tensor products with fermionic formulae.
In Section \ref{subsec:representation} we explain how the 
finite tensor products are embedded into the crystals 
of integrable highest weight modules.

In Section \ref{sec:FF} we define fermionic formulae labeled by
the tensor product 
$W=W^{(k_1)}_{s_1} \ot \cdots \ot W^{(k_d)}_{s_d}$.
There are two versions denoted by 
$M_l(W,\la,q)$ (\ref{eq:mm}) and $\tilde{M}_l(W,\la,q)$ (\ref{eq:Mtilde}).
Compared to $M_l$, $\tilde{M}_l$ involves 
unphysical contributions {}from negative vacancy numbers 
in the Bethe ansatz. 
However as emphasized in \cite{HKOTY1}, 
they conjecturally coincide and 
the latter enjoys a skew invariance under the Weyl group action
(Conjecture \ref{conj:MM}).
We present the spinon character formulae 
derived in the infinite tensor product limit 
in Propositions \ref{pr:spinon} and \ref{pr:spinon2}.
Proposition \ref{pr:dualspinon} describes the other limit 
where  $q$ is replaced with $q^{-1}$.

In Section \ref{sec:N} we present a conjectural 
$q$-series formula $N_l(\la,q)$ 
for level $l$ vacuum string functions \cite{Kac}.
It extends the nontwisted cases \cite{KNS}, where
a Bethe ansatz interpretation has been proposed.
We point out the relation (\ref{eq:MN}) to the  
fermionic formula $M_l$, although its meaning is yet to be 
clarified.
So far the conjecture is known valid 
for $A^{(1)}_1$ \cite{LP} and
$A^{(1)}_n$ \cite{G,HKKOTY} as well as $l=1$ case \cite{Kac}.
The $q$-series exhibits an asymptotic behavior 
consistent with \cite{Kac}
under the dilogarithm sum rule (\ref{eq:dilog}).
In particular the sum rule for the twisted cases 
$X^{(r)}_N \; (r > 1)$ 
can be obtained {}from the nontwisted one 
$X^{(1)}_N$ by a simple unfolding procedure due to
(\ref{eq:hehe}).
We remark that a generalization of $N_l(\la,q)$ related to  
tensor product of several vacuum modules 
has been proved for $A^{(1)}_n$ and 
conjectured for the other nontwisted cases 
in the section 6.1 of \cite{HKKOTY}.

Section \ref{sec:qsys} is devoted to the $Q$-system.
It is a possible character identity among 
$Q^{(a)}_j = \mbox{ch } W^{(a)}_j$
having a form of 2d Toda equations and 
essentially governs the fermionic formulae at $q=1$.
It was proposed firstly in \cite{KR2,Ki3} for nontwisted cases.
We introduce the twisted cases $X^{(r)}_N$ by folding 
the $Q$-system of the underlying simply laced 
case $X^{(1)}_N$ by an order $r$ Dynkin diagram automorphism.
Explicit solutions are presented 
for $A^{(2)}_{2n}, A^{(2)}_{2n-1}, D^{(2)}_{n+1}$ and 
$D^{(3)}_4$ in terms of characters of various 
classical subalgebras.
The results given there, Appendix \ref{app:BD} and \cite{KR2,HKOTY1},
cover all the non-exceptional affine Lie algebras $\geh$, and 
exhaust the possible choices of their classical simple Lie subalgebras
obtained by removing a vertex in the  Dynkin diagrams.
We find that only the choice of the classical subalgebra 
$\overset{\circ}{\geh}$ (Table \ref{tab:geh-bar}) 
fits the requirements (A) and (C) 
of Theorem \ref{th:completeness}.
The identity like (\ref{eq:completeness}) therein  
is sometimes called ``combinatorial completeness 
of the string hypothesis in the Bethe ansatz".
For nontwisted cases $U_q(X^{(1)}_N)$, 
the reason is that 
$M_\infty(W,\la,1)$ naturally emerges 
{}from a formal counting of the number of solutions to the 
Bethe equation \cite{OW} for rational vertex models 
with Yangian symmetry $Y(X_N)$\footnote{
Although it is not $M_\infty$  but ${\tilde M}_\infty$ 
that enters (\ref{eq:completeness}), 
we admit $M_\infty = \tilde{M}_\infty$ here according to 
Conjecture \ref{conj:MM}.}.
However for twisted cases, we have {\em not} succeeded in deriving 
our $q=1$ fermionic formula 
{}from the relevant 
Bethe equation \cite{RW} under a string hypothesis.
Consequently we do not yet know an interpretation of 
the identity (\ref{eq:completeness}) as a 
combinatorial completeness for the rational limit of 
trigonometric $U'_q(X^{(r)}_N)$ vertex models.
See the end of Section \ref{sec:complete} for related remarks.
In order to find the $q=1$ fermionic formula with $r>1$, 
we set out by postulating the $Q$-system first 
instead of investigating the Bethe equation.
Once the $Q$-system is specified, 
one can execute parallel calculations with the section 8.2 of \cite{HKOTY1}.
Our $M_\infty(W,\la,1)$ has come up {}from such analyses 
simultaneously with Theorem \ref{th:completeness}.
To invent the $q$-analogue is then straightforward by 
presuming Conjecture \ref{conj:X=M}, i.e., 
$M_\infty(W,\la,q) = q^{-D} X_\infty(B,\la,q)$ for some $D$.

Appendix \ref{app:mlist} provides a list of 
the basic fermionic formulae
$M_\infty(W^{(a)}_s,\la,q^{-1})$ for twisted cases.
The results are analytically obtained for non-exceptional algebras.
(See  appendix A of  \cite{HKOTY1} for 
a similar list  for nontwisted cases.)
Appendix \ref{app:regII} gives the explicit form of 
Proposition \ref{pr:dualspinon}.
Appendix \ref{app:BD} contains the 
solution of $B^{(1)}_n$ $Q$-system in terms of 
$D_n$ characters. 
Appendix D  provides the lists  of
isomorphism $B^{1,s_1} \ot B^{1,s_2} \simeq B^{1,s_2} \ot B^{1,s_1}$
$(1 \le s_2 \le s_1 \le 2)$ and energy function for 
$\geh = B^{(1)}_3, C^{(1)}_2, D^{(1)}_4, A^{(2)}_3, A^{(2)}_4, D^{(2)}_3$ 
and $D^{(3)}_4$.
They form sufficient data to compute the inhomogeneous 1dsums 
over the tensor products of these crystals.

Let us close the introduction with  
two remarks.
Solvable RSOS models at a generic primitive root of unity 
give rise to yet further 1dsums associated with complicated 
energy functions. 
A typical example is \cite{FB} corresponding to $A^{(1)}_1$, and 
as noted in Section \ref{subsec:sb}, 
intricate arithmetic structures have been worked out on 
the fermionic formula and relation to 
fractional level representations \cite{BMS,BMSW,FW2}.
The crystal theory has not yet been adapted to such situations.
It may be interesting to seek an analogous 
generalization of our fermionic formula.

According to Section \ref{subsec:representation},
the limit of our 1dsums 
$\lim_{L \rightarrow \infty}q^{-D}X_l(B^{\ot L},\la,q)$ 
($D$: normalization)  
yields the $q$-series having an affine Lie algebraic meaning.
Such  phenomena correspond, though rather roughly, 
to the so called ``regime III" or
``antiferromagnetic regime"
in the language of solvable RSOS models or vertex models.
However, those models indicate that the other 
regimes are no less interesting.
In particular, the analysis of the 
opposite ``regime II" or ``ferromagnetic regime" 
renders the following question:
what is the representation theoretical meaning of another $q$-series 
$\lim_{L \rightarrow \infty}q^{-D'}X_l(B^{\ot L},\la,q^{-1})$ 
($D'$: another normalization)?
Compared with ``regime III", such issues in ``regime II" 
are yet poorly understood\footnote{
If $\geh = A^{(1)}_n$ and $l < \infty$,  the question can be answered 
indirectly through the level-rank duality.}.
See Proposition \ref{pr:dualspinon}, Appendix \ref{app:regII} 
and (\ref{eq:MN}).
For $l = \infty$ the problem is relevant to another topic, 
the soliton cellular automata associated with crystal bases
\cite{HKOTY2}.

\medskip

\section*{Acknowledgements}
The authors thank V. Chari, T. Nakanishi, A. Schilling, M. Shimozono 
and J. Suzuki for useful discussions and Y. Yamada 
for the collaboration in their previous work \cite{HKOTY1}. 
A.K. and M.O. are partially supported by Grant-in-Aid for 
Scientific Research {}from the Ministry of Education, Culture, Sports,
Science and Technology of Japan.
Z.T. is supported by JSPS Research Fellowships for Young Scientists.

\section{Conjectural family of crystals}\label{sec:cfc}

In this section we prepare necessary notations for the affine Lie algebra $\geh$, quantum 
affine algebra $U_q(\geh)$ and its crystals. We also present a conjectural family of 
crystals revealed by the Bethe ansatz.

\subsection{\mathversion{bold} Affine Lie algebra $X^{(r)}_N$}\label{subsec:affine}

Let $\geh$ be a Kac-Moody Lie algebra of affine type $X^{(r)}_N= 
A^{(1)}_n (n \ge 1), B^{(1)}_n (n \ge 3), C^{(1)}_n (n \ge 2),
D^{(1)}_n (n \ge 4), E^{(1)}_n (n=6,7,8), F^{(1)}_4, G^{(1)}_2,
A^{(2)}_{2n}(n\ge1),A^{(2)}_{2n-1}(n\ge2),D^{(2)}_{n+1}(n\ge2),
E^{(2)}_6$ and $D^{(3)}_4$.
As is well known, $X^{(r)}_N$ is realized as the canonical
central extension of the loop algebra based on the pair 
$(X_N, \sigma)$ of finite-dimensional simple Lie algebra 
$X_N$  and its Dynkin diagram automorphism $\sigma$ of order $r=1,2,3$.
The vertices of the Dynkin diagram of $\geh = X^{(r)}_N$ are labeled 
as in Table \ref{tab:Dynkin}, which is quoted {}from
TABLE Aff1-3 in \cite{Kac}.

{\unitlength=.95pt
\begin{table}
\caption{Dynkin diagrams for $X^{(r)}_N$.
The enumeration of the nodes with
$I = \{0,1,\ldots, n\}$ is specified under or the right side of the nodes.
In addition, the numbers $t_i$ (resp. $t^\vee_i$) defined in 
(\ref{eq:ttdef}) are attached {\em above} the nodes for $r=1$ (resp. $r>1$)
if and only if $t_i \neq 1$ (resp. $t^\vee_i \neq 1$).}
\label{tab:Dynkin}
\begin{tabular}[t]{rl}
$A_1^{(1)}$:&
\begin{picture}(26,20)(-5,-5)
\multiput( 0,0)(20,0){2}{\circle{6}}
\multiput(2.85,-1)(0,2){2}{\line(1,0){14.3}}
\put(0,-5){\makebox(0,0)[t]{$0$}}
\put(20,-5){\makebox(0,0)[t]{$1$}}
\put( 6, 0){\makebox(0,0){$<$}}
\put(14, 0){\makebox(0,0){$>$}}
\end{picture}
\\
&
\\
\begin{minipage}[b]{4em}
\begin{flushright}
$A_n^{(1)}$:\\$(n \ge 2)$
\end{flushright}
\end{minipage}&
\begin{picture}(106,40)(-5,-5)
\multiput( 0,0)(20,0){2}{\circle{6}}
\multiput(80,0)(20,0){2}{\circle{6}}
\put(50,20){\circle{6}}
\multiput( 3,0)(20,0){2}{\line(1,0){14}}
\multiput(63,0)(20,0){2}{\line(1,0){14}}
\multiput(39,0)(4,0){6}{\line(1,0){2}}
\put(2.78543,1.1142){\line(5,2){44.429}}
\put(52.78543,18.8858){\line(5,-2){44.429}}
\put(0,-5){\makebox(0,0)[t]{$1$}}
\put(20,-5){\makebox(0,0)[t]{$2$}}
\put(80,-5){\makebox(0,0)[t]{$n\!\! -\!\! 1$}}
\put(100,-5){\makebox(0,0)[t]{$n$}}
\put(55,20){\makebox(0,0)[lb]{$0$}}
\end{picture}
\\
&
\\
\begin{minipage}[b]{4em}
\begin{flushright}
$B_n^{(1)}$:\\$(n \ge 3)$
\end{flushright}
\end{minipage}&
\begin{picture}(126,40)(-5,-5)
\multiput( 0,0)(20,0){3}{\circle{6}}
\multiput(100,0)(20,0){2}{\circle{6}}
\put(20,20){\circle{6}}
\multiput( 3,0)(20,0){3}{\line(1,0){14}}
\multiput(83,0)(20,0){1}{\line(1,0){14}}
\put(20,3){\line(0,1){14}}
\multiput(102.85,-1)(0,2){2}{\line(1,0){14.3}} 
\multiput(59,0)(4,0){6}{\line(1,0){2}} 
\put(110,0){\makebox(0,0){$>$}}
\put(0,-5){\makebox(0,0)[t]{$1$}}
\put(20,-5){\makebox(0,0)[t]{$2$}}
\put(40,-5){\makebox(0,0)[t]{$3$}}
\put(100,-5){\makebox(0,0)[t]{$n\!\! -\!\! 1$}}
\put(120,-5){\makebox(0,0)[t]{$n$}}
\put(25,20){\makebox(0,0)[l]{$0$}}

\put(120,13){\makebox(0,0)[t]{$2$}}
\end{picture}
\\
&
\\
\begin{minipage}[b]{4em}
\begin{flushright}
$C_n^{(1)}$:\\$(n \ge 2)$
\end{flushright}
\end{minipage}&
\begin{picture}(126,20)(-5,-5)
\multiput( 0,0)(20,0){3}{\circle{6}}
\multiput(100,0)(20,0){2}{\circle{6}}
\multiput(23,0)(20,0){2}{\line(1,0){14}}
\put(83,0){\line(1,0){14}}
\multiput( 2.85,-1)(0,2){2}{\line(1,0){14.3}} 
\multiput(102.85,-1)(0,2){2}{\line(1,0){14.3}} 
\multiput(59,0)(4,0){6}{\line(1,0){2}} 
\put(10,0){\makebox(0,0){$>$}}
\put(110,0){\makebox(0,0){$<$}}
\put(0,-5){\makebox(0,0)[t]{$0$}}
\put(20,-5){\makebox(0,0)[t]{$1$}}
\put(40,-5){\makebox(0,0)[t]{$2$}}
\put(100,-5){\makebox(0,0)[t]{$n\!\! -\!\! 1$}}
\put(120,-5){\makebox(0,0)[t]{$n$}}

\put(20,13){\makebox(0,0)[t]{$2$}}
\put(40,13){\makebox(0,0)[t]{$2$}}
\put(100,13){\makebox(0,0)[t]{$2$}}
\end{picture}
\\
&
\\
\begin{minipage}[b]{4em}
\begin{flushright}
$D_n^{(1)}$:\\$(n \ge 4)$
\end{flushright}
\end{minipage}&
\begin{picture}(106,40)(-5,-5)
\multiput( 0,0)(20,0){2}{\circle{6}}
\multiput(80,0)(20,0){2}{\circle{6}}
\multiput(20,20)(60,0){2}{\circle{6}}
\multiput( 3,0)(20,0){2}{\line(1,0){14}}
\multiput(63,0)(20,0){2}{\line(1,0){14}}
\multiput(39,0)(4,0){6}{\line(1,0){2}}
\multiput(20,3)(60,0){2}{\line(0,1){14}}
\put(0,-5){\makebox(0,0)[t]{$1$}}
\put(20,-5){\makebox(0,0)[t]{$2$}}
\put(80,-5){\makebox(0,0)[t]{$n\!\! - \!\! 2$}}
\put(103,-5){\makebox(0,0)[t]{$n\!\! -\!\! 1$}}
\put(25,20){\makebox(0,0)[l]{$0$}}
\put(85,20){\makebox(0,0)[l]{$n$}}
\end{picture}
\\
&
\\
$E_6^{(1)}$:&
\begin{picture}(86,60)(-5,-5)
\multiput(0,0)(20,0){5}{\circle{6}}
\multiput(40,20)(0,20){2}{\circle{6}}
\multiput(3,0)(20,0){4}{\line(1,0){14}}
\multiput(40, 3)(0,20){2}{\line(0,1){14}}
\put( 0,-5){\makebox(0,0)[t]{$1$}}
\put(20,-5){\makebox(0,0)[t]{$2$}}
\put(40,-5){\makebox(0,0)[t]{$3$}}
\put(60,-5){\makebox(0,0)[t]{$4$}}
\put(80,-5){\makebox(0,0)[t]{$5$}}
\put(45,20){\makebox(0,0)[l]{$6$}}
\put(45,40){\makebox(0,0)[l]{$0$}}
\end{picture}
\\
&
\\
$E_7^{(1)}$:&
\begin{picture}(126,40)(-5,-5)
\multiput(0,0)(20,0){7}{\circle{6}}
\put(60,20){\circle{6}}
\multiput(3,0)(20,0){6}{\line(1,0){14}}
\put(60, 3){\line(0,1){14}}
\put( 0,-5){\makebox(0,0)[t]{$0$}}
\put(20,-5){\makebox(0,0)[t]{$1$}}
\put(40,-5){\makebox(0,0)[t]{$2$}}
\put(60,-5){\makebox(0,0)[t]{$3$}}
\put(80,-5){\makebox(0,0)[t]{$4$}}
\put(100,-5){\makebox(0,0)[t]{$5$}}
\put(120,-5){\makebox(0,0)[t]{$6$}}
\put(65,20){\makebox(0,0)[l]{$7$}}
\end{picture}
\\
&
\\
\end{tabular}
\begin{tabular}[t]{rl}
$E_8^{(1)}$:&
\begin{picture}(146,40)(-5,-5)
\multiput(0,0)(20,0){8}{\circle{6}}
\put(100,20){\circle{6}}
\multiput(3,0)(20,0){7}{\line(1,0){14}}
\put(100, 3){\line(0,1){14}}
\put( 0,-5){\makebox(0,0)[t]{$0$}}
\put(20,-5){\makebox(0,0)[t]{$1$}}
\put(40,-5){\makebox(0,0)[t]{$2$}}
\put(60,-5){\makebox(0,0)[t]{$3$}}
\put(80,-5){\makebox(0,0)[t]{$4$}}
\put(100,-5){\makebox(0,0)[t]{$5$}}
\put(120,-5){\makebox(0,0)[t]{$6$}}
\put(140,-5){\makebox(0,0)[t]{$7$}}
\put(105,20){\makebox(0,0)[l]{$8$}}
\end{picture}
\\
&
\\
$F_4^{(1)}$:&
\begin{picture}(86,20)(-5,-5)
\multiput( 0,0)(20,0){5}{\circle{6}}
\multiput( 3,0)(20,0){2}{\line(1,0){14}}
\multiput(42.85,-1)(0,2){2}{\line(1,0){14.3}} 
\put(63,0){\line(1,0){14}}
\put(50,0){\makebox(0,0){$>$}}
\put( 0,-5){\makebox(0,0)[t]{$0$}}
\put(20,-5){\makebox(0,0)[t]{$1$}}
\put(40,-5){\makebox(0,0)[t]{$2$}}
\put(60,-5){\makebox(0,0)[t]{$3$}}
\put(80,-5){\makebox(0,0)[t]{$4$}}
\put(60,13){\makebox(0,0)[t]{$2$}}
\put(80,13){\makebox(0,0)[t]{$2$}}
\end{picture}
\\
&
\\
$G_2^{(1)}$:&
\begin{picture}(46,20)(-5,-5)
\multiput( 0,0)(20,0){3}{\circle{6}}
\multiput( 3,0)(20,0){2}{\line(1,0){14}}
\multiput(22.68,-1.5)(0,3){2}{\line(1,0){14.68}}
\put( 0,-5){\makebox(0,0)[t]{$0$}}
\put(20,-5){\makebox(0,0)[t]{$1$}}
\put(40,-5){\makebox(0,0)[t]{$2$}}
\put(30,0){\makebox(0,0){$>$}}
\put(40,13){\makebox(0,0)[t]{$3$}}
\end{picture}
\\
&
\\
$A^{(2)}_2$:&
\begin{picture}(26,20)(-5,-5)
\multiput( 0,0)(20,0){2}{\circle{6}}
\multiput(2.958,-0.5)(0,1){2}{\line(1,0){14.084}}
\multiput(2.598,-1.5)(0,3){2}{\line(1,0){14.804}}
\put(0,-5){\makebox(0,0)[t]{$0$}}
\put(20,-5){\makebox(0,0)[t]{$1$}}
\put(10,0){\makebox(0,0){$<$}}
\put(0,13){\makebox(0,0)[t]{$2$}}
\put(20,13){\makebox(0,0)[t]{$2$}}
\end{picture}
\\
&
\\
\begin{minipage}[b]{4em}
\begin{flushright}
$A_{2n}^{(2)}$:\\$(n \ge 2)$
\end{flushright}
\end{minipage}&
\begin{picture}(126,20)(-5,-5)
\multiput( 0,0)(20,0){3}{\circle{6}}
\multiput(100,0)(20,0){2}{\circle{6}}
\multiput(23,0)(20,0){2}{\line(1,0){14}}
\put(83,0){\line(1,0){14}}
\multiput( 2.85,-1)(0,2){2}{\line(1,0){14.3}} 
\multiput(102.85,-1)(0,2){2}{\line(1,0){14.3}} 
\multiput(59,0)(4,0){6}{\line(1,0){2}} 
\put(10,0){\makebox(0,0){$<$}}
\put(110,0){\makebox(0,0){$<$}}
\put(0,-5){\makebox(0,0)[t]{$0$}}
\put(20,-5){\makebox(0,0)[t]{$1$}}
\put(40,-5){\makebox(0,0)[t]{$2$}}
\put(100,-5){\makebox(0,0)[t]{$n\!\! -\!\! 1$}}
\put(120,-5){\makebox(0,0)[t]{$n$}}
\put(0,13){\makebox(0,0)[t]{$2$}}
\put(20,13){\makebox(0,0)[t]{$2$}}
\put(40,13){\makebox(0,0)[t]{$2$}}
\put(100,13){\makebox(0,0)[t]{$2$}}
\put(120,13){\makebox(0,0)[t]{$2$}}
\put(120,13){\makebox(0,0)[t]{$2$}}
\end{picture}
\\
&
\\
\begin{minipage}[b]{4em}
\begin{flushright}
$A_{2n-1}^{(2)}$:\\$(n \ge 3)$
\end{flushright}
\end{minipage}&
\begin{picture}(126,40)(-5,-5)
\multiput( 0,0)(20,0){3}{\circle{6}}
\multiput(100,0)(20,0){2}{\circle{6}}
\put(20,20){\circle{6}}
\multiput( 3,0)(20,0){3}{\line(1,0){14}}
\multiput(83,0)(20,0){1}{\line(1,0){14}}
\put(20,3){\line(0,1){14}}
\multiput(102.85,-1)(0,2){2}{\line(1,0){14.3}} 
\multiput(59,0)(4,0){6}{\line(1,0){2}} 
\put(110,0){\makebox(0,0){$<$}}
\put(0,-5){\makebox(0,0)[t]{$1$}}
\put(20,-5){\makebox(0,0)[t]{$2$}}
\put(40,-5){\makebox(0,0)[t]{$3$}}
\put(100,-5){\makebox(0,0)[t]{$n\!\! -\!\! 1$}}
\put(120,-5){\makebox(0,0)[t]{$n$}}
\put(25,20){\makebox(0,0)[l]{$0$}}

\put(120,13){\makebox(0,0)[t]{$2$}}
\end{picture}
\\
&
\\
\begin{minipage}[b]{4em}
\begin{flushright}
$D_{n+1}^{(2)}$:\\$(n \ge 2)$
\end{flushright}
\end{minipage}&
\begin{picture}(126,20)(-5,-5)
\multiput( 0,0)(20,0){3}{\circle{6}}
\multiput(100,0)(20,0){2}{\circle{6}}
\multiput(23,0)(20,0){2}{\line(1,0){14}}
\put(83,0){\line(1,0){14}}
\multiput( 2.85,-1)(0,2){2}{\line(1,0){14.3}} 
\multiput(102.85,-1)(0,2){2}{\line(1,0){14.3}} 
\multiput(59,0)(4,0){6}{\line(1,0){2}} 
\put(10,0){\makebox(0,0){$<$}}
\put(110,0){\makebox(0,0){$>$}}
\put(0,-5){\makebox(0,0)[t]{$0$}}
\put(20,-5){\makebox(0,0)[t]{$1$}}
\put(40,-5){\makebox(0,0)[t]{$2$}}
\put(100,-5){\makebox(0,0)[t]{$n\!\! -\!\! 1$}}
\put(120,-5){\makebox(0,0)[t]{$n$}}

\put(20,13){\makebox(0,0)[t]{$2$}}
\put(40,13){\makebox(0,0)[t]{$2$}}
\put(100,13){\makebox(0,0)[t]{$2$}}
\end{picture}
\\
&
\\
$E_6^{(2)}$:&
\begin{picture}(86,20)(-5,-5)
\multiput( 0,0)(20,0){5}{\circle{6}}
\multiput( 3,0)(20,0){2}{\line(1,0){14}}
\multiput(42.85,-1)(0,2){2}{\line(1,0){14.3}} 
\put(63,0){\line(1,0){14}}
\put(50,0){\makebox(0,0){$<$}}
\put( 0,-5){\makebox(0,0)[t]{$0$}}
\put(20,-5){\makebox(0,0)[t]{$1$}}
\put(40,-5){\makebox(0,0)[t]{$2$}}
\put(60,-5){\makebox(0,0)[t]{$3$}}
\put(80,-5){\makebox(0,0)[t]{$4$}}

\put(60,13){\makebox(0,0)[t]{$2$}}
\put(80,13){\makebox(0,0)[t]{$2$}}
\end{picture}
\\
&
\\
$D_4^{(3)}$:&
\begin{picture}(46,20)(-5,-5)
\multiput( 0,0)(20,0){3}{\circle{6}}
\multiput( 3,0)(20,0){2}{\line(1,0){14}}
\multiput(22.68,-1.5)(0,3){2}{\line(1,0){14.68}}
\put( 0,-5){\makebox(0,0)[t]{$0$}}
\put(20,-5){\makebox(0,0)[t]{$1$}}
\put(40,-5){\makebox(0,0)[t]{$2$}}
\put(30,0){\makebox(0,0){$<$}}

\put(40,13){\makebox(0,0)[t]{$3$}}
\end{picture}
\\
&
\\
\end{tabular}
\end{table}}
%
\begin{table}[h]
\caption{}\label{tab:geh-bar}
\begin{center}
\begin{tabular}{c|cccccc}
$\geh$ & $X^{(1)}_N$ & $A^{(2)}_{2n}$ & $A^{(2)}_{2n-1}$ & $D^{(2)}_{n+1}$ 
& $E^{(2)}_6$ & $D^{(3)}_4$ \\
\hline
$\overset{\circ}{\geh}$ & $X_N$ & $C_n$ & $C_n$ & $B_n$ & $F_4$ & $G_2$ \\
$\geh_{\ol{0}}$ &         $X_N$ & $B_n$ & $C_n$ & $B_n$ & $F_4$ & $G_2$
\end{tabular}
\end{center}
\end{table}

Among  finite-dimensional subalgebras of $\geh$, the following
$\overset{\circ}{\geh}$ and $\geh_{\ol{0}}$ will play an essential role.
Here $\overset{\circ}{\geh}$ is obtained by removing the $0$-th vertex
{}from the Dynkin diagram of $\geh$. 
See the section 6.3 of \cite{Kac}.
On the other hand $\geh_{\ol{0}}$ is the subalgebra of $X_N$ fixed by 
$\sigma$. 
See the section 8.3 of \cite{Kac}.
Note that $\overset{\circ}{\geh} = \geh_{\ol{0}}$ except 
$\geh = A^{(2)}_{2n}$.

Table \ref{tab:Dynkin} fixes the enumeration of the simple roots of 
$\overset{\circ}{\geh}$ by the restriction to $\{1, \ldots, n\}$.
As for $\geh_{\ol{0}}$ we naturally employ the same enumeration 
with $\overset{\circ}{\geh}$ when $\geh \neq A^{(2)}_{2n}$, 
since they are the same subalgebra.
As for the conventions on $\geh_{\ol{0}}$ for $\geh=A^{(2)}_{2n}$, 
see Section \ref{subsec:convention}.
Let $\alpha_i,h_i,
\La_i$ ($i \in I=\{0,1, \ldots, n\}$) be the simple roots, simple coroots,
fundamental weights of $\geh$. 
Let $\delta$ and $c$ denote the generator of imaginary roots and the 
canonical central element, respectively. 
Recall that 
$\delta=\sum_{i\in I}a_i\alpha_i, c = \sum_{i \in I}a^\vee_i h_i$ and 
the Kac label $a_i$ is given explicitly by 
\[
\begin{array}{rll}
\delta&=\alpha_0+\cdots + \alpha_n
&\mbox{for }\geh=A^{(1)}_{n},\\
&=\alpha_0+\alpha_1+2\alpha_2+\cd+2\alpha_n
&\mbox{for }\geh=B^{(1)}_{n},\\
&=\alpha_0+2\alpha_1+\cd+2\alpha_{n-1}+\alpha_n
&\mbox{for }\geh=C^{(1)}_{n},\\
&=\alpha_0+\alpha_1+2\alpha_2+\cd+2\alpha_{n-2}
+\alpha_{n-1}+\alpha_n
&\mbox{for }\geh=D^{(1)}_{n},\\
&=\alpha_0+\alpha_1+2\alpha_2+3\alpha_3+2\alpha_4+\alpha_5+2\alpha_6
&\mbox{for }\geh=E^{(1)}_{6},\\
&=\alpha_0+2\alpha_1+3\alpha_2+4\alpha_3+3\alpha_4+2\alpha_5+
\alpha_6+2\alpha_7
&\mbox{for }\geh=E^{(1)}_{7},\\
&=\alpha_0+2\alpha_1+3\alpha_2+4\alpha_3+5\alpha_4+6\alpha_5+
4\alpha_6+2\alpha_7+3\alpha_8
&\mbox{for }\geh=E^{(1)}_{8},\\
&=\alpha_0+2\alpha_1+3\alpha_2+4\alpha_3+2\alpha_4
&\mbox{for }\geh=F^{(1)}_{4},\\
&=\alpha_0+2\alpha_1+3\alpha_2
&\mbox{for }\geh=G^{(1)}_{2},\\
&=2\alpha_0+2\alpha_1+\cd+2\alpha_{n-1}+\alpha_n
&\mbox{for }\geh=A^{(2)}_{2n},\\
&=\alpha_0+\alpha_1+2\alpha_2+\cd+2\alpha_{n-1}+\alpha_n\qquad
&\mbox{for }\geh=A^{(2)}_{2n-1},\\
&=\alpha_0+\alpha_1+\cd+\alpha_{n-1}+\alpha_n
&\mbox{for }\geh=D^{(2)}_{n+1},\\
&=\alpha_0+2\alpha_1+3\alpha_2+2\alpha_3+\alpha_4
&\mbox{for }\geh=E^{(2)}_{6},\\
&=\alpha_0+2\alpha_1+\alpha_2
&\mbox{for }\geh=D^{(3)}_{4}.
\end{array}
\]
The dual Kac label $a^\vee_i$ is equal to $a_i$ of  
$\geh^\vee$ corresponding to the transposed Cartan matrix.
We let $h^\vee = \sum_{i \in I}a^\vee_i$ denote the dual Coxeter number.
Note that $h^\vee$ for $X^{(r)}_N$ is independent of $r$ 
(Remark 6.1 in \cite{Kac}).

For $i \in I$  we set
\begin{equation}\label{eq:ttdef}
t_i = \max(\frac{a_i}{a^\vee_i},a^\vee_0),
\qquad
t^\vee_i = \max(\frac{a^\vee_i}{a_i},a_0).
\end{equation}
They are given in Table \ref{tab:Dynkin}.
By the definition $t_i, t^\vee_i \in \{1,2,3\}$, and 
$(t_i,t^\vee_i)$ for $\geh$ is equal to $(t^\vee_i,t_i)$ for $\geh^\vee$.
In the sequel, $t^\vee_i$ and $t_i$ are used only for 
$i \in I \setminus \{0\}$.
For $1 \le a \le n$ we have
\begin{equation*}
t^\vee_a = 1 \quad \mbox{if } r = 1,\qquad
t_a = 1 \quad \mbox{if } r > 1,
\end{equation*}
and especially for $\geh = A^{(2)}_{2n}$, 
$(t_a,t^\vee_a) = (1,2)$ for all $1 \le a \le n$.
We remark that 
\begin{equation}
ra_0\sum_{a=1}^n \frac{1}{t^\vee_a} = N
\quad \mbox{for any } X^{(r)}_N. \label{eq:wa}
\end{equation}

Let
$P=\bigoplus_{i \in I}{\mathbb Z}\La_i\bigoplus{\mathbb Z}\delta$ 
be the weight lattice.
We define the following subsets of $P$: $P^+=\sum_{i \in I}
{\mathbb Z}_{\ge 0}\La_i$, $P^+_l=\{\la\in P^+\mid \langle\la,c\rangle=l\}$,
$\ol{P}=\sum_{i = 1}^n{\mathbb Z}\ol{\La}_i$, $\ol{P}^+=\sum_{i=1}^{n}
{\mathbb Z}_{\ge 0}\ol{\La}_i$. 
Here $\ol{\La}_i=\La_i-\langle \La_i,c \rangle 
\La_0$ is the classical
part of $\La_i$. This map $\ol{\phantom{\La}}$ is extended to
a map on $P$ so that it is ${\mathbb Z}$-linear and 
$\ol{\delta} = 0$. 
We also consider the classical weight lattice 
$P_{cl}=P/{\mathbb Z}\delta$. 
In this paper we canonically identify $P_{cl}$ with 
$\bigoplus_{i\in I}{\mathbb Z}\La_i \subset P$.
(See the section 3.1 of \cite{KMN1} for a precise treatment.)
We let $\ol{Q} = \bigoplus_{i=1}^n {\mathbb Z} \alpha_i$ denote
the classical root lattice. 
Let $(\cdot|\cdot)$ be the standard bilinear
form on $P$ normalized by $(\delta|\la)=\langle c,\la\rangle$ 
for any $\la\in P$.
This normalization agrees with that in \cite{Kac}.

\subsection{Quantum affine algebra and crystals}\label{subsec:quantum}

Let $U_q(\geh)$ be the quantum affine algebra associated to $\geh$. 
We follow the section 2.1 of
\cite{KMN1} for its definitions. $U_q(\overset{\circ}{\geh})$ turns out a subalgebra of $U_q(\geh)$
in an obvious manner. 
$U'_q(\geh)$ is the one without the generator $q^d$, where $d$ is the degree
operator in $\geh$.

We next review terminology for crystals. 
A crystal basis $B$ of a $U_q(\geh)$ ($U'_q(\geh),U_q(\overset{\circ}{\geh})$)-module 
can be regarded as a set of basis vectors of the module at $q=0$. 
On $B$ one has the action of modified Chevalley generators 
(Kashiwara operators)
\[
{\tilde e}_i,{\tilde f}_i:\;
B\longrightarrow B\sqcup\{0\}.
\]
For $b\in B$ we set
$\veps_i(b)=\max\{k\ge0\mid{\tilde e}_i^kb\ne0\},
\vphi_i(b)=\max\{k\ge0\mid{\tilde f}_i^kb\ne0\}$.
If $B_1$ and $B_2$ are crystals, the crystal structure on the tensor product
$B_1\ot B_2$ is given by
\begin{eqnarray}
{\tilde e}_i(b_1\ot b_2)&=&
\begin{cases}
{\tilde e}_ib_1\ot b_2&\mbox{if }\vphi_i(b_1)\ge\veps_i(b_2),\\
b_1\ot {\tilde e}_ib_2&\mbox{otherwise},
\end{cases}\label{eq:tensor_e}\\
{\tilde f}_i(b_1\ot b_2)&=&
\begin{cases}
{\tilde f}_ib_1\ot b_2&\mbox{if }\vphi_i(b_1)>\veps_i(b_2),\\
b_1\ot {\tilde f}_ib_2&\mbox{otherwise}.
\end{cases}\label{eq:tensor_f}
\end{eqnarray}

We mainly use two categories of crystals. The first one contains the crystal
basis $B(\la)$ of the irreducible integrable $U_q(\geh)$-module 
$L(\la)$ with highest weight
$\la\in P^+$. $B(\la)$ is a $P$-weighted crystal.
The other one contains a crystal basis $B$ of a finite-dimensional
$U'_q(\geh)$-module. As opposed to $B(\la)$, $B$ is a finite set and 
$P_{cl}$-weighted. We shall call it a finite crystal. 
(See \cite{HKKOT} for the precise definition.)
For a finite crystal $B$, we set
$\veps(b)=\sum_i\veps_i(b)\La_i$, $\vphi(b)=\sum_i\vphi_i(b)\La_i$,
$\wt b = \vphi(b) - \veps(b)$, 
and introduce the level of $B$ by
\[
\lev B=\min\{\langle c,\veps(b)\rangle\mid b\in B\}.
\]
We further set $B_{\min}=\{b\in B\mid\langle c,\veps(b)\rangle=\lev B\}$,
and call an element of $B_{\min}$ minimal.
Let $B_1,B_2$ be finite crystals of level $l_1,l_2$ . Then $B_1\ot B_2$ turns out a finite crystal
of level $\max(l_1,l_2)$.
A $U'_q(\geh)$-module can be viewed as a $U_q(\overset{\circ}{\geh})$-module.
For the latter the irreducible representation with 
highest weight $\la \in \ol{P}^+$ will be denoted by $V(\la)$.

We also need the notion of perfect crystals \cite{KMN1,KMN2}. 
A perfect crystal is a finite crystal with an additional condition.
Let $B$ be a finite crystal of level $l$. $B$ is said to be perfect if
the maps $\veps,\vphi:B_{\min}\longrightarrow P^+_l$ are bijections.
Set $\sigma=\veps\circ\vphi^{-1}$.
$\sigma$\footnote{Distinguish this with $\sigma$ in $(X_N,\sigma)$ which appeared in 
Section \ref{subsec:affine}.} 
is related to a Dynkin diagram automorphism and can be extended on 
$P^+$.

\subsection{\mathversion{bold} Family $W^{(k)}_s$}\label{subsec:family}

It is known \cite{CP1,CP2} that the irreducible finite-dimensional 
$U'_q(\geh)$-modules are characterized by $n$-tuple of polynomials 
$\{P_a(u)\}_{a=1}^n$ called Drinfel'd polynomials. 
Set ${\hat q}_k=q^{(\alpha_k|\alpha_k)/2}$.
We define an irreducible finite-dimensional $U'_q(\geh)$-module 
$W^{(k)}_s$ ($1\le k\le n,s\ge1$) by its Drinfel'd polynomials\footnote{
In the twisted cases, we are afraid that our normalization might not
agree with that in \cite{CP2}.}:
\[
P_a(u)=\left\{
\begin{array}{ll}
(1-{\hat q}_k^{s-1}u)(1-{\hat q}_k^{s-3}u)
\cdots (1-{\hat q}_k^{-s+1}u)\quad&(a=k),\\
1&(a\neq k).
\end{array}\right.
\] 
It is called a Kirillov-Reshetikhin module.

\begin{conjecture} \label{conj:fd-module}
The finite-dimensional $U'_q(\geh)$-module $W^{(k)}_s$ has the following features:
\begin{itemize}
\item[(1)] $W^{(k)}_s$ has a crystal basis $B^{k,s}$. $B^{k,s}$ is a 
           finite crystal of level $\lceil \frac{s}{t_k} \rceil$\footnote{
                   The symbol $\lceil x \rceil$ denotes the smallest integer not less than $x$.}.
                   Moreover, it is perfect if $\frac{s}{t_k}$ is an integer, and not perfect if not an integer.
\item[(2)] As a $U_q(\overset{\circ}{\geh})$-module, $W^{(k)}_s$ decomposes itself into
\[
W^{(k)}_s=\bigoplus_{\la\in \ol{P}^+}M_\infty(W^{(k)}_s,\la,q=1)\,V(\la),
\]
where $M_l(W^{(k)}_s,\la,q)$ is defined in (\ref{eq:mm}).

\item[(3)] Set $Q^{(k)}_s=\ch W^{(k)}_s$, then $Q^{(k)}_s$ satisfies
           the $Q$-system (\ref{eq:qsys}).

\end{itemize}
\end{conjecture}

{}From (\ref{eq:mm}) the decomposition (2) of $W^{(k)}_s$ has the form
\[
W^{(k)}_s=V(s\ol{\La}_k)\oplus\cdots,
\]
where $\cdots$ contains irreducible modules with highest weights
strictly lower than $s\ol{\La}_k$ only. 
A list of such decompositions is available in Appendix A (by setting
$q=1$). 

Several remarks are in order.

\begin{remark}
Although they did not specify the Drinfel'd polynomials, 
parts (2) and (3) (for nontwisted cases) are essentially 
the  Kirillov-Reshetikhin conjecture \cite{KR2}.
In the nontwisted cases, part (2) of the conjecture was solved 
recently \cite{C1} for all cases if $\geh$ is classical and 
several cases if exceptional.
See also \cite{C2}.
\end{remark}

\begin{remark}
In \cite{Ka3} an irreducible finite-dimensional $\Uqp$-module $W(\varpi_i)$
(fundamental representation) having a crystal basis is constructed. 
It should agree with our $W^{(i)}_1$. 
\end{remark}

\begin{remark}
It is noted that $t_k>1$ occurs only when $\geh=B^{(1)}_n(k=n), C^{(1)}_n(k\neq n),
F^{(1)}_4(k=3,4), G^{(1)}_2(k=2)$. In particular, $t_k=1$ for any $k$, if $\geh$ is simply laced 
or of twisted type. In these cases any $B^{k,s}$ is conjectured to be perfect.
\end{remark}

\begin{remark} \label{rem:coherent}
There is a notion of coherent family of perfect crystals \cite{KKM}. We also conjecture
that $\{B^{k,t_kl}\}_{l\ge1}$ forms a coherent family of perfect crystals for any $k$.
\end{remark}

\begin{remark}
For twisted cases the following crystals are known explicitly.
\[
\begin{array}{ll}
\mbox{\cite{KMN2}} &A^{(2)}_{2n}:B^{1,s},\,A^{(2)}_{2n-1}:B^{1,s},\,
                    D^{(2)}_{n+1}:B^{1,s}B^{n,s},\\
\mbox{\cite{JnMO}} &A^{(2)}_{2n}:B^{k,1},\,A^{(2)}_{2n-1}:B^{k,1}.
\end{array}
\]
See also \cite{KKM} for the explicit crystal structure of $B^{1,s}$.
In \cite{JnMO} one should notice that the labeling of the Dynkin diagrams is
opposite.
See also \cite{HKOTY1} Remark 2.3 for nontwisted cases.
\end{remark}

\section{Paths and one dimensional sums}\label{sec:paths}

Here we define a set of paths $\P(\bp,B)$ and the one dimensional sum $X_l(B,\la,q)$.
We then present several conjectures related to them.

\subsection{Energy function and Yang-Baxter equation}\label{sec:energy}

We review the energy function and Yang-Baxter equation for finite crystals.
Let $B$ be a finite crystal. Introduce an indeterminate $z$ and set
\[
\Aff(B)=\{z^\gamma b\mid\gamma\in\Z,b\in B\}.
\]
It is called the affinization of $B$. 
On $\Aff(B)$ one can define the action of $\et{i}$ and $\ft{i}$ by 
\[
\et{i}(z^\gamma b)=z^{\gamma+\delta_{i0}}(\et{i}b),\quad
\ft{i}(z^\gamma b)=z^{\gamma-\delta_{i0}}(\ft{i}b).
\]
Let $B_1$ and $B_2$ be finite crystals. Then there exists an isomorphism 
of crystals called the combinatorial $R$ \cite{KMN1}:
\[
\begin{array}{rcl}
R\;:\quad\Aff(B_1)\ot\Aff(B_2)&\longrightarrow&\Aff(B_2)\ot\Aff(B_1)\\
z^{\gamma_1}b_1\ot z^{\gamma_2}b_2&\longmapsto&
z^{\gamma_2+H(b_1\ot b_2)}b'_2\ot z^{\gamma_1-H(b_1\ot b_2)}b'_1.
\end{array}
\]
$H(b_1\ot b_2)$ is called the energy function and uniquely determined up to
an additive constant by the following rule:
\begin{eqnarray*}
H(\tilde{e}_i(b_1\ot b_2))&=H(b_1\ot b_2)+1
&\mbox{ if }i=0,\vphi_0(b_1)\geq\veps_0(b_2),\vphi_0(b'_2)\geq\veps_0(b'_1),
\\
&=H(b_1\ot b_2)-1
&\mbox{ if }i=0,\vphi_0(b_1)<\veps_0(b_2),\vphi_0(b'_2)<\veps_0(b'_1),
\\
&\hspace{-6mm}=H(b_1\ot b_2)&\mbox{ otherwise}.
\end{eqnarray*}
If we want to emphasize $B_1\ot B_2$, we write $H_{B_1B_2}$ for $H$.
Since $R^2=\id$, we have $H_{B_1B_2}(b_1\ot b_2)=H_{B_2B_1}(b'_2\ot b'_1)$.
Note that if $B_1=B_2$, one has $b'_2=b_1,b'_1=b_2$.

Now let $B_1,B_2$ and $B_3$ be finite crystals. Then the Yang-Baxter equation
holds on $\Aff(B_1)\ot\Aff(B_2)\ot\Aff(B_3)$.
\begin{equation} \label{eq:YBE}
(\id\ot R)(R\ot\id)(\id\ot R)=(R\ot\id)(\id\ot R)(R\ot\id).
\end{equation}
For $b_i\in B_i(i=1,2,3)$ let $b'_1,b'_2,{\hat b}_3,{\check b}_3$ be defined by
\[
\begin{array}{rcl}
B_1\ot B_2&\simarrow&B_2\ot B_1\\
b_1\ot b_2&\mapsto&b'_2\ot b'_1
\end{array}
\begin{array}{rcl}
B_2\ot B_3&\simarrow&B_3\ot B_2\\
b_2\ot b_3&\mapsto&{\hat b}_3\ot b''_2
\end{array}
\begin{array}{rcl}
B_1\ot B_3&\simarrow&B_3\ot B_1\\
b'_1\ot b_3&\mapsto&{\check b}_3\ot b''_1
\end{array}
\]
under the isomorphism of crystals. One can check (\ref{eq:YBE}) implies 
\begin{equation} \label{eq:used-later}
H(b_2\ot b_3)+H(b_1\ot{\hat b}_3)=H(b'_1\ot b_3)+H(b'_2\ot{\check b}_3).
\end{equation}

Let $\tau$ be a Dynkin diagram automorphism, {\em i.e.}, a permutation on letters {}from 
$\{0,1,\cd,n\}$\footnote{
This $\tau$ should be distinguished {}from that in Section 
\ref{subsec:qsys}.}. 
By setting $\tau(\sum_i m_i\La_i)=\sum_i m_i\La_{\tau(i)}$ one can
also define an automorphism on $P^+$. Next let $B$ 
be a perfect crystal. Since $\tau$ is a Dynkin
diagram automorphism, one necessarily has a unique bijection, 
denoted again by $\tau$, {}from
$B$ onto itself satisfying 
\begin{equation} \label{eq:auto-B}
\et{\tau(i)}(\tau(b))=\tau(\et{i}b),\quad \ft{\tau(i)}(\tau(b))=\tau(\ft{i}b).
\end{equation}
Here one should understand $\tau(0)=0$. Note that (\ref{eq:auto-B}) implies 
$\veps(\tau(b))=\tau(\veps(b)),\vphi(\tau(b))=\tau(\vphi(b))$.

We end this subsection with

\begin{conjecture} \label{conj:combR}
Let $B_1,B_2$ be perfect crystals of level $l_1,l_2$ ($l_1\ge l_2$). 
Under the isomorphism
of crystals $B_1\ot B_2\simarrow B_2\ot B_1$, 
any element $b_1\ot b_2$ of $(B_1\ot B_2)_{\min}$
is mapped to $\sigma_1(b_2)\ot \vphi_{B_1}^{-1}(\vphi_{B_1}(b_1)+\wt b_2)$, 
where $\sigma_1$ is the 
Dynkin diagram automorphism associated to the perfect crystal 
$B_1$ (see the end of Section \ref{subsec:quantum}) and $\vphi_{B_1}$ is the
map $\vphi$ {}from $B_1$.
\end{conjecture}

\begin{remark} \label{rem:injective}
If each $B_i$ belongs to a coherent family of perfect crystals (This coherent
family may depend on $i$. See also Remark \ref{rem:coherent}.) and 
the map $\veps\times\vphi:B_i\rightarrow
(P^+)^2;b\mapsto(\veps(b),\vphi(b))$ is injective for any $i$, then the 
conjecture can be shown to be true using the theory of crystals
with core \cite{KK}. Such case covers all coherent families of perfect crystals listed in
\cite{KKM}. The conjecture is proven for any case $B_i=B^{k_i,s_i}$ for
type $A^{(1)}_n$ (\cite{SS1} Theorem 7.3), although the injectivity
of $\veps\times\vphi$ is no longer valid.
\end{remark}

\subsection{Set of paths} \label{subsec:paths}

Let $B$ be a finite crystal of level $l$. {}From $B$ we construct a 
subset of $\cd\ot B\ot\cd\ot B$ called a set of paths. 
First we fix a reference path 
$\bp=\cd\ot\bb_j\ot\cd\ot\bb_2\ot\bb_1$. For any $j$, $\veps(\bb_j)$ should
have level $l$, and satisfy 
\[
\vphi(\bb_{j+1})=\veps(\bb_j).
\]
Set 
\[
\P(\bp,B)=\{p=\cd\ot b_j\ot\cd\ot b_2\ot b_1\mid 
b_j\in B,b_J=\bb_J\mbox{ for }J\gg1\}.
\]
An element of $\P(\bp,B)$ is called a path. For a path $p\in\P(\bp,B)$ 
we can define the energy $E(p)$ and weight $\wt p$ by 
\begin{eqnarray*}
E(p)&=&\sum_{j=1}^\infty j(H(b_{j+1}\ot b_j)-H(\bb_{j+1}\ot\bb_j)),\\
\wt p&=&\vphi(\bb_1)+\sum_{j=1}^\infty(\wt b_j-\wt\bb_j)-
(E(p)/a_0)\delta,
\end{eqnarray*}
where $a_0$ is the $0$-th Kac label.

Suppose $B$ is perfect of level $l$. 
Take $\la\in P^+_l$ and define $\bb_j$ to be an element of $B_{\min}$ such
that $\veps(\bb_j)=\sigma^j\la$. We denote the corresponding reference path
by $\bp^{(\la)}$. It is known in \cite{KMN1} that we have 
the following isomorphism of $P$-weighted crystals.
\[
\P(\bp^{(\la)},B)\simeq B(\la).
\]

Next let $B_i$ be a perfect crystal of level $l_i$ for $i=1,\cd,d$. One may
assume $B_i=B^{k_i,t_{k_i}l_i}$ in Conjecture
\ref{conj:fd-module} with some
$k_i$ ($1\le k_i\le n$). We also assume
$l_1\ge l_2\ge\cd\ge l_d$. Set $B=B_1\ot B_2\ot\cd\ot B_d$. For $\la_i\in
P^+_{l_i-l_{i+1}}$ ($i=1,\cd,d;l_{d+1}=0$) define a reference path 
$\bp^{(\la_1,\cd,\la_d)}=\cd\ot\bb_j\ot\cd\ot\bb_2\ot\bb_1$ ($\bb_j\in B$) by
\begin{eqnarray}
\bb_j&=&\veps_{B_1}^{-1}(\sigma_1^j\la_1+\sigma_1^j\sigma_2^j\la_2+\cd+
\sigma_1^j\sigma_2^j\cd\sigma_d^j\la_d) \label{eq:def-boldb}\\
&&\quad\ot\veps_{B_2}^{-1}(\sigma_1^{j-1}\sigma_2^j\la_2+
\sigma_1^{j-1}\sigma_2^j\sigma_3^j\la_3+\cd+\sigma_1^{j-1}\sigma_2^j\cd
\sigma_d^j\la_d) \nonumber\\
&&\qquad\ot\cd \nonumber\\
&&\qquad\quad\ot\veps_{B_d}^{-1}(\sigma_1^{j-1}\sigma_2^{j-1}\cd
\sigma_{d-1}^{j-1}\sigma_d^j\la_d). \nonumber
\end{eqnarray}
Here $\veps_{B_i}$ signifies the map $\veps$ on $B_i$ and $\sigma_i$ is the 
automorphism on $P^+$ corresponding to $B_i$.

\begin{conjecture} \label{conj:tensor-prod-th}
With the notations above, we have the following isomorphism of $P$-weighted 
crystals.
\[
\P(\bp^{(\la_1,\cd,\la_d)},B)\simeq B(\la_1)\ot\cd\ot B(\la_d).
\]
\end{conjecture}

\begin{remark}
This conjecture is intimately related to the following. There exists an isomorphism of 
$P_{cl}$-weighted crystals
\[
B(\la_1)\ot\cd\ot B(\la_d)\simeq
(B(\sigma_1\la_1)\ot\cd\ot B(\sigma_1\cd\sigma_d\la_d))\ot B.
\]
Let $u_\la$ be the highest weight element in $B(\la)$. Under this isomorphism, 
$u_{\la_1}\ot\cd\ot u_{\la_d}$ is mapped to 
$(u_{\sigma_1\la_1}\ot\cd\ot u_{\sigma_1\cd\sigma_d\la_d})\ot\bb_1$. Thus the reference 
path $\bp^{(\la_1,\cd,\la_d)}$ corresponds to the highest weight element in 
$B(\la_1)\ot\cd\ot B(\la_d)$.
\end{remark}

\begin{remark}\label{rem:One}
One can show Remark \ref{rem:coherent} and Conjecture \ref{conj:combR} imply this conjecture 
in a similar manner to the $(A^{(1)}_n, B^{1,l_i})$ case in \cite{HKKOT}.
Since Conjecture \ref{conj:combR} is true for all coherent families of perfect crystals listed in
\cite{KKM} (Remark \ref{rem:injective}), the above conjecture is true for these cases.
\end{remark}

This conjecture can be generalized to the case when each $B_i$ is of the form $B^{k_i,s_i}$
but not necessarily perfect. We prepare some notations. Set $t_{\max}=\max_{1\le i\le n}t_i$.
Later in this section we consider an affine Lie algebra $\geh$ such that $t_{\max}>1$. 
To such $\geh$ we associate a different algebra $\geh^{(i)}$ ($1\le i<t_{\max}$) by
\begin{eqnarray*}
\geh=B^{(1)}_n&\Longrightarrow&\geh^{(1)}=D^{(1)}_{n+1},\\
\geh=C^{(1)}_n&\Longrightarrow&\geh^{(1)}=A^{(1)}_{2n-1},\\
\geh=F^{(1)}_4&\Longrightarrow&\geh^{(1)}=E^{(1)}_6,\\
\geh=G^{(1)}_2&\Longrightarrow&\geh^{(1)}=B^{(1)}_3, \geh^{(2)}=D^{(1)}_4.
\end{eqnarray*}
Notice that there is a natural inclusion
\begin{equation} \label{eq:incl}
\begin{array}{ll}
\geh\hookrightarrow\geh^{(1)}&\quad\mbox{for }\geh=B^{(1)}_n,C^{(1)}_n,F^{(1)}_4,\\
\geh\hookrightarrow\geh^{(1)}\hookrightarrow\geh^{(2)}&\quad\mbox{for }\geh=G^{(1)}_2.
\end{array}
\end{equation}
Next consider an integrable $\geh^{(i)}$-module $L^{(i)}(\mu^{(i)})$ with highest weight $\mu^{(i)}$.
Note that $\mu^{(i)}$ is a weight 
with respect to the algebra $\geh^{(i)}$. Through the inclusion
(\ref{eq:incl}) one can decompose $L^{(i)}(\mu^{(i)})$ into irreducible integrable $\geh$-modules
as
\[
L^{(i)}(\mu^{(i)})=\bigoplus_\nu L(\nu)^{\ot m_\nu},
\]
where $\nu$ runs through the weight lattice $P$ of $\geh$ and $m_\nu$ is the multiplicity.
Now define $B^{(i)}(\mu^{(i)})$ by 
\[
B^{(i)}(\mu^{(i)})=\bigoplus_\nu B(\nu)^{\ot m_\nu}.
\]
Recall that $B(\nu)$ is the crystal basis of $L(\nu)$.

Consider the crystal $B=B^{k_1,s_1}\ot B^{k_2,s_2}\ot\cd\ot B^{k_d,s_d}$.
To the datum $(k_1,s_1),\cd,(k_d,s_d)$, associate a set 
\[
A=\{s_i/t_{k_i}\mid i=1,\cd,d\}.
\]
We decompose it into disjoint components.
\begin{eqnarray*}
A&=&\bigsqcup_jA[l_j],\\
A[l]&=&\{s_i/t_{k_i}\in A\mid \lceil s_i/t_{k_i}\rceil=l\}.
\end{eqnarray*}
Let $\ol{d}$ be the number of disjoint components. We arrange $l_1,\cd,l_{\ol{d}}$ so that
$l_1>\cd>l_{\ol{d}}$ and set $l_{\ol{d}+1}=0$. For any $i$ ($1\le i\le\ol{d}$) we introduce
the following crystal:

\medskip
$\bullet$ $t_{\max}=2$,
\begin{eqnarray*}
B[i]&=&B(\la_i)\ot B^{(1)}(\mu^{(1)}_i)\\
&&(\mbox{$\la_i$ : level $l_i-l_{i+1}-1$ $\geh$-weight, $\mu^{(1)}_i$ : level 1 $\geh^{(1)}$-weight})\\
&&\qquad\qquad\qquad\mbox{ if }l_i-1/2\in A[l_i],\\ 
&=&B(\la_i)\quad(\mbox{$\la_i$ : level $l_i-l_{i+1}$ $\geh$-weight})\\
&&\qquad\qquad\qquad\mbox{otherwise}.
\end{eqnarray*}

$\bullet$ $t_{\max}=3$,
\begin{eqnarray*}
B[i]&=&B(\la_i)\ot B^{(2)}(\mu^{(2)}_i)\\
&&(\mbox{$\la_i$ : level $l_i-l_{i+1}-1$ $\geh$-weight, $\mu^{(2)}_i$ : level 1 $\geh^{(2)}$-weight})\\
&&\qquad\qquad\qquad\mbox{ if both }l_i-1/3\mbox{ and }l_i-2/3\in A[l_i],\\ 
&=&B(\la_i)\ot B^{(1)}(\mu^{(1)}_i)\\
&&(\mbox{$\la_i$ : level $l_i-l_{i+1}-1$ $\geh$-weight, $\mu^{(1)}_i$ : level 1 $\geh^{(1)}$-weight})\\
&&\qquad\qquad\qquad\mbox{ if either }l_i-1/3\mbox{ or }l_i-2/3\in A[l_i],\\ 
&=&B(\la_i)\quad(\mbox{$\la_i$ : level $l_i-l_{i+1}$ $\geh$-weight})\\
&&\qquad\qquad\qquad\mbox{otherwise}.
\end{eqnarray*}
Choosing such datum $(\hat{\la},\hat{\mu})
=\{\la_i,\mu^{(t)}_i\mid i=1,\cd,\ol{d},1\le t<t_{\max}\}$, 
we set 
\[
{\mathcal B}(\hat{\la},\hat{\mu})=B[1]\ot B[2]\ot\cd\ot B[\ol{d}].
\]
By definition, it is the crystal basis of a level $l_1$ integrable $U_q(\geh)$-module,
which is not irreducible in general. If $\la_i$ is level 0, $B(\la_i)$ should be considered 
as the crystal basis of the trivial module and hence can be omitted.

\begin{conjecture} \label{conj:tensor-prod-th-np}
For any datum $(\hat{\la},\hat{\mu})$ as above, there exists 
a reference path $\bp^{(\hat{\la},\hat{\mu})}$, and we have the following
isomorphism of $P$-weighted crystals.
\[
\P(\bp^{(\hat{\la},\hat{\mu})},B)\simeq{\mathcal B}(\hat{\la},\hat{\mu}).
\]
\end{conjecture}

We give examples of ${\mathcal B}(\hat{\la},\hat{\mu})$.

\begin{example}
\begin{itemize}

\item[(i)] $\geh=C^{(1)}_3$ ($t_1=t_2=2,t_3=1,t_{\max}=2,\geh^{(1)}=A^{(1)}_5$)
\begin{eqnarray*}
&&B=B^{2,5}\ot B^{1,5}\ot B^{1,4}\ot B^{2,2}\ot B^{1,1}.\\
&&A=A[3]\sqcup A[2]\sqcup A[1],\\
&&A[3]=\{5/2\},A[2]=\{2\},A[1]=\{1,1/2\}.\\
&&{\mathcal B}(\hat{\la},\hat{\mu})=B^{(1)}(\mu^{(1)}_1)\ot B(\la_2)\ot B^{(1)}(\mu^{(1)}_3),
\end{eqnarray*}
where $\la_2$ is a level 1 $C^{(1)}_3$-weight and $\mu^{(1)}_1,\mu^{(1)}_3$ are level 1 $A^{(1)}_5$-weights.

\item[(ii)] $\geh=G^{(1)}_2$ ($t_1=1,t_2=3,t_{\max}=3,\geh^{(1)}=B^{(1)}_3,\geh^{(2)}=D^{(1)}_4$)
\begin{eqnarray*}
&&B=B^{2,8}\ot B^{2,7}\ot B^{1,1}\ot B^{2,1}.\\
&&A=A[3]\sqcup A[1],\\
&&A[3]=\{8/3,7/3\},A[1]=\{1,1/3\}.\\
&&{\mathcal B}(\hat{\la},\hat{\mu})=B(\la_1)\ot B^{(2)}(\mu^{(2)}_1)\ot B^{(1)}(\mu^{(1)}_2),
\end{eqnarray*}
where $\la_1$ is a level 1 $G^{(1)}_2$-weight, $\mu^{(2)}_1$ level 1 $D^{(1)}_4$-weight and
$\mu^{(1)}_2$ level 1 $B^{(1)}_3$-weight.
\end{itemize}
\end{example}

\def\bc{\boldsymbol{c}}

Except for the case when $B^{k_i,s_i}$ are all perfect (\ref{eq:def-boldb}),
we do not know how to construct the reference path 
$\bp^{(\hat{\la},\hat{\mu})}$ in general. We only give the answer for the case
of $\geh=C^{(1)}_n, B=B^{1,s_1}\ot\cd\ot B^{1,s_d}$ here. Let us first recall
the crystal $B^{1,s}$. As a set, it reads as 
\begin{equation*}
\begin{split}
B^{1,s}&=\{ (x_1,\dots ,x_{n},\bar{x}_{n},\dots,\bar{x}_1) \in 
{\mathbb Z}^{2n} |\,
x_i, \bar{x}_i \ge 0,\, \\ &\quad\quad
\sum_{i=1}^{n} (x_i+\bar{x}_i) \le s,\,
\sum_{i=1}^{n} (x_i+\bar{x}_i) \equiv s \pmod{2} \}.
\end{split}
\end{equation*}
For more details, see Example 2.5 of \cite{HKOTY1}.

Now consider the crystal $B=B^{1,s_1}\ot\cd\ot B^{1,s_d}$ 
($s_1>s_2>\cd>s_d>s_{d+1}=0$). 
Let $l_j=\lceil s_j/2\rceil$ and fix a datum
\begin{equation} \label{eq:la-mu}
(\hat{\la},\hat{\mu})=(\la_j,\hat{\Lambda}_{\mu_i}\mid
1\le j\le d,i\in\{j\mid s_j\mbox{ is odd }\}),
\end{equation}
where $\la_j$ is a level $l_j-l_{j+1}$ (resp. $l_j-l_{j+1}-1$) 
dominant integral weight of $C^{(1)}_n$ if $s_j$ is even (resp. odd) and
$\hat{\La}_{\mu_i}$ is a fundamental weight of $A^{(1)}_{2n-1}$. 
We need to define two maps $\pi,\iota$. $\pi$ is a map from 
$A^{(1)}_{2n-1}$ fundamental weights $\hat{\Lambda}_{\mu}$ 
($\mu=0,\dots,2n-1$) to $C^{(1)}_n$ fundamental weights given by
$\pi(\hat{\Lambda}_{\mu})=\Lambda_\mu$ ($0\le\mu\le n$), 
$=\Lambda_{2n-\mu}$ ($n<\mu\le 2n-1$). For a $C^{(1)}_n$ weight 
$\la=\sum_{i=0}^n\la^i\La_i$, define $\iota(\la)=
(\lambda^1,\cd,\lambda^n,\lambda^n,\cd,\lambda^1)$.
Let $e_k$ be the $k$-th unit row-vector of dimension $2n$.
If $2n<k<4n$, $e_k$ should be understood as $e_{k-2n}$. 
Now define an element $\bc^{(j)}$ of $(B^{1,s_j})^{\ot 2n}$ by
\begin{eqnarray*}
\bc^{(j)}&=&\bc^{(j)}_1\otimes\cdots\otimes\bc^{(j)}_{2n}\quad
(\bc^{(j)}_i\in B^{1,s_j}),\\
\bc^{(j)}_i&=&\left\{\begin{array}{ll}
\iota(\sum_{j\le k\le d}\la_k+\sum_{j+1\le k\le d,s_k: odd}
\pi(\hat{\La}_{\mu_k}))&\mbox{ if $s_j$ is even}\\
\iota(\sum_{j\le k\le d}\la_k+\sum_{j+1\le k\le d,s_k: odd}
\pi(\hat{\La}_{\mu_k}))+e_{\mu_j+i}&\mbox{ if $s_j$ is odd}.
\end{array}\right.
\end{eqnarray*}
Next applying crystal isomorphisms 
$B^{1,s_i}\ot B^{1,s_j}\simarrow B^{1,s_j}\ot B^{1,s_i}$
as many times as needed, one defines an element
$\bb_{2n}\ot\cd\ot\bb_1$ of $B^{\ot 2n}$ by
\begin{eqnarray}
(B^{1,s_1})^{\ot 2n}\ot(B^{1,s_2})^{\ot 2n}\ot\cd\ot(B^{1,s_d})^{\ot 2n}
&\simarrow&B^{\ot 2n} \label{eq:C-iso}\\
\bc^{(1)}\ot\bc^{(2)}\ot\cd\ot\bc^{(2n)}
&\mapsto&\bb_{2n}\ot\cd\ot\bb_1. \nonumber
\end{eqnarray}
Then $\bb_{2n}\ot\cd\ot\bb_1$ gives the first $2n$ components of our
reference path $\bp^{(\hat{\la},\hat{\nu})}$ corresponding to the datum
(\ref{eq:la-mu}). The other components should be determined so that
they are periodic with period $2n$. To calculate the isomorphism 
(\ref{eq:C-iso}) explicitly, one can use an algorithm based on 
insertion scheme \cite{HKOT1}. See also Appendix \ref{app:example}.

We give examples of the above mentioned $\bb_{2n}\ot\cd\ot\bb_1$
for $n=2$, i.e. $\geh=C^{(1)}_2$.

\begin{example}
\begin{itemize}

\item[(i)] $B=B^{1,4}\ot B^{1,1}$, 
$(\hat{\la},\hat{\mu})=(\la_1=\La_0,\la_2=0,\hat{\La}_{\mu_2}=\hat{\La}_0)$
\begin{eqnarray*}
&&\bc^{(1)}=(0,0,0,0)\otimes(0,0,0,0)\otimes(0,0,0,0)\otimes(0,0,0,0)\\
&&\bc^{(2)}=(1,0,0,0)\otimes(0,1,0,0)\otimes(0,0,1,0)\otimes(0,0,0,1)\\
&&\bb_1=(1,0,0,1)\otimes(0,0,0,1)\quad\bb_2=(0,1,1,0)\otimes(0,0,1,0)\\
&&\bb_3=(1,0,0,1)\otimes(0,1,0,0)\quad\bb_4=(0,0,0,0)\otimes(1,0,0,0)
\end{eqnarray*}

\item[(ii)] $B=B^{1,5}\ot B^{1,3}$, 
$(\hat{\la},\hat{\mu})=(\la_1=0,\la_2=\La_2,
\hat{\La}_{\mu_1}=\hat{\La}_3,\hat{\La}_{\mu_2}=\hat{\La}_1)$
\begin{eqnarray*}
&&\bc^{(1)}=(1,1,1,2)\otimes(2,1,1,1)\otimes(1,2,1,1)\otimes(1,1,2,1)\\
&&\bc^{(2)}=(0,2,1,0)\otimes(0,1,2,0)\otimes(0,1,1,1)\otimes(1,1,1,0)\\
&&\bb_1=(0,1,2,0)\otimes(1,1,1,0)\quad\bb_2=(1,2,1,1)\otimes(0,1,1,1)\\
&&\bb_3=(1,2,2,0)\otimes(0,1,2,0)\quad\bb_4=(1,1,1,2)\otimes(0,2,1,0)
\end{eqnarray*}

\end{itemize}
\end{example}

\subsection{One dimensional sums}\label{subsec:sums}

Let $B_i$ ($i=1,\cd,d$) be finite crystals, and consider the tensor
product $B_1\ot\cd\ot B_d$.
We define $b^{(i)}_j$ ($i<j$) by
\begin{eqnarray*}
&&
\begin{array}{ccccc}\hspace{-5mm}
B_i\ot\cd\ot B_{j-1}\ot B_j&\simarrow&
B_i\ot\cd\ot B_j\ot B_{j-1}&\simarrow&\cd\\
b_i\ot\cd\ot b_{j-1}\ot b_j&\mapsto&
b_i\ot\cd\ot b^{(j-1)}_j\ot b'_{j-1}&\mapsto&\cd
\end{array}\\
&&\hspace{5cm}
\begin{array}{ccc}
\cd&\simarrow&B_j\ot B_i\ot\cd\ot B_{j-1}\\
\cd&\mapsto&b^{(i)}_j\ot b'_i\ot\cd\ot b'_{j-1},
\end{array}
\end{eqnarray*}
and set $b^{(i)}_i=b_i$. 
Suppose $B_i=B^{k_i,s_i}$ ($1\le i\le d$) and let $l_i$ be the level of 
the crystal $B^{k_i,s_i}$. Take $b_i^\natural\in B^{k_i,s_i}$ such
that $\vphi(b_i^\natural)=l_i\La_0$. 
Even if $B^{k_i,s_i}$ is not perfect, we expect such $b_i^\natural$ exists uniquely.
It is indeed true for the $C^{(1)}_n$-crystal 
$B^{1,s}$ with odd $s$ presented in \cite{HKKOT}.
For later use we also define $b^{ext}_i\in B^{k_i,s_i}$ to be the highest
weight element of the $U_q(\overset{\circ}{\geh})$-crystal
$B(s_i\Lab_{k_i})$ contained in $B^{k_i,s_i}$.
For a finite path $p=b_1\ot\cd\ot b_d\in
B^{k_1,s_1}\ot\cd\ot B^{k_d,s_d}$, define
\begin{equation}\label{eq:subtle}
D(p)=\sum_{1\le i<j\le d}H(b_i\ot b^{(i+1)}_j)+
\sum_{1\le j\le d}H(b_j^\natural\ot b^{(1)}_j).
\end{equation}
Fix a positive integer $l$ such that $l\ge\max_i l_i$.
We define a (restricted)
one dimensional sum by 
\begin{equation}\label{eq:defX}
X_l(B^{k_1,s_1}\ot\cd\ot B^{k_d,s_d},\la,q)=\mathop{{\sum}^*}_p q^{D(p)},
\end{equation}
where the sum $\sum_p^*$ is taken over all 
$p=b_1\ot\cd\ot b_d\in B^{k_1,s_1}\ot\cd\ot B^{k_d,s_d}$ satisfying $\sum_{j=1}^d\wt b_j=\la$
and
\begin{equation} \label{eq:hwc}
{\tilde e}_i^{l\delta_{i0}+1}(b_1\ot\cd\ot b_d)=0\quad\forall i.
\end{equation}
(\ref{eq:hwc}) can be rewritten 
as local conditions.
\[
{\tilde e}_i^{\langle h_i,l\La_0+\wts b_1+\cd+\wts b_{j-1}\rangle+1}b_j=0
\quad\forall i,1\le j\le d.
\]
We formally consider the $l=\infty$ case, where the condition (\ref{eq:hwc}) for $i=0$ is vacant.
$X_\infty(B,\la,q)$ was called classically restricted in \cite{HKOTY1}.
It is easy to see {}from the definition that
\begin{equation} \label{eq:count-mult}
X_\infty(B^{k_1,s_1}\ot\cd\ot B^{k_d,s_d},\la,1)
=[W^{(k_1)}_{s_1}\ot\cd\ot W^{(k_d)}_{s_d}:\la].
\end{equation}
Here $W^{(k)}_s$ is the corresponding finite-dimensional $U'_q(\geh)$-module 
of $B^{k,s}$ and 
\[
[M:\mu]=\dim\C\langle v\in M\mid\wt v=\mu,e_iv=0\,\forall i\neq0\rangle,
\]
where $e_i$ denotes the root vector corresponding to the simple root $\alpha_i$.

\begin{proposition}
$X_l(B^{k_1,s_1}\ot\cd\ot B^{k_d,s_d},\la,q)$ does not depend on the order of 
the crystals $B^{k_1,s_1},\cd,B^{k_d,s_d}$ in the tensor product.
\end{proposition}

\begin{proof}
We set $B_i=B^{k_i,s_i}(i=1,\cd,d)$. It is enough to check that $X_l$ remains
unchanged under the switch of $B_i$ and $B_{i+1}$ in the tensor product. 
Suppose $b_i\ot b_{i+1}$ is mapped to $b'_{i+1}\ot b'_i$ under the isomorphism
$B_i\ot B_{i+1}\simeq B_{i+1}\ot B_i$. 
For $p=b_1\ot\cd\ot b_d\in B_1\ot\cd\ot B_d$ define $p'=b_1\ot\cd\ot b'_{i+1}\ot
b'_i\ot\cd\ot b_d\in B_1\ot\cd\ot B_{i+1}\ot B_i\ot\cd\ot B_d$. Note that 
$\wt p=\wt p'$ and $p$ satisfies (\ref{eq:hwc}) if and only if so does $p'$.
Hence it is enough to show $D(p)=D(p')$. Set
\[
D_j(p)=\sum_{k=1}^{j-1}H(b_k\ot b^{(k+1)}_j),\quad
{\hat D}(p)=\sum_{j=1}^d H(b^\natural_j\ot b^{(1)}_j).
\]
Then $D(p)=\sum_{j=1}^d D_j(p)+{\hat D}(p)$. Since ${\hat D}(p)={\hat D}(p')$
is evident, we are left to show the following:
\begin{itemize}
\item[(i)] $D_j(p)=D_j(p')\qquad(j\neq i,i+1)$,
\item[(ii)] $D_i(p)+D_{i+1}(p)=D_i(p')+D_{i+1}(p')$.
\end{itemize}
(i) for $j<i$ is trivial. To check (i) for $j>i+1$ use (\ref{eq:YBE}),
in particular, (\ref{eq:used-later}) with $b_1\rightarrow b_i,
b_2\rightarrow b_{i+1},b_3\rightarrow b^{(i+2)}_j$. For (ii) notice that
$D_{i+1}(p)=D_i(p')+H(b_i\ot b_{i+1}),D_{i+1}(p')=D_i(p)+H(b'_{i+1}\ot b'_i)$.
\end{proof}

\begin{conjecture} \label{conj:X=M}
We have
\begin{itemize}
\item[(1)] $q^{-D^{ext}}X_\infty(B^{k_1,s_1}\ot\cd\ot B^{k_d,s_d},\la,q)
        =M_\infty(W^{(k_1)}_{s_1}\ot\cd\ot W^{(k_d)}_{s_d},\la,q)$

        \hspace{8.5cm}$\quad\mbox{for }\la\in\ol{P}^+$,
\item[(2)] $q^{-D^{ext}}X_l(B^{k_1,s_1}\ot\cd\ot B^{k_d,s_d},0,q)
        =M_l(W^{(k_1)}_{s_1}\ot\cd\ot W^{(k_d)}_{s_d},0,q)$,
\end{itemize}
where $D^{ext}=D(b^{ext}_1\ot\cd\ot b^{ext}_d)$ and $M_l$ is defined 
in (\ref{eq:mm}).
\end{conjecture}
At present a proof is available only for  $A^{(1)}_n$ case \cite{KSS,SS1}.

\subsection{Representation theoretical meaning}\label{subsec:representation}

We explain a representation theoretical meaning of our one dimensional sum
$X_l(B,\la,q)$. Set $B=B_1\ot\cd\ot B_d$. We first assume all $B_i$'s  are perfect.
Let $l_i$ be the level of $B_i$. Assume $l_1\ge l_2\ge\cd\ge l_d\ge l_{d+1}=0$.
Write $b^*=b_1\ot\cd\ot b_d$ for an element of $B$ and fix $\bb_2$ as in (\ref{eq:def-boldb})
with $\la_i=\sigma_1^{-1}\cd\sigma_i^{-1}(l_i-l_{i+1})\La_0$. 

\begin{proposition} \label{prop:H-D}
If Conjecture \ref{conj:combR} is true, then $H_{BB}(\bb_2\ot b^*)-D(b^*)$ is independent of
$b^*$.
\end{proposition}

We need three lemmas.

\begin{lemma}[\cite{HKKOT} Proposition 3.1] \label{lem:1}
Let $B_1,B_2$ be finite crystals.
Set $B=B_1\ot B_2$, then
\begin{eqnarray*}
H_{BB}((b_1\ot b_2)\ot(b'_1\ot b'_2))
&=&H_{B_1B_2}(b_1\ot b_2)+H_{B_1B_1}(\bt_1\ot b'_1)\\
&&+H_{B_2B_2}(b_2\ot \bt'_2)+H_{B_1B_2}(b'_1\ot b'_2).
\end{eqnarray*}
Here $\bt_1,\bt'_2$ are defined as
\begin{eqnarray*}
B_1\ot B_2&\simarrow&B_2\ot B_1\\
b_1\ot b_2&\mapsto&\bt_2\ot\bt_1\\
b'_1\ot b'_2&\mapsto&\bt'_2\ot\bt'_1.
\end{eqnarray*}
\end{lemma}

\begin{lemma} \label{lem:2}
Consider the tensor product of finite crystals $B=B_1\ot B_2\ot\cd\ot B_d$.
Set $B_{>1}=B_2\ot\cd\ot B_d$. For an element $b_1\ot\cd\ot b_d\in B$,
define $b_1^{\langle j\rangle}$ ($1\le j\le d$) by
\[
\begin{array}{ccc}
B_1\ot B_2\ot\cd\ot B_j&\simarrow&B_2\ot\cd\ot B_j\ot B_1\\
b_1\ot b_2\ot\cd\ot b_j&\mapsto&\bt_2\ot\cd\ot\bt_j\ot b_1^{\langle j\rangle}.
\end{array}
\]
We understand $b_1^{\langle 1\rangle}=b_1$.
Then we have 
\[
H_{B_1B_{>1}}(b_1\ot(b_2\ot\cd\ot b_d))=\sum_{2\le j\le d}H_{B_1B_j}(b_1^{\langle j-1\rangle}\ot b_j).
\]
\end{lemma}

\begin{lemma} \label{lem:3}
We retain the notations in the previous lemma. Consider isomorphisms of crystals
\[
\begin{array}{lcl}
B_1\ot B_2\ot\cd\ot B_{j-1}\ot B_j&\simarrow&B_2\ot\cd\ot B_{j-1}\ot B_j\ot B_1\\
b_1\ot b_2\ot\cd\ot b_{j-1}\ot b_j&\mapsto&\bt_2\ot\cd\ot \bt_{j-1}\ot\bt_j\ot b_1^{\langle j\rangle}\\
&&\quad\cd\\
&\simarrow&B_2\ot\cd\ot B_{i-1}\ot B_j\ot B_i\ot\cd\ot B_{j-1}\ot B_1\\
&\mapsto&\bt_2\ot\cd\ot \bt_{i-1}\ot\bt_j^{(i)}\ot\cd\\
&&\quad\cd\\
&\simarrow&B_j\ot B_2\ot\cd\ot B_{j-1}\ot B_1\\
&\mapsto&\bt_j^{(2)}\ot\cd
\end{array}
\]
Then we have $\bt_j^{(2)}=b_j^{(1)}$ and $\sum_{i=2}^{j-1}H_{B_iB_j}(\bt_i\ot\bt_j^{(i+1)})
+H_{B_1B_j}(b_1^{\langle j-1\rangle}\ot b_j)=\sum_{i=1}^{j-1}H_{B_iB_j}(b_i\ot b_j^{(i+1)})$.
\end{lemma}

The last two lemmas can be shown by using the Yang-Baxter equation (\ref{eq:YBE}).

\medskip
\noindent{\it Proof of Proposition \ref{prop:H-D}.}\;
The notations in the lemmas are retained. 
Write $\bb_2={\bar b}_1\ot\cd\ot{\bar b}_d,{\bar b}_i\in B_i$. Let 
${\bar b}_i\ot{\bar b}_{i+1}\ot\cd\ot{\bar b}_d$ be mapped to 
${\bar b}'_{i+1}\ot\cd\ot{\bar b}'_d\ot{\bar b}'_i$ under the isomorphism
$B_i\ot B_{i+1}\ot\cd\ot B_d\simarrow B_{i+1}\ot\cd\ot B_d\ot B_i$.
Then Conjecture \ref{conj:combR} implies ${\bar b}'_i=b_i^\natural$.

We prove by induction on $d$. Using Lemma \ref{lem:1} and the induction hypothesis, 
one has
\begin{eqnarray*}
H_{BB}(\bb_2\ot b^*)&=&H_{B_1B_1}(b_1^\natural\ot b_1)+\sum_{2\le i<j\le d}H_{B_iB_j}(\bt_i\ot \bt_j^{(i+1)})\\
&&+\sum_{2\le j\le d}H_{B_jB_j}(b_j^\natural\ot \bt_j^{(2)})+H_{B_1B_{>1}}(b_1\ot(b_2\ot\cd\ot b_d)).
\end{eqnarray*}
Here in the equality we omitted $b^*$-independent terms.
Lemma \ref{lem:2} and \ref{lem:3} now complete the proof.
\qed

Proposition \ref{prop:H-D} implies that under the embedding
\[
\begin{array}{ccl}
B&\hookrightarrow&\P(\bp^{(\la_1,\cd,\la_d)},B)\simeq B(\la_1)\ot\cd\ot B(\la_d)\\
b^*&\mapsto&\cd\ot\bb_j\ot\cd\ot\bb_2\ot b^*,
\end{array}
\]
the weight is given by $l_1\La_0+\wt b^*-(D(b^*)-D(\bb_1))\delta$.
Here $\la_i=\sigma_1^{-1}\cd\sigma_i^{-1}(l_i-l_{i+1})\La_0$
and $\bb_j$ is the $j$-th component of the reference path $\bp^{(\la_1,\cd,\la_d)}$. 
Let $\gamma$ be the minimal integer such that $\sigma_i^\gamma=\id$ for all 
$i=1,\cd,d$. For a positive integer $L$ such that $L\equiv0$ (mod $\gamma$), we
consider a similar embedding
\[
\begin{array}{ccl}
B^{\ot L}&\hookrightarrow&\P(\bp^{(\la^0_1,\cd,\la^0_d)},B)\simeq B(\la^0_1)\ot\cd\ot B(\la^0_d)\\
b^L\ot\cd\ot b^1&\mapsto&\cd\ot\bb_j\ot\cd\ot\bb_{L+1}\ot b^L\ot\cd\ot b^1.
\end{array}
\]
This time $\la^0_i=(l_i-l_{i+1})\La_0$ and $\bb_j$ is redefined as the $j$-th component of
the reference path $\bp^{(\la^0_1,\cd,\la^0_d)}$. Under the embedding the weight is given
by $l_1\La_0+\sum_{j=1}^L\wt b^j-(D(b^L\ot\cd\ot b^1)-D(\bb_L\ot\cd\ot\bb_1))\delta$. 
In view of (\ref{eq:hwc}), one has
\begin{eqnarray} 
&&\lim_{L\rightarrow\infty\atop L\equiv0\,(mod\,\gamma)}
q^{-D(\bb_L\ot\cd\ot\bb_1)}X_\infty(B^{\ot L},\la,q) \nonumber\\
&&\hspace{2cm}
=\sum_j[\bigotimes_{i=1}^dL((l_i-l_{i+1})\La_0):l_1\La_0+\la-j\delta]q^{ja_0}.
\label{eq:lim-Xinf}
\end{eqnarray}
Note that in $[\bigotimes_{i=1}^d
L(\la^0_i):\mu]$, $\mu$ is a weight in $P$ 
as opposed to $\la$ in (\ref{eq:count-mult}) in $\ol{P}$.
Let $l$ be an integer such that $l\ge l_1$. We also have 
\begin{equation} 
\begin{gathered}
\lim_{L\rightarrow\infty\atop L\equiv0\,(mod\,\gamma)}
q^{-D(\bb_L\ot\cd\ot\bb_1)}X_l(B^{\ot L},\la,q)\\
=\sum_j[[L((l-l_1)\La_0)\ot\bigotimes_{i=1}^dL((l_i-l_{i+1})\La_0):
l\La_0+\la-j\delta]]q^{ja_0}.
\end{gathered}\label{eq:lim-Xl}
\end{equation}
Here $[[M,\mu]]=\dim\C\langle v\in M\mid\wt v=\mu,e_iv=0\,\forall i\rangle$.

Generalization of such interpretation is straightforward to the case when $B_i$ is not
necessarily perfect. In view of Conjecture \ref{conj:tensor-prod-th-np},
let ${\mathcal L}(\hat{\la},\hat{\mu})$ be the corresponding integrable $U_q(\geh)$-module
of ${\mathcal B}(\hat{\la},\hat{\mu})$. One only needs to replace $\bigotimes_{i=1}^d
L((l_i-l_{i+1})\La_0)$ with ${\mathcal L}(\hat{\la},\hat{\mu})$.

\section{\mathversion{bold} Fermionic formula}\label{sec:FF}

\subsection{\mathversion{bold} Conventions on $\geh_{\ol{0}}$}
\label{subsec:convention}

{}From Section \ref{subsec:affine} we remind that 
$\{\ol{\La}_a\}_{a=1}^n$ and $\{\alpha_a\}_{a=1}^n$ 
are the sets of fundamental weights and simple roots of 
$\overset{\circ}{\geh}$, and 
$\ol{P} = \Z \ol{\La}_1 \oplus \cd \oplus \Z \ol{\La}_n$ denotes its 
weight lattice.
$\{\alpha_a\}$ has been normalized as 
$(\alpha_a \vert \alpha_a) = 2r$ if $\alpha_a$ is a long root.
See Exercise 6.1 of \cite{Kac}.
In Sections \ref{sec:FF}, \ref{sec:N} and \ref{sec:qsys} 
we find it convenient to further introduce  
the simple roots $\{ \tilde{\alpha}_a\}_{a=1}^n$ 
and the fundamental weights 
$\{\tilde{\La}_a\}_{a=1}^n$  of $\geh_{\ol{0}}$, and let 
$\tilde{P} = \Z \tilde{\La}_1 \oplus \cd \oplus \Z \tilde{\La}_n$.
We equip the bilinear form $(\cdot \vert \cdot )'$ on 
$\tilde{P}$ which is again 
normalized as 
$(\tilde{\alpha}_a\vert\tilde{\alpha}_a )' = 2r$ 
if $\tilde{\alpha}_a$ is a long root of $\geh_{\ol{0}}$.
(For $A^{(2)}_2$, we understand that the unique simple root 
$\tilde{\alpha}_1$ of $\geh_{\ol{0}}=B_1$ is short. See below.)
Actually the situation $\overset{\circ}{\geh} \neq \geh_{\ol{0}}$
happens only for  $\geh = A^{(2)}_{2n}$, where 
$\overset{\circ}{\geh}=C_n$ and $\geh_{\ol{0}}=B_n$.
Comparing the normalizations in the above and 
Section \ref{subsec:affine}, we find that 
$\alpha_a, \ol{\La}_a$ and $(\cdot \vert \cdot )$ may be identified 
with $\tilde{\alpha}_a, \tilde{\La}_a$ and $(\cdot \vert \cdot )'$ 
if $\geh \neq A^{(2)}_{2n}$.
In order to deal with the $A^{(2)}_{2n}$ case uniformly,
we introduce a $\Z$-linear map 
$\iota : \ol{P} \rightarrow \tilde{P}$ by
\begin{equation}\label{eq:iota}
\iota(\ol{\La}_a) = 
\epsilon_a \tilde{\La}_a \quad 1 \le a \le n,
\end{equation}
where $\epsilon_a$ is specified as
\begin{equation}\label{eq:epsdef}
\epsilon_a = \begin{cases}
2 & \text{if $\geh=A^{(2)}_{2n}$ and $a=n$}\\
1 & \text{otherwise}
\end{cases}.
\end{equation}
Note that (\ref{eq:iota}) induces 
$\iota(\alpha_a) = \epsilon_a \tilde{\alpha}_a$.
One can check 
$(\iota(\alpha_b)\vert \iota(\alpha_b))' 
= a_0(\alpha_b \vert \alpha_b)$ for any $1 \le b \le n$.
For example in $A^{(2)}_{2n}$ case,
the both sides are equal to $8$ if $b=n$ and $4$ otherwise.
Especially for $A^{(2)}_2 \, (n=1)$,
we have 
$(\tilde{\alpha}_1\vert\tilde{\alpha}_1)' =2$ and 
$(\alpha_1 \vert \alpha_1)=4$.
Now that the definitions being clear, 
in the rest of the paper 
(mostly in the rest of this section and 
Sections \ref{sec:N} and \ref{sec:qsys}) 
we shall simply write $(\cdot \vert\cdot)'$ 
also as $(\cdot\vert\cdot)$.
This will not cause a confusion because 
the entries shall be limited to the self-explanatory objects 
$\alpha_a, \ol{\La}_a, 
\tilde{\alpha}_a, \tilde{\La}_a, \iota(\ol{\La}_a)$, etc.
We will also make abbreviations 
$\vert \mu \vert^2 = (\mu\vert\mu)$ for $\mu$ 
belonging to either $\overset{\circ}{\geh}$ or 
$\geh_{\ol{0}}$ root systems.
With these conventions one finds
\begin{equation*}
\tilde{\alpha}_a = \sum_{b=1}^n \frac{\epsilon_bt_b}{t^\vee_b}
(\tilde{\alpha}_a\vert\tilde{\alpha}_b)\tilde{\La}_b,
\qquad
{\alpha}_a = \sum_{b=1}^n \frac{a_0t_b}{\epsilon_bt^\vee_b}
({\alpha}_a\vert{\alpha}_b)\ol{\La}_b,
\end{equation*}
which  tells that the matrices 
$\left(\frac{\epsilon_a t_a}{t^\vee_a}(\tilde{\alpha}_a\vert\tilde{\alpha}_b)
\right)_{1 \le a,b \le n}$
and
$\left(\frac{a_0 t_a}{\epsilon_at^\vee_a}({\alpha}_a\vert{\alpha}_b)
\right)_{1 \le a,b \le n}$
coincide with the Cartan matrices of 
$\geh_{\ol{0}}$ and $\overset{\circ}{\geh}$,
respectively for any $\geh$.
In particular, 
\begin{equation}\label{eq:account}
\frac{2t^\vee_b}{t_b\vert \tilde{\alpha}_b\vert^2}
= \frac{a_0t_b\vert\alpha_b\vert^2}{2t^\vee_b} = \epsilon_b\quad 
1 \le b \le n
\end{equation}
is valid.
\subsection{\mathversion{bold} $M_l(W,\la,q)$ and 
$\tilde{M}_l(W,\la,q)$}\label{subsec:Ml}

For $m \in {\mathbb Z}_{\ge 0}$ and $p \in {\mathbb Z}$,
we define 
\begin{equation}\label{eq:qbinomial}
\left[\begin{array}{c} p + m \\ m \end{array} \right]_q  = 
\frac{(q^{p+1};q)_\infty(q^{m+1};q)_\infty}
{(q;q)_\infty(q^{p+m+1};q)_\infty}, 
\end{equation}
where $(x;q)_\infty = \prod_{j=0}^\infty(1-xq^j)$.
The quantity (\ref{eq:qbinomial}) is the 
$q$-binomial coefficient for $p \in \Z_{\ge 0}$, 
vanishing for $-m \le p \le -1$, 
and is equal to 
$(-q^{p+(m+1)/2})^m\left[\begin{array}{c} -p-1 \\ 
m \end{array} \right]_q$
for $p \le -m-1$.
In the $q \rightarrow 1$ limit it becomes
\begin{equation*}
\left[\begin{array}{c} p + m \\ m \end{array} \right]_1 = 
\frac{\Gamma(p+m+1)}{\Gamma(p+1) \Gamma(m+1)}.
\end{equation*}
We shall also use the notation 
$(q)_k = (q;q)_\infty/(q^{k+1};q)_\infty$ for $k \in {\mathbb Z}_{\ge0}$.

For any $l \in \Z_{\ge 1}$, we introduce the index set 
\begin{align*}
H_l &= \{(a,i) \mid 1 \le a \le n, 1 \le i \le t_al\},\\
\ol{H}_l &= \{(a,i) \mid 1 \le a \le n, 1 \le i \le t_al-1\}.
\end{align*}
We will also use $H_l$ formally with $l = \infty$.

Given  $\{ \nu^{(a)}_j \in {\mathbb Z}_{\ge 0} \mid 
(a,j) \in H_l \}$, we set 
$W = \bigotimes_{(a,j) \in H_l}
\bigl(W^{(a)}_j\bigr)^{\otimes \nu^{(a)}_j}$.
When $l = \infty$ we assume that $\sum_{j \ge 1}\nu^{(a)}_j < \infty$ 
for all $1\le a \le  n$.

Let $\la = \la_1 \ol{\La}_1 + \cd + \la_n \ol{\La}_n$ be an element of 
$\ol{P}$.
We define the fermionic form $M_l(W,\la,q)$ for $\geh$ by
\begin{eqnarray}
M_l(W,\lambda,q) & = & \sum_{\{m \}} q^{c(\{m\})}
\prod_{(a,i) \in H_l}
\left[ \begin{array}{c} p^{(a)}_i +  m^{(a)}_i
 \\   m^{(a)}_i \end{array} \right]_{q_a},
 \label{eq:mm}\\
c(\{m\}) & = & \frac{1}{2}
\sum_{(a,j),(b,k) \in H_l} (\tilde{\alpha}_a \vert \tilde{\alpha}_b) 
\mbox{min}(t_bj,t_ak) m^{(a)}_j m^{(b)}_k \label{eq:mc}\\
&& \qquad\quad -\sum_{a=1}^n t^\vee_a 
\sum_{1 \le j,k \le t_al} \mbox{min}(j, k)\nu^{(a)}_j
m^{(a)}_k,\nonumber\\
p^{(a)}_i& = & \sum_{j = 1}^{t_al}\nu^{(a)}_j\mbox{min}(i,j)
- \frac{1}{t^\vee_a}\sum _{(b,k) \in H_l}  
(\tilde{\alpha}_a \vert \tilde{\alpha}_b) \mbox{min}(t_bi,t_ak) m^{(b)}_k,
\label{eq:mp}
\end{eqnarray}
where $q_a$ is given by
\begin{equation}\label{eq:mqt}
q_a = q^{t^\vee_a}.
\end{equation}
The sum $\sum_{\{ m \}}$ is taken over
$\{ m^{(a)}_i \in {\mathbb Z}_{\ge 0} \mid (a,i) \in H_l \}$
satisfying 
\begin{eqnarray}
&p^{(a)}_j \ge 0 \mbox{ for all } (a,j) \in H_l, \label{eq:p-positive1}\\
&\sum_{(a,i) \in H_l} i m^{(a)}_i \tilde{\alpha}_a = 
\iota\left(\sum_{(a,i) \in H_l} i\nu^{(a)}_i \ol{\La}_a - \lambda
\right).\label{eq:wmc}
\end{eqnarray}
{}From (\ref{eq:mp}) it is easy to see that  (\ref{eq:wmc}) is equivalent to 
$p^{(a)}_{t_al} = \la_a$ 
for any $\geh$.
In view of this and the constraint (\ref{eq:p-positive1})  one has 
\begin{equation}\label{eq:mzero}
M_l(W,\la,q) = 0\quad\mbox{unless } \la \in 
\Bigl(\sum_{(a,i) \in H_l} i \nu^{(a)}_i \ol{\La}_a
 - \sum_{a=1}^n {\mathbb Z}_{\ge 0} \iota^{-1}(\tilde{\alpha}_a)\Bigr) \cap 
\ol{P}^+.
\end{equation}
It is easy to see the properties
\begin{equation*}
p^{(a)}_i = -\frac{1}{t^\vee_a}\frac{\partial c(\{m\})}{\partial m^{(a)}_i},
\qquad
M_l(W,\la,q) \in q^{-\epsilon}\mathbb{Z}_{\ge 0}[q^{-a_0}],
\end{equation*}
where $\epsilon \equiv 
\sum_{a=1}^na(\sum_{i=1}^li\nu^{(a)}_i-\la_a)$
mod $a_0\Z$.

Except (\ref{eq:MN}) in the next section, 
we will actually consider 
$M_l(W,\la,q)$ only for $\la = 0$ when  $l < \infty$.
Using $p^{(a)}_{t_al} = 0$ and 
eliminating  $m^{(a)}_{t_al}$ by 
\begin{equation}\label{eq:eliminate}
\sum_{(a,i) \in H_l} i m^{(a)}_i \tilde{\alpha}_a = 
\sum_{(a,i) \in H_l} i\nu^{(a)}_i \iota(\ol{\La}_a),
\end{equation}
one can rewrite the fermionic form for $l < \infty$ as 
\begin{eqnarray}
M_l(W,0,q^{-1}) & = & \sum_{\{m \}} q^{\ol{c}(\{m\})}
\prod_{(a,i) \in \ol{H}_{l}}
\left[ \begin{array}{c} p^{(a)}_i +  m^{(a)}_i
 \\   m^{(a)}_i \end{array} \right]_{q_a},
 \label{eq:mm2}\\
\ol{c}(\{m\}) & = & \frac{1}{2}
\sum_{(a,j),(b,k) \in \ol{H}_{l}} (\tilde{\alpha}_a \vert \tilde{\alpha}_b) 
K^{(t_at_bl)}_{t_bj,t_ak} m^{(a)}_j m^{(b)}_k \label{eq:mc2}\\
&& \qquad\quad + \frac{1}{2l}\vert \sum_{(a,i)\in H_l} 
i\nu^{(a)}_i \iota(\ol{\La}_a) \vert^2 \nonumber\\
p^{(a)}_i& = & \sum_{j = 1}^{t_al-1}K^{(t_al)}_{i,j}\nu^{(a)}_j
- \frac{1}{t^\vee_a}\sum _{(b,k) \in \ol{H}_{l}} (\tilde{\alpha}_a \vert \tilde{\alpha}_b) 
K^{(t_at_bl)}_{t_bi,t_ak} m^{(b)}_k,
\label{eq:mp2}
\end{eqnarray}
where the quantity $K^{(l)}_{j k}$ is defined by
\begin{equation}\label{eq:kdef}
K^{(l)}_{j k} = \mbox{min}(j,k)-\frac{jk}{l}.
\end{equation}
The sum $\sum_{\{ m \}}$ is now taken over
$\{ m^{(a)}_i \in {\mathbb Z}_{\ge 0} \mid (a,i) \in \ol{H}_{l} \}$
satisfying 
$p^{(a)}_i \ge 0$ for all $(a,i) \in \ol{H}_{l}$ and 
the condition that $m^{(a)}_{t_al}$ determined {}from 
(\ref{eq:eliminate}) is a non-negative integer.
When $l=\infty$, a similar calculation leads to
\begin{align}
&M_\infty(\bigotimes_{a,j}W^{(a)\; \otimes \nu^{(a)}_j}_j,
\sum_{a,j}j\nu^{(a)}_j\ol{\La}_a-\mu,q^{-1}) \nonumber \\
&=  \sum_{\{m \}} q^{\hat{c}(\{m\})}
\prod_{1\le a \le n, i \ge 1}
\left[ \begin{array}{c} p^{(a)}_i +  m^{(a)}_i
 \\   m^{(a)}_i \end{array} \right]_{q_a},
 \label{eq:mm3}\\
&\hat{c}(\{m\})  =  \frac{1}{2}
\sum_{1 \le a,b \le n} \sum_{j,k \ge 1}(\tilde{\alpha}_a \vert \tilde{\alpha}_b) 
\min(t_bj,t_ak)m^{(a)}_j m^{(b)}_k,\label{eq:mc3}\\
&p^{(a)}_i =  \sum_{j \ge 1}\nu^{(a)}_j\mbox{min}(i,j)
- \frac{1}{t^\vee_a}\sum_{b=1}^n\sum_{k\ge 1}  
(\tilde{\alpha}_a \vert \tilde{\alpha}_b) \mbox{min}(t_bi,t_ak) m^{(b)}_k.
\label{eq:mp3}
\end{align}
Here the sum $\sum_{\{ m \}}$ is taken over
$\{ m^{(a)}_i \in {\mathbb Z}_{\ge 0} \mid 1 \le a \le n, i \ge 1 \}$
satisfying 
$p^{(a)}_i \ge 0$ for all $1\le a \le n, i \ge 1$ and 
$\sum_{a=1}^n\sum_{j\ge 1} j m^{(a)}_j \tilde{\alpha}_a = \iota(\mu)$.

\medskip

Let us present the examples {}from $\geh$ of rank 1.
\begin{example}
Take $\geh=A^{(1)}_1$ and $W = \bigotimes_j W_j^{\otimes \nu_j}$.
(We understand $W_j = W^{(1)}_j$.)  Then we have 
\begin{eqnarray*}
M_\infty(W,\lambda_1\ol{\La}_1,q^{-1}) & = & \sum q^{\sum_{i,j \ge 1}\min(i,j)m_im_j}
\prod_{i \ge 1}
\left[ \begin{array}{c} p_i +  m_i
 \\   m_i \end{array} \right]_{q},\\
p_i& = & \sum_{j\ge 1}(\nu_j - 2m_j)\min(i,j).
\end{eqnarray*}
Here $\lambda_1$ is a non-negative integer congruent to $\sum_jj\nu_j$ mod 2.
The sum is taken over $m_1, m_2, \ldots  \in \Z_{\ge 0}$ such that 
$\forall p_j \ge 0$ and 
$\sum_{j\ge 1}jm_j = (\sum_j j \nu_j-\lambda_1)/2$.
For $W = W_1^{\otimes 2} \ot W_2$, the non-zero ones are given by
\begin{equation*}
M_\infty(W,4\ol{\La}_1,q^{-1}) = 1,\quad
M_\infty(W,2\ol{\La}_1,q^{-1}) = q+q^2,\quad
M_\infty(W,0,q^{-1}) = q^2.
\end{equation*}
\end{example}

\begin{example}
Take $\geh=A^{(2)}_2$ and $W = \bigotimes_j W_j^{\otimes \nu_j}$.
(We understand $W_j = W^{(1)}_j$.)  Then we have 
\begin{eqnarray*}
M_\infty(W,\lambda_1\ol{\La}_1,q^{-1}) & = & \sum q^{\sum_{i,j \ge 1}\min(i,j)m_im_j}
\prod_{i \ge 1}
\left[ \begin{array}{c} p_i +  m_i
 \\   m_i \end{array} \right]_{q^2},\\
p_i& = & \sum_{j\ge 1}(\nu_j - m_j)\min(i,j).
\end{eqnarray*}
Here $\lambda_1$ is a non-negative integer.
The sum is taken over $m_1, m_2, \ldots  \in \Z_{\ge 0}$ such that 
$\forall p_j \ge 0$ and 
$\sum_{j\ge 1}jm_j = \sum_j j \nu_j-\lambda_1$.
For $W = W_1^{\otimes 2} \ot W_2$, the non-zero ones are give by
\begin{alignat*}{2}
&M_\infty(W,4\ol{\La}_1,q^{-1}) = 1,& &\\
&M_\infty(W,3\ol{\La}_1,q^{-1}) = q+q^3+q^5, &\quad
&M_\infty(W,2\ol{\La}_1,q^{-1}) = q^2 + 2q^4 + 2q^6 + q^8, \\
&M_\infty(W,\ol{\La}_1,q^{-1}) = q^3 + 2q^5 + 2q^7 + 2q^9, &\quad
&M_\infty(W,0,q^{-1}) = q^4 + q^6 + 2q^8 + q^{10}. 
\end{alignat*}
\end{example}

It is natural to consider a relative of $M_l$ 
defined by\footnote{$\tilde{M}_l$ here was denoted by $N_l$ in eq.(4.16) of 
\cite{HKOTY1}.}
\begin{equation}\label{eq:Mtilde}
\tilde{M}_l(W,\la,q) = (\ref{eq:mm}) - (\ref{eq:wmc}) 
\text{ {\em without} the constraint } (\ref{eq:p-positive1}).
\end{equation}

\subsection{\mathversion{bold} Symmetry and recursion relation}
\label{subsec:symmetry}

Compared with $M_l$ in (\ref{eq:mzero}),
$\tilde{M}$ introduced in (\ref{eq:Mtilde}) has a nice symmetry 
as emphasized in \cite{HKOTY1} when $r=1$.
In fact our computer experiments indicate\footnote{
The former equality in Conjecture \ref{conj:MM} 
reduces to eq.(4.22) in \cite{HKOTY1} 
when $\la=0$. We withdraw the claim in the last sentence on 
p.260 of \cite{HKOTY1}.}
\begin{conjecture}\label{conj:MM}
\begin{align*}
M_l(W,\la,q) &= \tilde{M}_l(W,\la,q),\quad \text{for $l \le \infty$ and }
\la \in \ol{P}^+,\\
\tilde{M}_\infty(W,w(\la+\ol{\rho})-\ol{\rho},q) &= 
\det w \ \tilde{M}_\infty(W,\la,q) 
\quad \mbox{ for any } \la \in \ol{P},
\end{align*}
where $\ol{\rho} = \ol{\La}_1 + \cd + \ol{\La}_n$ and 
$w$ is any element of the Weyl group of 
$\overset{\circ}{\geh}$.
\end{conjecture}
If $\geh$ is non-exceptional or $D^{(3)}_4$, i.e.,
$\geh = A^{(1)}_n, B^{(1)}_n, C^{(1)}_n, D^{(1)}_n, 
A^{(2)}_{2n}, A^{(2)}_{2n-1}, D^{(2)}_{n+1}$ or $D^{(3)}_4$,
the last property can indeed be proved at $q=1$.
It is based on the argument in Remark 8.8 of \cite{HKOTY1}
together with the existence of the Weyl group invariant solution
to the $Q$-system established in 
Theorem 7.1 in \cite{HKOTY1} and 
Theorem \ref{th:domino} in this paper.

Consider the nontwisted simply laced algebras 
$\geh=X^{(1)}_N = A^{(1)}_n, D^{(1)}_n\; (N=n)$ and $E^{(1)}_{6}$.
Let $\sigma$ be the order $r$ automorphism of the 
$X_N$ Dynkin diagram as 
in Section \ref{subsec:affine}.
We have $r=3$ if $\geh = D^{(1)}_4$ and $r=2$ otherwise.
We shall let $\sigma$ act on the index set $\{1,\ldots, N\}$ as well as 
on the weight lattice of $\overset{\circ}{\geh}= X_N$.
{}From the definition of the fermionic formula, it is 
straightforward to see
\begin{equation}\label{eq:symmetry}
M_l(\bigotimes_{(a,j) \in H_l} \left(W^{(a)}_j \right)^{\ot \nu^{(a)}_j},\la,q) %
= M_l(\bigotimes_{(a,j) \in H_l} 
\left(W^{(\sigma(a))}_j\right)^{\ot \nu^{(a)}_j},
\sigma(\la),q),
\end{equation}
which also implies (\ref{eq:small2}).

\medskip

Next we present a recursion relation of 
the fermionic form $M_l(W,\la,q)$ with respect to $W$.
For $a,b \in \{1,\ldots,n\}$ and $j,k \in \Z$, we set
(cf. eq.(A.6) in \cite{KN})
\begin{equation}\label{eq:Bdef}
\begin{gathered}
B_{bk,aj}=
2\min(t_ak,t_bj)-\min(t_ak,t_b(j+1))
-\min(t_ak,t_b(j-1)) \\
=\begin{cases}
2\delta_{k,2j}+\delta_{k,2j+1}+\delta_{k,2j-1}
& (t_b,t_a)=(2,1)\\
3\delta_{k,3j}+2\delta_{k,3j+1}+2\delta_{k,3j-1}&
 (t_b,t_a)=(3,1)\\
\qquad\qquad
+\delta_{k,3j+2}+\delta_{k,3j-2}&\\
t_b\delta_{t_ak,t_b j}
& \text{otherwise}.
\end{cases}
\end{gathered}
\end{equation}
\begin{proposition}\label{pr:recursion}
Fix $(a,j) \in \ol{H}_l$ arbitrarily.
Given any $\la \in \ol{P}^+$ and 
$W = \bigotimes_{(b,k) \in H_l}
(W^{(b)}_k)^{\otimes \nu^{(b)}_k}$, set
\begin{eqnarray*}
W_1 & = & W^{(a)}_{j} \otimes W^{(a)}_{j} \otimes W, \\
W_2 & = & W^{(a)}_{j + 1} 
\otimes W^{(a)}_{j - 1} \otimes W, \\
W_3 & = &
\bigotimes_{(b,k) \in H_l}
W^{(b)\,\otimes
(2\delta_{ab}\delta_{jk}-(\tilde{\alpha}_a\vert\tilde{\alpha}_b)B_{aj,bk}/t^\vee_b)}_k
\otimes W.
\end{eqnarray*}
Then we have
\begin{align*}
M_l(W_1,\lambda,q) & =  M_l(W_2,\lambda, q) 
+ q_{a}^{-\theta} M_l(W_3,\lambda,q), \\
\theta & = (2-\epsilon_a^{-1}) j 
+ \sum_{k = 1}^{t_al} \nu^{(a)}_k \min(j, k).
\end{align*}
The same recursion relation holds also for 
$\tilde{M}_l(W,\la,q)$.
\end{proposition}
When $j = 1$ we understand that 
$W_2 = W^{(a)}_2 \otimes W$.
The proof is similar to that for Theorem 6.1 in \cite{HKOTY1}
for $r=1$. 
In view of (\ref{eq:completeness}), $W_1,W_2$ and $W_3$ with trivial $W$ 
correspond to the three terms in the $Q$-system (\ref{eq:qsys}).

\subsection{\mathversion{bold} Limit of fermionic formulae}
\label{subsec:limit}

Recall that our fermionic formulae are connected to the 1dsums 
by Conjecture \ref{conj:X=M}, where the latters are 
indeed affine Lie algebra characters in the infinite lattice limit as in 
(\ref{eq:lim-Xinf}) and (\ref{eq:lim-Xl}).
Motivated by this fact, we here 
investigate the behavior of 
our fermionic formulae, either of argument $q$ or $q^{-1}$, 
when the quantum space data $\{ \nu^{(a)}_j \}$ gets large.
It turns out that they tend to an intriguing form\footnote{
Convergence in the sense of the $q$-adic topology.}.
Roughly, they break into a sum of products in which each factor
is again identified with a fermionic form
associated with ``smaller" levels or algebras.
It will be interesting to seek 
physical interpretations of such phenomena.
See \cite{ANOT,BLS2,NY2} and 
reference therein.

\begin{proposition}[Spinon character formula]\label{pr:spinon}
Assume that $l/t_b \in \Z$.
For any $1 \le b \le n$ we have
\begin{eqnarray}
\lim
q^{-c_0} M_\infty(W^{(b)\otimes L}_l,\lambda,q)
& = & \sum_{\zeta}\frac{M_\infty(W(\zeta),\lambda,q^{-1}) 
M_{l/t_b}(W(\zeta),0,q^{-1})}
{(q_1)_{\zeta_1} \cdots (q_n)_{\zeta_n}},\nonumber \\
W(\zeta) & = & \bigotimes_{a=1}^n W^{(a) \otimes \zeta_a}_1.\label{eq:Wzeta}
\end{eqnarray}
Here $\lim$ means the limit 
$L \rightarrow \infty$  so that 
$L(\iota(\ol{\La}_a) \vert \iota(\ol{\La}_b))/(t^\vee_at_b) 
\in {\mathbb Z}$ for 
$1 \le \forall a \le n$.
$c_0 = -L^2l\vert \iota(\ol{\La}_b) \vert^2/(2t_b) \in \Z$.
The sum is taken over $\zeta_1, \ldots, \zeta_n \in {\mathbb Z}_{\ge 0}$
such that $\sum_{a=1}^n \zeta_a \iota(\ol{\La}_a) \in 
\oplus_{a=1}^n \Z_{\ge 0} \tilde{\alpha}_a$.
\end{proposition}

Under Conjecture \ref{conj:fd-module} (1), 
the condition $l/t_b \in \Z$ is equivalent to assuming 
that the crystal $B^{b,l}$ is perfect.
Thus Proposition \ref{pr:spinon} covers 
all the twisted cases and  Theorem 5.4 in \cite{HKOTY1}, 
to which the proof is similar.
The quantity $c_0$ is the minimum  of the quadratic form $c(\{m\})$ 
in (\ref{eq:mc}) at $m^{(a)}_i = 
\delta_{i, t_al/t_b}L(\iota(\ol{\La}_a) \vert \iota(\ol{\La}_b))/(t^\vee_at_b)$.

For any $\geh$, one can extend  Proposition \ref{pr:spinon} to the situation 
where $W^{(b)}_l$ is replaced by a tensor product 
$W^{(k_1)}_{s_1} \ot \cdots \ot W^{(k_d)}_{s_d}$.
To present it,  we introduce a generalization of
(\ref{eq:mm2}) depending on a subset $S \subset \ol{H}_l$.
\begin{align}
&M^S_l(\bigotimes_{(a,j) \in H_l}
\bigl(W^{(a)}_j\bigr)^{\otimes \nu^{(a)}_j},q^{-1}) \nonumber\\
& = \sum_{\{m \}} q^{\ol{c}(\{m\})}
\prod_{(a,i) \in \ol{H}_{l}\setminus S}
\left[ \begin{array}{c} p^{(a)}_i +  m^{(a)}_i
 \\   m^{(a)}_i \end{array} \right]_{q_a}
\prod_{(a,i) \in S}\frac{1}{(q_a)_{p^{(a)}_i}},
\label{eq:MS}
\end{align}
where 
$\ol{c}(\{m\})$ and $p^{(a)}_i$ are given by the same formulae as
(\ref{eq:mc2}) and (\ref{eq:mp2}).
The sum $\sum_{\{ m \}}$ is now taken over
$\{ m^{(a)}_i \in {\mathbb Z}_{\ge 0} \mid (a,i) \in \ol{H}_{l}\setminus S \}$
and 
$\{ m^{(a)}_i \in {\mathbb Z} \mid (a,i) \in S \}$
satisfying 
$p^{(a)}_i \ge 0$ for all $(a,i) \in \ol{H}_{l}$ and 
the condition that $m^{(a)}_{t_al}$ determined {}from 
(\ref{eq:eliminate}) is a non-negative integer.
(\ref{eq:MS}) reduces to (\ref{eq:mm2}) when $S = \emptyset$.

Given the quantum space 
$(W^{(k_1)}_{s_1} \ot \cdots \ot W^{(k_d)}_{s_d})^{\ot L}$
and the associated fermionic form (\ref{eq:mm2}),
let $\{m_{cr} \} = \{m^{(a)}_{i,cr} (\in {\mathbb R}) \}$ be the 
critical point of the corresponding quadratic form (\ref{eq:mc2}).
We exclusively consider those 
$L$ such that $\forall m^{(a)}_{i,cr} \in {\mathbb Z}$.
With the help of  Remark 5.3\footnote{
The first line on the RHS therein should be corrected as 
$\{(a,\frac{t_as}{t_r}) \mid 1 \le a \le n \}$ 
for $\frac{s}{t_r} \in {\mathbb Z}$.} 
in \cite{HKOTY1} (see also eq.(3.5b-d) in \cite{Ku}), 
we find for any $\geh$ that 
\begin{align*}
\{(a,i) \mid m^{(a)}_{i,cr} > 0 \} & =  
S^{(k_1)}_{s_1} \cup \cdots \cup S^{(k_d)}_{s_d},\\
S^{(k)}_s &= \begin{cases}
\{ (a,\frac{t_as}{t_k}) \mid 1 \le a \le n \}\quad
\mbox{ for } \frac{s}{t_k} \in {\mathbb Z},\\[5pt]
\{ (a,\frac{t_a(s-s_0)}{t_k}),
 (a,\frac{t_a(s-s_0+t_k)}{t_k}) \mid 1 \le a \le n \}\\
\cup \{(a,s) \mid 1 \le a \le n, t_a = t_k(>1)\} \quad 
\mbox{ for } \frac{s}{t_k} \not\in {\mathbb Z},
\end{cases}
\end{align*}
where $s_0 \equiv s$ mod $t_k{\mathbb Z}$ and $1 \le s_0 \le t_k-1$.
It is easy to check that 
$\{(a,i) \mid m^{(a)}_{i,cr} > 0 \} \subset H_l$ 
for $l$ defined by 
$l = \max\{\lceil s_i/t_{k_i} \rceil \mid 1 \le i \le d \}$.

\begin{proposition}\label{pr:spinon2}
With $L$ and $l$ as above, we have the following for any $\geh$.
\begin{align}
&\lim_{L\rightarrow \infty} q^{-\ol{c}(\{m_{cr}\})}
M_\infty((W^{(k_1)}_{s_1} \ot \cdots \ot W^{(k_d)}_{s_d})^{\ot L},\la,q) \nonumber \\
&= \sum_\zeta \frac{M_\infty(W(\zeta),\la,q^{-1})
M^S_l(W(\zeta),q^{-1})}
{(q_1)_{\zeta_1} \cdots (q_n)_{\zeta_n}},\nonumber \\
&S = \{(a,i) \mid m^{(a)}_{i,cr} > 0 \} \setminus 
\{ (a,t_al) \mid 1 \le a \le n \}, \nonumber
\end{align}
where the sum $\sum_\zeta$ is over $\zeta_1, \ldots, \zeta_n \in 
{\mathbb Z}_{\ge 0}$ such that 
$\sum_{a=1}^n \zeta_a \iota(\ol{\La}_a) \in 
\bigoplus_{a=1}^n {\mathbb Z}_{\ge 0}{\tilde \alpha}_a$, and $W(\zeta)$ is
defined by (\ref{eq:Wzeta}).
\end{proposition}
Earlier Proposition \ref{pr:spinon} corresponds to the case 
$S = \emptyset$.
(Not $M^S_{l/t_b}$ here but $M^S_l$ as above.)
\begin{remark}\label{rem:rspinon}
When $\la=0$, Proposition \ref{pr:spinon2} admits a level truncation.
Let $l' (\ge l)$ be an integer. 
Under the same setting as Proposition \ref{pr:spinon2}, we have 
\begin{align*}
&\lim_{L\rightarrow \infty} q^{-\ol{c}(\{m_{cr}\})}
M_{l'}((W^{(k_1)}_{s_1} \ot \cdots \ot W^{(k_d)}_{s_d})^{\ot L},0,q) \\
&= \sum_\zeta \frac{M_{l'-l}(W(\zeta),0,q^{-1})
M^S_l(W(\zeta),q^{-1})}
{(q_1)_{\zeta_1} \cdots (q_n)_{\zeta_n}}.
\end{align*}
\end{remark}
In terms of appropriate summation variables\footnote{
Typically, one switches {}from the summation over 
$\{m^{(a)}_i\}$ to $\{p^{(a)}_i\}$.},
Remark \ref{rem:rspinon} reproduces the many expressions 
argued in \cite{KKMM2,BLS2,NY2,HKKOTY} for simply-laced $\geh$.
In the next two examples we let 
$(C_{i j})_{i,j=1}^{l-1}$ denote the Cartan matrix of $A_{l-1}$. 
\begin{example}
Take $\geh=A^{(1)}_1$ and $1 \le s \le l-1$. Then we have
\begin{align*}
&\lim_{\begin{subarray}{c} L \rightarrow \infty \\
L\in 2\Z  \end{subarray}} q^{L^2s/4}M_\infty(W^{(1) \ot L}_s,0,q) \\
& = 
\sum_{\{p\}} q^{\frac{1}{4}\sum_{i,j \in \ol{H}_l}C_{ij}p_ip_j}\frac{1}{(q)_{p_s}}
\prod_{\begin{subarray}{c} 1 \le i \le l-1  \\
i \neq s \end{subarray}}
\left[ \begin{array}{c} p_i -\frac{1}{2}\sum_{j=1}^{l-1}C_{ij}p_j
 \\   p_i \end{array} \right]_{q},
\end{align*}
where the sum $\sum_{\{p\}}$  extends over 
$p_1, \ldots, p_{l-1} \in 2\Z_{\ge 0}$ 
such that $\sum_{j=1}^{l-1}C_{ij}p_j \le 0$ for all $1 \le i \le l-1$ 
except $i=s$.
\end{example}
The RHS in the above agrees with the section 4.1 of \cite{KKMM2} for example.
\begin{example}
Take $\geh=A^{(2)}_2$ and $1 \le s \le l-1$. Then we have 
\begin{align*}
&\lim_{L \rightarrow \infty } q^{L^2s}M_\infty(W^{(1) \ot L}_s,0,q) \\
& = 
\sum_{\{p\}} q^{\sum_{i,j \in \ol{H}_l}C_{ij}p_ip_j}\frac{1}{(q^2)_{p_s}}
\prod_{\begin{subarray}{c} 1 \le i \le l-1  \\
i \neq s \end{subarray}}
\left[ \begin{array}{c} p_i - \sum_{j=1}^{l-1}C_{ij}p_j
 \\   p_i \end{array} \right]_{q^2},
\end{align*}
where the sum $\sum_{\{p\}}$  extends over 
$p_1, \ldots, p_{l-1} \in \Z_{\ge 0}$ such that 
$\sum_{j=1}^{l-1}C_{ij}p_j \le 0$ for all $1 \le i \le l-1$ 
except $i=s$.
\end{example}

\medskip

The spinon character formulae 
in Propositions \ref{pr:spinon}, \ref{pr:spinon2} and Remark \ref{rem:rspinon}
describe the limit of fermionic forms with argument $q$.
Now we turn to the other limit where $q$ is replaced by $q^{-1}$.
We first prepare some notations.
To save the space we set
\begin{equation}\label{eq:mcal}
{\mathcal M}^{(a)}_{q,\nu}(\geh \vert \mu) = 
M_\infty(\bigotimes_{j\ge 1}W^{(a)\; \otimes \nu_j}_j,
\sum_{j\ge 1}j\nu_j\ol{\La}_a-\mu,q).
\end{equation}
Here we have exhibited the dependence on the 
affine Lie algebra $\geh$ explicitly.
The index $1 \le a \le n$ specifies 
a vertex of the Dynkin diagram of $\overset{\circ}{\geh}$, 
$\nu =\{ \nu_k \}_{k\ge 1}$ is a sequence of non-negative integers and 
$\mu$ is an element of 
$\oplus_{c=1}^n \Z_{\ge 0}{\alpha}_c/\epsilon_c$.
In the sequel we let $D(X)$ denote the Dynkin diagram of $X$.
Let $\geh$ be an affine Lie algebra of rank $n$ and 
$\{t^\vee_i \mid i \in I\}$ be the associated data 
(\ref{eq:ttdef}).
Fix an index $1 \le b \le n$, and 
remove the two vertices $0$  and $b$ and all the edges connected to them 
{}from the Dynkin diagram $D(\geh)$.
The result is a disjoint union of connected diagrams $D_a$
labeled by the indices $1 \le a \le n$ 
satisfying $({\alpha}_a \vert {\alpha}_b) < 0$.
To each of them we assign an affine Lie subalgebra 
$\geh^\ast \subset \geh$ uniquely determined by the following condition:
\begin{equation}\label{eq:fool}
\begin{gathered}
\geh^\ast = A^{(\epsilon_a)}_{\epsilon_an'} \quad \text{if }
D_a = D(A_{n'}) \text{ for some } n',\\
(D(\overset{\circ}{\geh}{}^\ast), 
\{t^\vee_c \, (c \neq 0) \text{ of } \geh^\ast \})
=
(D_a, \{t^\vee_c \text{ of } 
\geh \mid c \text{ belonging to } D_a \})\quad \text{otherwise}.
\end{gathered}
\end{equation}
See Table \ref{tab:geh-bar} for the definition of the 
classical subalgebra $\overset{\circ}{\geh}{}^\ast \subset \geh^\ast$.
The symbol $\epsilon_a$ in  (\ref{eq:epsdef}) actually becomes 2
only when $\geh=A^{(2)}_{2n}, (a,b)=(n,n-1)$ hence $n'=1$ in the above.
When regarding $D(\overset{\circ}{\geh}{}^\ast)$ as a subdiagram of 
$D(\geh)$, the labeling of vertices within the former along Table \ref{tab:geh-bar}
can become different {}from the restriction of the latter.
So we let $a^\ast (\neq 0)$ denote the label within $\geh^\ast$ 
corresponding to the $a$ for $\geh$.
Let further $\mu^\ast$ be the projection of $\mu$ 
to the lattice $\oplus_{c=1}^n \Z_{\ge 0}{\alpha}_c/\epsilon_c$ 
for $\overset{\circ}{\geh}{}^\ast$ according to 
$\overset{\circ}{\geh}{}^\ast \subset \overset{\circ}{\geh}$\footnote{
Depending on the labelings in the two algebras, 
one can have 
$(\alpha_i\text{ for } \overset{\circ}{\geh})^\ast = 
\alpha_{i'} \text{ for }\overset{\circ}{\geh}{}^\ast$ with $i \neq i'$ 
in general. Concrete formulae are available in 
Appendix \ref{app:regII} for each $\geh$.}.

\begin{proposition}\label{pr:dualspinon}
Fix $1 \le b \le n,\; l \in \Z_{\ge 1}$ and 
$\mu = \mu_1 {\alpha}_1 + \cd + \mu_n{\alpha}_n$. 
With the above definitions of $\geh, \, a^\ast$ and $\mu^\ast$, we have 
\begin{align}
&\lim_{L \rightarrow \infty}
M_\infty(W^{(b) \ot L}_l,Ll\ol{\La}_b-\mu,q^{-1})\nonumber \\
&= \sum_{\nu} \frac{q^{t^\vee_b \psi(\nu)/\epsilon_b}}
{\prod_{i\ge 1}(q_b)_{\nu_i}}
\prod_{\begin{subarray}{c} 1 \le a \le n  \\
({\alpha}_a\vert {\alpha}_b)<0 \end{subarray}}
{\mathcal M}^{(a^\ast)}_{q^\ast,\nu^\ast}(\geh^\ast\vert  \mu^\ast),
\label{eq:mpro} \\
&\psi(\nu) = \sum_{i,j \ge 1} \min(i,j) \nu_i\nu_j, \label{eq:nuphi}\\
&\nu^\ast = \{\nu^\ast_j \}_{j \ge 1}, \quad 
\nu^\ast_j = -\frac{(\tilde{\alpha}_a\vert \tilde{\alpha}_b)}{t^\vee_a}
\sum_{k\ge 1} \nu_k B_{b  k, a j} \; \; (\in \Z_{\ge 0}),
\label{eq:nuast}\\
&q^\ast = \begin{cases}
q_a & \text{ if } \geh^\ast = A^{(1)}_{n'} \text{ for some } n',\\
q & \text{otherwise},
\end{cases} \label{eq:qstar}
\end{align}
where the sum in (\ref{eq:mpro}) 
is taken over $\nu_1, \nu_2, \ldots \in \Z_{\ge 0}$ 
such that 
\begin{equation}\label{eq:nuhani}
\sum_{i \ge 1}i\nu_i = \epsilon_b\mu_b.
\end{equation}
\end{proposition}
Explicit forms of (\ref{eq:mpro}) are available in Appendix \ref{app:regII}.
Note that the results are totally independent of $l$.
\begin{proof}
We start {}from the expression (\ref{eq:mm3}).
In the limit $L \rightarrow \infty$, one has 
$p^{(b)}_i \rightarrow \infty$ for all $i$. 
Thus by setting $\nu_i = m^{(b)}_i$, the factor 
$\prod_i\left[ \begin{array}{c} p^{(b)}_i +  m^{(b)}_i
 \\   m^{(b)}_i \end{array} \right]_{q_b}$ gives rise to 
$1/\prod_i(q_b)_{\nu_i}$.
The summand 
$\frac{1}{2}(\tilde{\alpha}_b \vert \tilde{\alpha}_b)
\sum_{j,k}\min(t_bj,t_bk)m^{(b)}_jm^{(b)}_k$ contained 
in the quadratic form (\ref{eq:mc3}) thereby becomes 
$t^\vee_b \psi(\nu)/\epsilon_b$ on account of (\ref{eq:account}).
For an $a$ such that $(\alpha_a \vert \alpha_b) < 0$ 
(hence $a \neq b$ and $(\tilde{\alpha}_a \vert \tilde{\alpha}_b) < 0$),
write $p^{(a)}_i$ in (\ref{eq:mp3}) as
\begin{equation*}
p^{(a)}_i = 
- \frac{(\tilde{\alpha}_a \vert \tilde{\alpha}_b)}{t^\vee_a}
\sum_{k\ge 1}  
 \mbox{min}(t_bi,t_ak) \nu_k 
- \frac{1}{t^\vee_a}\sum_{c \neq b}\sum_{k\ge 1}  
(\tilde{\alpha}_a \vert \tilde{\alpha}_c) \mbox{min}(t_ci,t_ak) m^{(c)}_k,
\end{equation*}
where the $c$-sum actually extends over only the vertices of $D_a$.
To read off the quantum space data $\nu^\ast$ for  $\geh^\ast$ 
{}from this, one identifies the first terms here with  
that for $p^{(a^\ast)}_i$ {}from (\ref{eq:mp3}),
leading to 
\begin{equation*}
- \frac{(\tilde{\alpha}_a \vert \tilde{\alpha}_b)}{t^\vee_a}
\sum_{k\ge 1}  
 \mbox{min}(t_ak,t_bi) \nu_k 
 = \sum_{j\ge 1} \min(i,j)\nu^\ast_j.
\end{equation*}
Solving this one finds (\ref{eq:nuast}).
By a direct case check, 
the other  $q$-binomial factors and the remaining 
terms in the quadratic form (\ref{eq:mc3}) 
can be matched with the contributions {}from the product of 
${\mathcal M}^{(a^\ast)}_{q^\ast, \nu^\ast}$'s.
\end{proof}
\medskip
%

\section{\mathversion{bold} $q$-series $N_l(\la,q)$}
\label{sec:N}
We keep the convention 
$\tilde{\alpha}_a \in \tilde{P}$ and $(.\vert.)$ 
explained in Section \ref{subsec:convention}.

For $l \in \Z_{\ge 1}$ and $\la \in \ol{P}$,
we introduce a $q$-series $N_l(\la,q)$ by
\begin{eqnarray}
N_l(\lambda,q) & = & 
\frac{1}{\prod_{a=1}^n(q_a)_\infty}
\sum_{\{m \}} \frac{q^{\phi(\{m\})}}
{\prod_{(a,j) \in \ol{H}_l}(q_a)_{m^{(a)}_j}}, \label{eq:Ndef}\\
\phi(\{m\}) & = & \frac{1}{2}
\sum_{(a,j),(b,k) \in \ol{H}_l} (\tilde{\alpha}_a \vert \tilde{\alpha}_b) 
K^{(t_at_bl)}_{t_bj, t_ak}m^{(a)}_jm^{(b)}_k,  \label{eq:bekiphi}
\end{eqnarray}
where $q_a$ and $K^{(l)}_{j,k}$ are defined in 
(\ref{eq:mqt}) and (\ref{eq:kdef}), respectively.
The sum $\sum_{\{m\}}$ is taken over 
$\{m^{(a)}_j \in \Z_{\ge 0} \mid (a,j) \in \ol{H}_l \}$
satisfying 
\begin{align*}
&\sum_{(a,j) \in \ol{H}_l} jm^{(a)}_j \tilde{\alpha}_a \equiv \iota(\la)
\quad \mod lM,\\
&M = \oplus_{a=1}^n \Z t_a \tilde{\alpha}_a.
\end{align*}
The sublattice $M \subset  \tilde{P}$
is essentially the one introduced in the section 6.5 of 
\cite{Kac}\footnote{This $M$ should not be confused with 
a fermionic form.}.
{}From the definition it follows that 
\begin{align}
N_l(\la,q) &= N_l(\la + l\beta,q) \quad \mbox{for } \beta \in M,
\label{eq:Mshift}\\
N_l(\la,q) & \in q^{\epsilon'}\Z_{\ge 0}[[q^{a_0}]], \quad 
\epsilon' \equiv -\frac{\vert\iota(\la)\vert^2}{2l} + \sum_{a=1}^na\la_a \; 
\mod a_0\Z.\nonumber
\end{align}
The $q$-series $N_l(\la,q)$ is related to the 
fermionic formula in the previous section via
\begin{equation}\label{eq:MN}
N_l(\la,q) = 
\lim q^{-\frac{\vert \iota(\la) \vert^2}{2l}}
M_l(W,\sum_{(a,j) \in H_l}j\nu^{(a)}_j \ol{\La}_a - \la,q^{-1}),
\end{equation}
where $M_l$ can also be replaced with $\tilde{M}_l$.
Here the quantum space data is given by 
$W = \bigotimes_{(a,j) \in H_l}
\bigl(W^{(a)}_j\bigr)^{\otimes \nu^{(a)}_j}$, and the limit is 
taken with respect to $\{\nu^{(a)}_j\}$ so that 
$\sum_{j \ge 1}\nu^{(a)}_j \rightarrow \infty$ 
for all $1\le a \le  n$.

Let us present the examples {}from $\geh$ of rank 1.

\begin{example}
Take $\geh=A^{(1)}_1$, for which the formula reads
\begin{equation*}
N_l(\la_1\ol{\La}_1,q) = 
\frac{1}{(q)_\infty}
\sum_{\{m \}} \frac{q^{\sum_{1 \le i,j \le l-1}K^{(l)}_{i,j}m_im_j}}
{\prod_{j=1}^{l-1}(q)_{m_j}},
\end{equation*}
where $\la_1 \in 2\Z$ and 
the sum is over $m_1, \ldots, m_{l-1} \in \Z_{\ge 0}$ such that 
$\sum_{j=1}^{l-1}jm_j \equiv \lambda_1/2 \mod l\Z$.
See (\ref{eq:kdef}) for $K^{(l)}_{i,j}$.
We have
\begin{align*}
&N_2(\la_1\ol{\La}_1,q) \\
&= 1+q+3q^2+5q^3+10q^4+16q^5+28q^6+43q^7+70q^8+\cdots
\quad \mbox{if } \la \in 4\Z,\\
&=q^{\frac{1}{2}}
(1+2q+4q^2+7q^3+13q^4+21q^5+35q^6+55q^7+86q^8+\cdots)
\quad \mbox{if } \la_1 \not\in 4\Z.\\
&N_3(\la_1\ol{\La}_1,q) \\
&= 1+q+3q^2+6q^3+12q^4+21q^5+39q^6+64q^7+108q^8+\cdots
\quad \mbox{if } \la_1 \in 6\Z,\\
&=q^{\frac{2}{3}}
(1+2q+5q^2+9q^3+18q^4+31q^5+55q^6+90q^7+149q^8+\cdots)
\quad \mbox{if } \la_1 \not\in 6\Z.
\end{align*}
\end{example}

\begin{example}
Take $\geh=A^{(2)}_2$, for which the formula reads
\begin{equation*}
N_l(\la_1\ol{\La}_1,q) = 
\frac{1}{(q^2)_\infty}
\sum_{\{m \}} \frac{q^{\sum_{1 \le i,j \le l-1}K^{(l)}_{i,j}m_im_j}}
{\prod_{j=1}^{l-1}(q^2)_{m_j}},
\end{equation*}
where $\la_1 \in \Z$ and 
the sum is over $m_1, \ldots, m_{l-1} \in \Z_{\ge 0}$ such that 
$\sum_{j=1}^{l-1}jm_j \equiv \lambda_1 \mod l\Z$.
Besides the summation condition, the RHS with the replacement 
$q \rightarrow q^{1/2}$  differs from 
$A^{(1)}_1$ case only by the factor $1/2$ in the quadratic form.
We have
\begin{align*}
&N_2(\la_1\ol{\La}_1,q) \\
&= 1+2q^2+4q^4+8q^6+15q^8+26q^{10}+44q^{12}+72q^{14}+115q^{16}+\cdots
\quad \mbox{if } \la_1 \in 2\Z,\\
&=q^{\frac{1}{2}}
(1+2q^2+5q^4+9q^6+17q^8+29q^{10}+50q^{12}+80q^{14}+129q^{16}+\cdots)
\quad \mbox{if } \la_1 \not\in 2\Z.\\
&N_3(\la_1\ol{\La}_1,q) \\
&= 1+2q^2+5q^4+12q^6+24q^8+46q^{10}+85q^{12}+150q^{14}+\cdots
\quad \mbox{if } \la_1 \in 3\Z,\\
&=q^{\frac{2}{3}}
(1+3q^2+7q^4+15q^6+30q^8+57q^{10}+103q^{12}+180q^{14}+\cdots)
\quad \mbox{if } \la_1 \not\in 3\Z.
\end{align*}
\end{example}

Our $q$-series $N_l(\la,q)$ is connected to representation theory via

\begin{conjecture}\label{con:c=N}
For $l \in \Z_{\ge 1},\, \la \in \ol{P}$ and 
some $\kappa \in \Z-\epsilon'/a_0$, one has
\[
\sum_{i \in \Z} \mbox{mult}_{L(l\La_0)}
(\la+l\La_0-i\delta) q^i = q^\kappa N_l(\la,q^{\frac{1}{a_0}}).
\]
\end{conjecture}

The LHS is a string function of the level $l$ vacuum module over 
$X^{(r)}_N$ up to a power of $q$.
The conjecture is valid for  
$A^{(1)}_1$ \cite{LP}, $A^{(1)}_n$ \cite{FS,G,HKKOTY},
and also for $l=1$ \cite{Kac}.
For general nontwisted $\geh$,  
the $q$-series $N_l(\la,q)$ 
appeared in \cite{KNS}.
See also the section 6.1 of \cite{HKKOTY} for a generalization of 
$N_l(\la,q)$ corresponding to the tensor product 
$\otimes_iL((l_i-l_{i+1})\La_0)$  for nontwisted cases.
It is natural to expect a similar generalization for 
twisted cases, but  we leave it to a future study.

Under  Conjecture \ref{con:c=N},
the property (\ref{eq:Mshift}) is consistent with 
eq.(12.7.9) in \cite{Kac}.
The asymptotic behavior of $N(\la,q^{1/a_0})$ as $q \rightarrow 1$ 
also matches the known result. 
To explain this let  
\[
L(x) = -\frac{1}{2}\int_0^x\left( \frac{\log(1-y)}{y} + 
\frac{\log y}{1-y} \right) dy
\]
be the Rogers dilogarithm function.
To the power series expansion \\
$q^{-\epsilon'/a_0}\left(\prod_{b=1}^n(q^{1/a_0}_b)_\infty\right)
N_l(\la,q^{1/a_0}) = \sum_k c_k q^k$, 
one can apply the saddle-point method of \cite{RS} to derive a crude estimate 
$\log c_k \sim 2\sqrt{ka_0\sum_{(a,j) \in \ol{H}_l}L(f^{(a)}_j)/t^\vee_a}
+ O(\log k)$ for $k \gg 1$.
Here $f^{(a)}_j$ is the one in  Conjecture \ref{con:qsys} (iii), 
by which the ``saddle point condition" has been solved as 
$q^{t^\vee_am^{(a)}_j/a_0} = 1 - f^{(a)}_j$  with the help of (\ref{eq:rqsys2}).
Combining the result with the well-known property 
$\lim_{q \rightarrow 1}(1-q)\log (q^{1/a_0}_a)_\infty^{-1}
=  \pi^2a_0/(6t^\vee_a)$, one arrives at
\begin{equation}\label{eq:Ndiv}
\lim_{q \rightarrow 1} (1-q) \log N_l(\la,q^{1/a_0}) = 
a_0\left( \frac{\pi^2}{6} \sum_{a=1}^n\frac{1}{t^\vee_a} + 
\sum_{(a,j) \in \ol{H}_l} \frac{L(f^{(a)}_j)}{t^\vee_a} \right).
\end{equation}
Applying (\ref{eq:wa}) and admitting  (\ref{eq:dilog}),
one finally  obtains
\begin{equation}\label{eq:div}
\lim_{q \rightarrow 1} (1-q) \log N_l(\la,q^{1/a_0}) = 
\frac{\pi^2}{6r}\left(
\frac{ l \dim X_N }{l + h^\vee}\right).
\end{equation}
This agrees with the behavior of the level $l$ string function 
as $q= \exp(-2\pi \beta) \rightarrow 1$ given by
 Theorem 13.14. a) and eq.(13.13.12) in \cite{Kac}.

\section{\mathversion{bold} $Q$-system and 
fermionic formula at $q=1$}\label{sec:qsys}

\subsection{\mathversion{bold} $Q$-system of type $\geh$}
\label{subsec:qsys}

Let $\{Q^{(a)}_j \mid 1 \le a \le n, j \in \Z_{\ge 0} \}$ 
be a set of infinitely many commutative variables
with $Q^{(a)}_0 = 1$.
We  call the simultaneous equation 
\begin{equation}\label{eq:qsys}
Q^{(a)\, 2}_{j} = Q^{(a)}_{j + 1} Q^{(a)}_{j - 1}
+  Q^{(a)\, 2}_{j}
\prod_{b=1}^n\prod_{k \ge 1}
Q^{(b)\, -(\tilde{\alpha}_a\vert\tilde{\alpha}_b)B_{aj,bk}/t^\vee_b}_k
\end{equation}
the (unrestricted) $Q$-system of type $\geh$.
Here $B_{aj,bk}$ is defined in (\ref{eq:Bdef}) and 
the second term on RHS actually contains finitely many factors 
with positive power.
For $\geh = X^{(1)}_n$, it appeared
in \cite{KR2,Ki3} as a possible character identity 
among the irreducible finite-dimensional modules
of the Yangian $Y(X_n)$ viewed as $X_n$ modules.
The explicit form  for the nontwisted cases is 
available in Appendix C in \cite{HKOTY1}.
Let us write down the twisted cases $r =2$ and $3$. 
Below we assume that $Q^{(a)}_j = 1$ whenever 
$a \not \in \{1, \ldots, n\}$. \\ 
\noindent 
${\mathfrak g}=A_{2n-1}^{(2)}:$
\begin{eqnarray*}
Q^{(a)^2}_j
  &=&Q_{j-1}^{(a)}Q_{j+1}^{(a)} + Q^{(a-1)}_j Q^{(a+1)}_j
\quad \text{for} \quad 1 \le a \le n-1 , \\
Q^{(n)^2}_{j}
&=& Q_{j-1}^{(n)} Q_{j+1}^{(n)}
+ Q^{(n-1)^2}_j. 
\end{eqnarray*}
\noindent
${\mathfrak g}=A_{2n}^{(2)}:$
\begin{eqnarray*}
Q^{(a)^2}_j
  &=&Q_{j-1}^{(a)}Q_{j+1}^{(a)} + Q^{(a-1)}_j Q^{(a+1)}_j
\quad \text{for} \quad 1 \le a \le n-1 , \\
Q^{(n)^2}_{j}
&=& Q_{j-1}^{(n)} Q_{j+1}^{(n)}
+ Q^{(n-1)}_jQ^{(n)}_j. 
\end{eqnarray*}
\noindent
${\mathfrak g}=D_{n+1}^{(2)}:$
\begin{eqnarray*}
&& Q^{(a)^2}_j
  =Q_{j-1}^{(a)}Q_{j+1}^{(a)} + Q^{(a-1)}_j Q^{(a+1)}_j
\quad \text{for} \quad 1 \le a \le n-2, \\
&& Q^{(n-1)^2}_{j}
= Q_{j-1}^{(n-1)} Q_{j+1}^{(n-1)}
+ Q^{(n-2)}_jQ^{(n)^{2}}_j, 
\\ 
&&Q^{(n)^2}_{j}
= Q_{j-1}^{(n)} Q_{j+1}^{(n)}
+ Q^{(n-1)}_j. 
\end{eqnarray*}
\noindent
${\mathfrak g}=E_{6}^{(2)}:$
\begin{eqnarray*}
Q^{(1)^2}_j
  &=&Q_{j-1}^{(1)}Q_{j+1}^{(1)} + Q^{(2)}_j, \\
Q^{(2)^2}_{j}
&=& Q_{j-1}^{(2)} Q_{j+1}^{(2)}
+ Q^{(1)}_jQ^{(3)}_j, 
\\ 
Q^{(3)^2}_{j}
&=& Q_{j-1}^{(3)} Q_{j+1}^{(3)}
+ Q^{(2)\,2}_jQ^{(4)}_j, \\ 
Q^{(4)^2}_{j}
&=& Q_{j-1}^{(4)} Q_{j+1}^{(4)}
+ Q^{(3)}_j.
\end{eqnarray*}
\noindent 
${\mathfrak g}=D_{4}^{(3)}:$
\begin{eqnarray*}
Q^{(1)^2}_j
  &=&Q_{j-1}^{(1)}Q_{j+1}^{(1)} + Q^{(2)}_j, \\
Q^{(2)^2}_{j}
&=& Q_{j-1}^{(2)} Q_{j+1}^{(2)}
+ Q^{(1)\,3}_j.
\end{eqnarray*}
Let us comment on an origin of the above proposal 
for twisted cases.
First, it is a reduction of the $T$-system in \cite{KS}
by forgetting the spectral parameter. 
Second, $X^{(r)}_N$ $Q$-system with $r > 1$ is obtained {}from the 
$X^{(1)}_N$ $Q$-system by imposing a symmetry 
under the order $r$ Dynkin diagram automorphism $\sigma$.
To explain this, we introduce the map 
$\tau$\footnote{The $\tau$ here should be distinguished {}from the one 
in (\ref{eq:auto-B}).} between the non 0 vertices of the Dynkin diagrams 
$\tau: \{1,2,\ldots,N\}$ for $X^{(1)}_N  \rightarrow 
\{1,\ldots, n\}$ for $X^{(r)}_N$ by
\begin{align*}
\tau(a) &= \min(a,N+1-a) \quad \mbox{ for } A^{(2)}_{2n}, A^{(2)}_{2n-1},\\
\tau(a) &= \min(a,n) \quad \mbox{ for } D^{(2)}_{n+1},\\
\tau(a)&= \min(a,6-a) \; (1 \le a \le 5), \;\;\tau(6) = 4
\quad \mbox{ for } E^{(2)}_{6}, \\
\tau(a) &= 1 \;(a=1,3,4), \;\;\tau(2) = 2 
\quad \mbox{ for } D^{(3)}_{4}.
\end{align*}
The map $\tau$ respects the Dynkin diagram symmetry  
$\sigma$ on $\{1,\ldots, N\}$ in the sense that 
$\tau \circ \sigma = \tau$.
Our $X^{(r)}_N$ type $Q$-system is obtainable from the 
$X^{(1)}_N$ type $Q$-system through the replacement
$Q^{(a)}_j \rightarrow Q^{(\tau(a))}_j$.

For a given $l \in \Z_{\ge 1}$, impose the
condition $Q^{(a)}_{t_al}= 1, Q^{(a)}_{t_al+1}=0$.
Then one can consider the truncated system of equations 
\begin{equation}\label{eq:rqsys}
Q^{(a)\, 2}_{j} = Q^{(a)}_{j + 1} Q^{(a)}_{j - 1}
+  Q^{(a)\, 2}_{j}
\prod_{(b,k) \in \ol{H}_l}
Q^{(b)\, -(\tilde{\alpha}_a\vert\tilde{\alpha}_b)B_{aj,bk}/t^\vee_b}_k\qquad
(a,j) \in \ol{H}_l,
\end{equation}
which close among $\{Q^{(a)}_j \mid (a,j) \in \ol{H}_l\}$.
We call it the level $l$ restricted $Q$-system of type $\geh$.
In terms of the 
variables $f^{(a)}_j = 1-Q^{(a)}_{j+1}Q^{(a)}_{j+1}/Q^{(a) 2}_j$,
it is expressed as 
\begin{equation}\label{eq:rqsys2}
\log f^{(a)}_j = \sum_{(b,k) \in \ol{H}_l}
\frac{(\tilde{\alpha}_a \vert \tilde{\alpha}_b)}{t^\vee_b}
K^{(t_at_bl)}_{t_bj, t_ak}\log(1-f^{(b)}_k).
\end{equation}

For any $\geh = X^{(r)}_N$ we have
\begin{conjecture}\label{con:qsys}

\noindent 
(i) The unrestricted $Q$-system (\ref{eq:qsys}) has a solution 
in terms of a linear 
combination of $\overset{\circ}{\geh}$-characters of the form
\[
Q^{(a)}_j = \sum_{\la \in 
\ol{P} \cap (j\ol{\La}_a -\sum_{b=1}^n \Z_{\ge 0}\alpha_b/\epsilon_b)}
d_\la \mbox{ch }V(\la),
\quad d_{j\ol{\La}_a} = 1, \;\; d_\la \in \Z_{\ge 0}.
\]

\noindent
(ii) Specialization of (i) 
to the $q$-dimension at $q = \exp(2\pi i/(l+h^\vee))$ has the properties 
$Q^{(a)}_{t_al}= 1, Q^{(a)}_{t_al+1}=0$ and 
$Q^{(a)}_j \in {\mathbb R}_{>0}$ for all $(a,j) \in \ol{H}_l$.
(Hence it yields a solution of the level $l$ restricted $Q$-system 
(\ref{eq:rqsys}).)

\noindent
(iii) In terms of $Q^{(a)}_j \in {\mathbb R}_{>0}$ in (ii), set 
$f^{(a)}_j = 1-Q^{(a)}_{j+1}Q^{(a)}_{j+1}/Q^{(a) 2}_j$.
(Hence it is a solution of (\ref{eq:rqsys2}) such that 
$0 < f^{(a)}_j < 1$ for all 
$(a,j) \in \ol{H}_l$.)
Then we have 
\begin{equation}\label{eq:dilog}
a_0r \sum_{(a,j) \in \ol{H}_l}
\frac{ L(f^{(a)}_j)}{t^\vee_a} = \frac{\pi^2}{6}\left(
\frac{ l \dim X_N }{l + h^\vee}-N\right).
\end{equation}
\end{conjecture}

Several remarks are in order.
First, (i) is indeed valid if $\geh$ is non-exceptional or $D^{(3)}_4$.
The explicit forms are available in \cite{KR2, HKOTY1} 
for the nontwisted cases $\geh = A^{(1)}_n, B^{(1)}_n, C^{(1)}_n, D^{(1)}_n$
and in Section \ref{subsec:qsol} of this paper 
for $\geh = A^{(2)}_{2n}, A^{(2)}_{2n-1}, D^{(2)}_{n+1}$ and $D^{(3)}_4$.
Thus (i) is actually an existence  conjecture for the other
$\geh$'s.
On the other hand for {\em any} $\geh$, 
the forthcoming Theorem \ref{th:completeness} 
asserts the uniqueness as well as a 
unified construction 
$Q^{(a)}_j = \sum_\la \tilde{M}(W^{(a)}_j,\la,1) \mbox{ch} V(\la)$ 
in terms of the fermionic form at $q=1$ 
under a further assumption (C) on the asymptotic behavior.
For the exceptional $\geh$'s other than $D^{(3)}_4$, 
it is yet to be confirmed if 
the above $Q^{(a)}_j$ indeed satisfies the $Q$-system.
Second, if Conjecture \ref{con:c=N} and \ref{con:qsys} (i)--(ii) are valid, 
then the dilogarithm sum rule 
(\ref{eq:dilog}) follows  {}from 
(\ref{eq:Ndiv}) and  Theorem 13.14. a) in \cite{Kac} by 
reversing the argument in the end of Section \ref{sec:N}.
In this sense  Conjecture \ref{con:c=N} is a $q$-analogue of 
(\ref{eq:dilog}).
Third, (\ref{eq:dilog}) with $r=1$ was originally conjectured 
in \cite{Ki3}, where $A^{(1)}_n$ case was proved.
(For $X^{(1)}_n$ non-simply laced, 
the range of the sum therein is misprinted and should be read as above.)
Our (\ref{eq:dilog}) for $X^{(r)}_N$ with $r>1$ 
is actually a consequence of the 
$X^{(1)}_N$ case.
To see this, note that 
$f^{\prime (a)}_j$ for $X^{(1)}_N$ enjoys the symmetry 
$f^{\prime (a)}_j = f^{\prime (\sigma(a))}_j$ under the 
order $r$ Dynkin diagram automorphism $\sigma$.
Therefore $f^{(a)}_j$ for $X^{(r)}_N$ is obtained as 
$f^{(a)}_j = f^{\prime (b)}_j$ for any $b$ such that $\tau(b) = a$.
In view of this, LHS of (\ref{eq:dilog}) 
is equal to $\sum_{a=1}^N\sum_{j=1}^{l-1}L(f^{\prime (a)}_j)$, 
reducing it to the $X^{(1)}_N$ case since 
$h^\vee$ is independent of $r$.
In particular, the factor $a_0r/t^\vee_a \in \{1, r \}$ 
has the interpretation
\begin{equation}\label{eq:hehe}
\frac{a_0r}{t^\vee_a} = \sharp\{b \in \{1,\ldots, N\} \mid \tau(b) = a \}
\quad \mbox{for } 1 \le a \le n.
\end{equation}
In view of (\ref{eq:wa}) and 
$f^{(a)}_{t_al} = 1$ under $Q^{(a)}_{t_al+1}=0$,
the sum rule (\ref{eq:dilog}) is also stated as 
\begin{equation*}
a_0r \sum_{(a,j) \in {H}_l}
\frac{ L(f^{(a)}_j)}{t^\vee_a} = \frac{\pi^2}{6}\left(
\frac{ l \dim X_N }{l + h^\vee}\right).
\end{equation*}

\subsection{\mathversion{bold} Solutions in terms of 
classical characters}\label{subsec:qsol}

Let us present solutions of the unrestricted $Q$-system 
(\ref{eq:qsys}) for 
$\geh = A^{(2)}_{2n}, A^{(2)}_{2n-1}$, $D^{(2)}_{n+1}, D^{(3)}_4$
in terms of the character of a classical subalgebra 
$\geh' \subset \geh$.
Similar result is available for 
$\geh = A^{(1)}_n, B^{(1)}_n, C^{(1)}_n, D^{(1)}_n$
\cite{KR2,HKOTY1} for the choice $\geh' = \overset{\circ}{\geh}$.
Here we shall include the result 
$(\geh, \geh') = (A^{(2)}_{2n},C_n), 
(A^{(2)}_{2n},B_n), (A^{(2)}_{2n-1},C_n)$, 
$(A^{(2)}_{2n-1},D_n), (D^{(2)}_{n+1},B_n)$ and 
$(D^{(3)}_4,G_2)$.
In case $\geh' = \overset{\circ}{\geh}$,
Appendix \ref{app:mlist} with $q=1$ yields
further conjectures on $E^{(2)}_6$.

Let us introduce the linear combination of characters. \\
\noindent $({\mathfrak g},\geh')=(A_{2n}^{(2)},C_{n})$:
\begin{eqnarray}
& \chi^{(a)}_j &= 
\sum \mbox{ch}V(k_1\ol{\La}'_1+ k_2\ol{\La}'_2 + \cdots + k_a\ol{\La}'_a) 
\quad \text{for} \quad 1 \le a \le n.
\label{eq:domino-aec}
\end{eqnarray}
\noindent $({\mathfrak g},\geh')=(A_{2n}^{(2)},B_{n})$:
\begin{eqnarray}
\chi^{(a)}_j = 
\begin{cases}
\sum \mbox{ch}V(k_1\ol{\La}'_1+ k_2\ol{\La}'_2 + \cdots + k_a\ol{\La}'_a) 
 & \text{if} \quad 1 \le a \le n-1, \\ 
\sum \mbox{ch}V(k_1\ol{\La}'_1+ k_2\ol{\La}'_2 + 
 \cdots +k_{n-1}\ol{\La}'_{n-1} \\ 
 \hspace{100pt} +2k_n\ol{\La}'_n) 
 & \text{if} \quad a=n.
\end{cases}
\label{eq:domino-aeb} 
\end{eqnarray}
\noindent $({\mathfrak g},\geh')
=(A_{2n-1}^{(2)},C_{n})$:
\begin{eqnarray}
&& \chi^{(a)}_j = \sum 
\mbox{ch}V(k_{b}\ol{\La}'_{b} + k_{b + 2}\ol{\La}'_{b + 2} 
+ \cdots +k_a\ol{\La}'_a)   \label{eq:domino-aoc} \\ 
&& \hspace{100pt} 
\text{for} \quad 
1 \le a \le n, \quad 
b \equiv a \; \; \mbox{mod 2}, \quad 
b=0 \; \; \text{or} \; \; 1. \nonumber 
\end{eqnarray}
\noindent $({\mathfrak g},\geh')=(A_{2n-1}^{(2)},D_{n})$: 
\begin{eqnarray}
\chi^{(a)}_j =
\begin{cases}
\sum \mbox{ch}V(k_{1}\ol{\La}'_{1} + k_{2}\ol{\La}'_{ 2} 
+ \cdots +k_a\ol{\La}'_a)  & \text{if} \quad 1 \le a \le n-2,  \\ 
\sum \mbox{ch}V(k_{1}\ol{\La}'_{1} + k_{2}\ol{\La}'_{ 2} 
+ \cdots +k_{n-2}\ol{\La}'_{n-2} \\ 
 \hspace{90pt} +k_{n-1}\ol{\omega}_{n-1})
 & \text{if} \quad a=n-1,  \\ 
\sum \mbox{ch}V(k_{1}\ol{\La}'_{1} + k_{2}\ol{\La}'_{ 2} 
+ \cdots +k_{n-2}\ol{\La}'_{n-2} \\ 
 \hspace{40pt} 
+(k_{n-1}+|k_{n}|)\ol{\omega}_{n-1}+k_{n}\ol{\omega}_{n}) 
& \text{if} \quad a=n.  
\end{cases} \label{eq:domino-aod}
\end{eqnarray}
$({\mathfrak g},\geh')=(D_{n+1}^{(2)},B_{n})$:
\begin{eqnarray}
\chi^{(a)}_j = 
\begin{cases}
\sum \mbox{ch}V(k_1\ol{\La}'_1+ k_2\ol{\La}'_2 + \cdots + k_a\ol{\La}'_a) 
 & \text{if} \quad 1 \le a \le n-1,\\
\mbox{ch}V(j\ol{\La}'_n) &  \text{if} \quad a = n.
\end{cases}
\label{eq:domino-db}
\end{eqnarray}
$({\mathfrak g},\geh')=(D_{4}^{(3)},G_2)$:
\begin{equation}\label{eq:domino-dg}
\begin{split}
\chi^{(1)}_j &= 
\sum_{k=0}^j \mbox{ch}V(k\ol{\La}'_1), \\
\chi^{(2)}_j &= 
\sum_{\begin{subarray}{c}
		k_1 + k_2 \le j\\
		k_1,k_2 \in \mathbb{Z}_{\ge 0}
	\end{subarray}}
\min(1+k_2,1+j-k_1-k_2 )(k_1+1) 
\mbox{ch}V(k_1\ol{\La}'_1 + k_2\ol{\La}'_2).
\end{split}
\end{equation}
Here $\mbox{ch}V(\lambda)$ denotes the irreducible 
$\geh'$ character with highest weight $\lambda$.
$\ol{\La}'_a\; (1 \le a \le n)$ are the fundamental weights 
of $\geh'$.
We have also used the convention 
$\ol{\La}'_0 = 0$ and 
$\ol{\omega}_{n-1} = \ol{\La}'_{n-1}+\ol{\La}'_{n}$, 
$\ol{\omega}_{n} = -\ol{\La}'_{n-1}+\ol{\La}'_{n}$.
The sums in (\ref{eq:domino-aec}) --(\ref{eq:domino-db}) 
are taken as follows.
\begin{align*}
(\ref{eq:domino-aec}), (\ref{eq:domino-db})_{a<n}:& \;
k_1, \ldots, k_a \in \Z_{\ge 0}, k_1 + \cdots + k_a \le j,\\
(\ref{eq:domino-aoc}):& \;
k_b, k_{b+2}, \ldots, k_a \in \Z_{\ge 0}, 
k_b + k_{b+2}+ \cdots + k_a = j,\\
(\ref{eq:domino-aeb}), (\ref{eq:domino-aod})_{a<n}:& \;
k_1, \ldots, k_a  \in \Z_{\ge 0}, k_1 + \cdots + k_a \le j,
k_b-j\delta_{a b}\in 2\Z \, (1 \le \forall b \le a),\\
(\ref{eq:domino-aod})_{a=n}:& \;
k_1, \ldots, k_{n-1}\in \Z_{\ge 0}, k_n \in \Z, 
k_1 + \cdots + k_{n-1} + \vert k_n \vert  \le j,\\
& \; k_b-j\delta_{n b} \in 2\Z \, (1 \le \forall b \le n).
\end{align*}
If one depicts the highest weights in 
the sums (\ref{eq:domino-aoc}), 
(\eqref{eq:domino-aec}, \eqref{eq:domino-db}${}_{a<n}$) 
and (\eqref{eq:domino-aeb}${}_{a<n}$, \eqref{eq:domino-aod}${}_{a<n-1}$)
with the Young diagrams as usual, they correspond to those
obtained {}from the $a \times j$ rectangle  by successively
removing $2\times 1$, $1 \times 1$ and $1\times 2$ blocks, respectively.

Let  $x_a = e^{\ol{\La}'_a}$ be a complex variable.
Then the $\geh'$ character 
$\mbox{ch} V(\la)$  belongs to 
${\mathbb Z}[x^{\pm 1}_1, \ldots, x^{\pm 1}_n]$.
By an argument similar to \cite{HKOTY1} one can show
\begin{theorem}\label{th:domino}
If $({\mathfrak g},\geh')=
(A_{2n-1}^{(2)},C_{n}),(A_{2n}^{(2)},C_{n}), (D_{n+1}^{(2)},B_{n})$ 
or $(D^{(3)}_4, G_2)$ (i.e., $\geh' = \overset{\circ}{\geh}$
hence $\ol{\La}'_a = \ol{\La}_a$), 
we have 
\begin{enumerate}
\renewcommand{\theenumi}{\Alph{enumi}}
\renewcommand{\labelenumi}{(\theenumi)}
\item 
$\chi^{(a)}_j = \sum_{\la \in \ol{P}^+} d_\la \mbox{ch}V(\la),
d_\la \in {\mathbb C},\ d_{j\ol{\La}_a} = 1$,\ 
$d_\la = 0$ unless $\la \in j\ol{\La}_a - 
\sum_{b=1}^n {\mathbb Z}_{\ge 0} \alpha_b/\epsilon_b,
\quad 1 \le a \le n, j \in {\mathbb Z}_{\ge 0}$,
\item
$Q^{(a)}_j = \chi^{(a)}_j$ solves the $Q$-system (\ref{eq:qsys}),
\item
$\lim_{j \rightarrow \infty}
\left(\frac{\chi^{(a)}_j}{\chi^{(a)}_{j+1}}\right)
= x^{-1}_a$ 
in the domain $\vert e^{{\alpha}_1} \vert, \ldots, 
\vert e^{{\alpha}_n} \vert > 1$.
\end{enumerate}
\end{theorem}
When $({\mathfrak g},\geh')=(A_{2n-1}^{(2)},D_{n})$ or 
$(A_{2n}^{(2)},B_{n})$, (B) is still valid for 
$\chi^{(a)}_j$ given in (\ref{eq:domino-aeb}) or  (\ref{eq:domino-aod}).
However, (A) and (C) no longer hold.
For example instead of (C), one has

\noindent 
$({\mathfrak g},\geh')=(A_{2n-1}^{(2)},D_{n})$: 
\begin{eqnarray}
\lim_{j \rightarrow \infty}
\left(\frac{\chi^{(a)}_j}{\chi^{(a)}_{j+1}}\right)= 
\begin{cases}
x^{-1}_{a}          & \text{if} \quad 1\le a \le n-2, \\ 
(x_{n-1}x_{n})^{-1} & \text{if} \quad a=n-1, \\ 
x^{-2}_{n-1}        & \text{if} \quad a=n \quad \text{and} \quad 
   \vert e^{\alpha'_{n-1}} \vert > \vert e^{\alpha'_{n}} \vert, \\ 
x^{-2}_{n}          & \text{if} \quad a=n \quad \text{and} \quad 
   \vert e^{\alpha'_{n-1}} \vert < \vert e^{\alpha'_{n}} \vert, 
\end{cases}
\label{eq:lim-aod}
\end{eqnarray}
\noindent 
$({\mathfrak g},\geh')=(A_{2n}^{(2)},B_{n})$: 
\begin{eqnarray}
\lim_{j \rightarrow \infty}
\left(\frac{\chi^{(a)}_j}{\chi^{(a)}_{j+1}}\right)= 
\begin{cases}
x^{-1}_{a} & \text{if} \quad 1\le a \le n-1, \\ 
x_{n}^{-2} & \text{if} \quad a=n , 
\end{cases}
\label{eq:lim-aeb}
\end{eqnarray}
in the domain $\vert e^{\alpha'_1} \vert, \ldots, 
\vert e^{\alpha'_n} \vert > 1$, where 
$\alpha'_a$ denotes a 
simple root of $\geh'$.

The fact that the $\geh'$ character 
$\chi^{(a)}_j$ satisfies the $Q$-system of type 
$\geh = X^{(2)}_N$ can be shown by
using the solution to the nontwisted cases $X^{(1)}_N$ \cite{KR2,HKOTY1}.
All one needs is to restrict the earlier solution 
in terms of $X_N$-characters to those for the subalgebras $\geh'$.
It is done by means of a class of identities known as 
Littlewood's lemma \cite{L}.
See for example the references listed in the sequel.

\noindent
$({\mathfrak g},{\mathfrak g}^{\prime})=(A_{2n-1}^{(2)},C_{n})$:
Theorem in \cite{KS} and eq. (1) in p.492 (or eq. (2) in p.507) of \cite{KT}. \\
$({\mathfrak g},{\mathfrak g}^{\prime})=(A_{2n-1}^{(2)},D_{n})$: 
Theorem in \cite{KS}, 
eq. (1.1) in p.471 of \cite{KT}
and eq. (5) in p.492 (or eq. (3) in p.507) of \cite{KT}. \\
$({\mathfrak g},{\mathfrak g}^{\prime})=(A_{2n}^{(2)},C_{n})$: 
Theorem in \cite{KS} and  eq. (1.5.1) in p.488 of \cite{KT}, 
the equation under eq. (1.5.2) in p.490 of \cite{KT} and 
 eq. (11.9;6) in p.238 of \cite{L}. \\ 
$({\mathfrak g},{\mathfrak g}^{\prime})=(A_{2n}^{(2)},B_{n})$: 
Theorem in \cite{KS}, the equation under eq. (1.1) in p.471 of \cite{KT}
 and  eq. (3) in p.492 (or eq.(1) in p.506) of \cite{KT}. \\
$({\mathfrak g},{\mathfrak g}^{\prime})=(D_{n+1}^{(2)},B_{n})$: 
Theorem 7.1 (B) for $D^{(1)}_{n+1}$ in \cite{HKOTY1} and a branching rule for 
 $O(2n+2)\downarrow O(2n+1)$ \cite{W} and 
$O(n)\downarrow SO(n)$ (p.471 of \cite{KT}).\\
$({\mathfrak g},{\mathfrak g}^{\prime})=(D_{4}^{(3)},G_{2})$: 
Theorem 7.1 (B) for $D^{(1)}_4$ in \cite{HKOTY1} and the 
branching rules $O(8) \downarrow O(7)$ \cite{W} and 
$O(7) \downarrow G_2$ \cite{Wy}.

The solution (\ref{eq:domino-aeb}) of the $A^{(2)}_{2n}$ type 
$Q$-system is relevant to a superalgebra 
$B^{(1)}(0 \vert s)$ \cite{T}.
In Appendix \ref{app:BD} 
we have also included a solution of $B^{(1)}_n$ $Q$-system 
in terms of $D_n$ characters.
%

\subsection{\mathversion{bold} 
Fermionic form and $Q$-system}\label{sec:complete}

In this subsection we 
consider the pair $\overset{\circ}{\geh} \subset \geh$ 
for general $\geh$.

\begin{theorem}\label{th:completeness}
Suppose that a linear combination of $\overset{\circ}{\geh}$ 
characters
\begin{equation*}
Q^{(a)}_j = \sum_{\la \in \ol{P}^+} d_\la \mbox{ch} V(\la) 
\quad 1 \le a \le n, \; j \in {\mathbb Z}_{\ge 0}
\end{equation*}
possesses the properties:
\begin{enumerate}
\renewcommand{\theenumi}{\Alph{enumi}}
\renewcommand{\labelenumi}{(\theenumi)}
\item 
$d_\la \in {\mathbb C},\ d_{j\ol{\La}_a} = 1$,\ 
$d_\la = 0$ unless $\la \in j\ol{\La}_a - 
\sum_{b=1}^n {\mathbb Z}_{\ge 0} \alpha_b/\epsilon_b$,
\item
$\{Q^{(a)}_j\}$ satisfies the $Q$-system (\ref{eq:qsys}),
\item
$\lim_{j \rightarrow \infty}
\left(\frac{Q^{(a)}_j}{Q^{(a)}_{j+1}}\right)
= e^{-\ol{\La}_a}$ 
in the domain $\vert e^{\alpha_1} \vert, \ldots, 
\vert e^{\alpha_n} \vert > 1$.
\end{enumerate}
Then for any finitely many non-negative integers 
$\nu^{(a)}_j$ we have 
\begin{equation}\label{eq:completeness}
\prod_{a=1}^n\prod_{j \ge 1}\left(Q^{(a)}_j\right)^{\nu^{(a)}_j} = 
\sum_{\la} \tilde{M}_\infty(W,\la,1) \mbox{ch}\,V(\la),
\end{equation}
where $\tilde{M}_\infty(W,\la,1)$ is defined by 
(\ref{eq:Mtilde}) with $W = \bigotimes_{a=1}^n\bigotimes_{j \ge 1}
(W^{(a)}_j)^{\otimes \nu^{(a)}_j}$.
The sum $\sum_{\la}$ runs over 
$\la \in 
(\sum_{a=1}^n \sum_{i \ge 1} i \nu^{(a)}_i \ol{\La}_a
 - \sum_{b=1}^n {\mathbb Z}_{\ge 0} \alpha_b/\epsilon_b ) \cap 
\ol{P}^+$. 
\end{theorem}
This theorem has analogous aspects to  
Theorem 8.1 in \cite{HKOTY1}, where $r=1$ case was 
treated.
In particular it asserts the uniqueness of 
$\{Q^{(a)}_j\}$ satisfying (A) -- (C), but does not 
assure the existence, i.e. it does not 
guarantee that the functions 
$\{ \sum_\la \tilde{M}_\infty(W^{(a)}_j,\la,1)\mbox{ch}V(\la)\}$
fulfill\ (B) and (C)\footnote{They do (A). 
Note also the consistency of the condition (B) 
with Proposition \ref{pr:recursion} at $q=1$.}.
In this paper, we have shown the existence  of such $\{Q^{(a)}_j\}$
for ${\mathfrak g} = A_{2n-1}^{(2)}, A_{2n}^{(2)}, D_{n+1}^{(2)}$ 
and $D^{(3)}_4$ in  Theorem \ref{th:domino}.

We omit the proof of  Theorem \ref{th:completeness} 
because it is essentially a repetition 
of the one in \cite{HKOTY1}, whose idea 
goes back to \cite{Ki1,Ki2}.
Another reason is that it is 
a corollary of more general results in  \cite{KN,KNT}.
In these works, 
{\em power series} solutions 
of a generalized  $Q$-system adapted to 
$e^{-j\ol{\La}_a}Q^{(a)}_j$ has been studied.
A proof of uniqueness, existence and two explicit 
constructions have been done
under a certain convergence condition.
(\ref{eq:completeness}) follows {}from them 
upon a further assumption of the Weyl group invariance.

For nontwisted cases $r=1$, the identity 
(\ref{eq:completeness}) can be interpreted 
as a formal (combinatorial) completeness 
of the Bethe ansatz for rational vertex models 
with Yangian $Y(X_N)$-symmetry \cite{KR2,Ki1,Ki2}.
They correspond to 
$q \rightarrow 1$ degeneration of 
trigonometric vertex models with $U'_q(X^{(1)}_N)$ 
symmetry, and the fermionic form $\tilde{M}_\infty(W,\la,1)$ 
arises naturally {}from the rational Bethe equation \cite{OW}
under a string hypothesis.
It is tempting to interpret (\ref{eq:completeness}) 
similarly for twisted cases $r > 1$
in terms of the $\overset{\circ}{\geh}$-symmetry 
of a rational limit of the vertex models 
associated with the twisted $U'_q(X^{(r)}_N)$.
However, at present we do not know a way, even formally, 
to derive $\tilde{M}_\infty(W,\la,1)$ with $r > 1$ {}from 
a rational limit of the Bethe equation \cite{RW} 
and a string hypothesis.
We remark that in a naive $q \rightarrow 1$ limit,
transfer matrix spectrum turns out to 
exhibit a much larger degeneracy than the one 
implied {}from the $\overset{\circ}{\geh}$-symmetry.
%

\appendix
\section{\mathversion{bold}
List of $M_\infty(W^{(a)}_s,\la,q^{-1})$}
\label{app:mlist}

Consider a formal linear combination
of the $\overset{\circ}{\geh}$ modules $V(\la)$ of highest weight $\la$:
$$
{\mathcal W}^{(a)}_s =
\sum_{\lambda \in \overline{P}^{+}}
M_\infty(W^{(a)}_s,\lambda,q^{-1}) V(\lambda).
$$
In this appendix, we give a list of ${\mathcal W}^{(a)}_s$ for
${\mathfrak g}=A^{(2)}_{2n-1},A^{(2)}_{2n},D^{(2)}_{n+1}$ 
and a conjecture 
(in some cases) for ${\mathfrak g}=E^{(2)}_{6},D^{(3)}_{4}$. 
When $q=1$ they reduce to 
(\ref{eq:domino-aec}), (\ref{eq:domino-aoc}) and (\ref{eq:domino-db})
for $A^{(2)}_{2n},A^{(2)}_{2n-1}$ and $D^{(2)}_{n+1}$, respectively.

\noindent
${\mathfrak g}=A^{(2)}_{2n-1}:$
\[
{{\mathcal W}^{(a)}_s}=
\sum_{\lambda} q^{\frac{1}{2}(\ol{\La}_n | s\ol{\La}_a - \lambda)}V(\lambda)
\quad 
(1 \le a \le n),
\]
where the sum $\sum_{\lambda}$ is taken over
$\lambda \in \{k_{e}\ol{\La}_{e}+k_{e+2}\ol{\La}_{e+2}+\dots
+k_{a}\ol{\La}_a\in\ol{P}^+ \, |\, k_{e}+k_{e+2}+\dots
+k_{a-2}+k_{a}=s\}$ with
$\ol{\La}_0=0$ and $e \equiv a \pmod{2},\;e=0$ or $1$.
In the computation of
$M_\infty(W^{(a)}_s,\,k_{e}\ol{\La}_{e}+k_{e+2}\ol{\La}_{e+2}+\dots
+k_{a}\ol{\La}_a,\,q^{-1})$,
the only choice of $\{m^{(c)}_j\}$ such that 
$\forall p^{(c)}_i \ge 0$ is the following:
\newline
if $a=2u$ is even, then
\[
\begin{cases}
m^{(2c-1)}_j=\sum_{b=1}^{\min(c-1,u)} \delta_{j,l_b}+
\sum_{b=1}^{\min(c,u)} \delta_{j,l_b} &
(c \ge 1,\, 2c-1 \le n-1,\,j \ge 1)\\
m^{(2c)}_j=2 \sum_{b=1}^{\min(c,u)} \delta_{j,l_b} &
(c \ge 1,\, 2c \le n-1,\,j \ge 1)\\
m^{(n)}_j=\sum_{b=1}^{u} \delta_{j,l_b} &
(j \ge 1)
\end{cases},
\]
and if $a=2u+1$ is odd, then
\[
\begin{cases}
m^{(1)}_j=0 & (j \ge 1)\\
m^{(2c)}_j=\sum_{b=1}^{\min(c-1,u)} \delta_{j,l_b}+
\sum_{b=1}^{\min(c,u)} \delta_{j,l_b} &
(c \ge 1,\, 2c \le n-1,\,j \ge 1)\\
m^{(2c+1)}_j=2 \sum_{b=1}^{\min(c,u)} \delta_{j,l_b} &
(c \ge 1,\, 2c+1 \le n-1,\,j \ge 1)\\
m^{(n)}_j=\sum_{b=1}^{u} \delta_{j,l_b} &
(j \ge 1)
\end{cases},
\]
where $l_b=s-(k_a+k_{a-2}+\dots +k_{2b+e})$ for 
the both cases.

\noindent
${\mathfrak g}=A^{(2)}_{2n}:$
\[
{{\mathcal W}^{(a)}_s}=
\sum_{\lambda} q^{(\ol{\La}_{n} | s\ol{\La}_a - \lambda)}V(\lambda) 
\quad 
(1 \le a \le n),
\]
where the sum $\sum_{\lambda}$ is taken over
$\lambda \in \{k_{1}\ol{\La}_{1}+\dots
+k_{a}\ol{\La}_a\in\ol{P}^+ \, |\, k_{1}+\dots +k_{a}\le s \}$.
In the computation of
$M_\infty(W^{(a)}_s,\,k_{1}\ol{\La}_{1}+\dots+k_{a}\ol{\La}_a,\,q^{-1})$,
the only choice of $\{m^{(c)}_j\}$ such that 
$\forall p^{(c)}_i \ge 0$ is the following:
\[
m^{(c)}_j=\sum_{b=1}^{\min(c,a)} \delta_{j,l_b} \quad
(1 \le c \le n,\,j \ge 1),
\]
where $l_b=s-(k_b+k_{b+1}+\dots +k_{a})$.

\noindent
${\mathfrak g}=D^{(2)}_{n+1}:$
\[
{{\mathcal W}^{(a)}_s}=
\begin{cases}
\sum_{\lambda} q^{(\ol{\La}_{n} | s\ol{\La}_a - \lambda)}V(\lambda) &
(1 \le a \le n-1) \\
V({s{\ol{\La} }_n}) & (a=n)
\end{cases},
\]
where the sum $\sum_{\lambda}$ is taken over
$\lambda \in \{k_{1}\ol{\La}_{1}+\dots
+k_{a}\ol{\La}_a\in\ol{P}^+ \, |\, k_{1}+\dots +k_{a}\le s \}$.
In the computation of
$M_\infty(W^{(a)}_s,\,k_{1}\ol{\La}_{1}+\dots+k_{a}\ol{\La}_a,\,q^{-1})$,
the only choice of $\{m^{(c)}_j\}$ such that 
$\forall p^{(c)}_i \ge 0$ is the same as $\geh = A^{(2)}_{2n}$ case.

\noindent $\mathfrak{g}=E^{(2)}_{6}$: 
\begin{eqnarray*}
&&{\mathcal W}^{(1)}_1=
q\,V(0) +V({{\ol{\La} }_1}),\\
&&{\mathcal W}^{(2)}_1=
q^3\,V(0) + (q+q^{2})\,V({{\ol{\La} }_1}) 
+ q\, \,V({{\ol{\La} }_4}) +V({{\ol{\La} }_2}),
\\ 
&&{\mathcal W}^{(3)}_1=
(q^4+q^6)\,V(0) + (2q^{3}+q^{4}+q^{5})\,V({{\ol{\La} }_1}) 
+ (2q^2+q^4)\, \,V({{\ol{\La} }_4}) \\ 
&& \quad \qquad + 
(q+q^{2}+q^3)\,V({{\ol{\La} }_2})+
q^{2}\,V(2{{\ol{\La} }_1}) 
+ q\,V({{\ol{\La} }_1}+{{\ol{\La} }_4})
+V({{\ol{\La} }_3}), \\ 
&&{\mathcal W}^{(3)}_2=
(q^8+q^{10}+q^{12})\,V(0) + 
(2q^{6}+4q^{8}+q^{10})\,V({{\ol{\La} }_4}) \\
&& \quad \qquad +(3q^{4}+3q^{6}+q^{8})\,V(2{{\ol{\La} }_4})
+(2q^{4}+q^{5}+6q^{6}+q^{7}+q^{8})\,V({{\ol{\La} }_3})\\ 
&&\quad \qquad+(2q^{2}+2q^{4})\,V({{\ol{\La} }_3}+{{\ol{\La} }_4})+
V(2{{\ol{\La} }_3}) \\ 
&& \quad \qquad+
(q^{5}+2q^{6}+5q^{7}+3q^{8}+2q^{9})\,V({{\ol{\La} }_2})\\ 
&& \quad \qquad +
(2q^{3}+3q^{4}+5q^{5}+2q^{6}+q^{7})\,V({{\ol{\La} }_2}+{{\ol{\La} }_4})\\ 
&& \quad \qquad+(q+q^{2}+q^{3})\,V({{\ol{\La} }_2}+{{\ol{\La} }_3})
+(q^{2}+q^{3}+3q^{4}+q^{5}+q^{6})\,V(2{{\ol{\La} }_2})\\
&&\quad \qquad +(2q^{7}+q^{8}+3q^{9}+q^{10}+q^{11})\,V({{\ol{\La} }_1})\\ 
&& \quad \qquad 
+(5q^{5}+3q^{6}+6q^{7}+q^{8}+q^{9})\,V({{\ol{\La} }_1}+{{\ol{\La} }_4})\\ 
&&\quad \qquad+(2q^{3}+q^{5})\,V({{\ol{\La} }_1}+2{{\ol{\La} }_4})
+(3q^{3}+2q^{4}+3q^{5})\,V({{\ol{\La} }_1}+{{\ol{\La} }_3})\\ 
&&\quad \qquad+q\,V({{\ol{\La} }_1}+{{\ol{\La} }_3}+{{\ol{\La} }_4})
+(2q^{4}+5q^{5}+5q^{6}+3q^{7}+q^{8})\,V({{\ol{\La} }_1}+{{\ol{\La} }_2})\\ 
&&\quad \qquad+
(q^{2}+2q^{3}+q^{4})\,V({{\ol{\La} }_1}+{{\ol{\La} }_2}+{{\ol{\La} }_4})\\ 
&& \quad \qquad 
+(4q^{6}+2q^{7}+4q^{8}+q^{9}+q^{10})\,V(2{{\ol{\La} }_1})\\
&&\quad \qquad+(4q^{4}+q^{5}+2q^{6})\,V(2{{\ol{\La} }_1}+{{\ol{\La} }_4})
+q^{2}\,V(2{{\ol{\La} }_1}+2{{\ol{\La} }_4})
+q^{2}\,V(2{{\ol{\La} }_1}+{{\ol{\La} }_3})\\ 
&&\quad \qquad+(q^{3}+q^{4}+q^{5})\,V(2{{\ol{\La} }_1}+{{\ol{\La} }_2})
+(2q^{5}+q^{6}+q^{7})\,V(3{{\ol{\La} }_1})\\ 
&&\quad \qquad+q^{3}\,V(3{{\ol{\La} }_1}+{{\ol{\La} }_4})
+q^{4}\,V(4{{\ol{\La} }_1}),
\\ 
&&{\mathcal W}^{(4)}_1=
q^{2}\,V(0) + q\,V({{\ol{\La} }_1})+V({{\ol{\La} }_4}).
\end{eqnarray*}
In addition we have a conjecture for
${{\mathcal W}^{(a)}_s}\, (a = 1,2,4)$: 
\begin{eqnarray*}
&&{{\mathcal W}^{(1)}_s}=
\sum_{k=0}^{s} q^{s-k} V(k\ol{\La}_1),
\\ 
&&{{\mathcal W}^{(2)}_s} = 
\sum_{
	\begin{subarray}{c}
		j_1 + j_2 + 2j_3 + j_4 \le s\\
		j_1,j_2,j_3,j_4 \in \mathbb{Z}_{\ge 0}
	\end{subarray}
	}
\min \left( 1+j_2,\,1+s-j_1- j_2-2j_3-j_4 \right)\,
q^{3s-2j_1-3j_2-4j_3-2j_4} \\
&& \qquad \qquad 
\times \begin{bmatrix} j_1+1 \\ 1
\end{bmatrix}_q
V \left( j_1 \ol{\La}_1 +j_2 \ol{\La}_2
+j_3 \ol{\La}_3 +j_4 \ol{\La}_{4} \right),\\
&& {{\mathcal W}^{(4)}_s}=
\sum_{\begin{subarray}{c}
		j_1 + j_4 \le s\\
		j_1,j_4 \in \mathbb{Z}_{\ge 0}
	\end{subarray}}
 q^{2s-j_{1}-2j_{4}}
V \left( j_{1} \ol{\La}_1 + j_{4} \ol{\La}_4 \right).
\end{eqnarray*}
The conjecture has been checked for 
$1 \le s \le 8$. 
In addition, we have observed that the support for 
$V \left( j_1 \ol{\La}_1 +j_2 \ol{\La}_2+j_3 \ol{\La}_3 +j_4 \ol{\La}_{4} \right)$  
in ${{\mathcal W}^{(3)}_s}$ is given by the condition  
$ j_1 +2j_2 +2j_3 +j_4 \le 2s, [(j_2 +1)/2] +j_3 +j_4 \le s, 
j_1,j_2,j_3,j_4 \in {\mathbb Z}_{\ge 0}$
at least for $1 \le s \le 7$\footnote{Unlike $\lceil x \rceil$ in Conjecture 
\ref{conj:fd-module} (1), the symbol $[ x ]$ here denotes the greatest integer 
not exceeding $x$.}. \\ 
\noindent 
$\mathfrak{g}=D^{(3)}_{4}$ : 
\begin{eqnarray*}
&&{{\mathcal W}^{(1)}_1}=
q\,V(0) + V({{\ol{\La} }_1}),
\\
&&{{\mathcal W}^{(2)}_1}=
q^3\,V(0)+(q+q^2)\,V({{\ol{\La} }_1})+V({{\ol{\La} }_2}).
\end{eqnarray*}
In addition we have a conjecture for ${\mathcal W}^{(a)}_s$:
\begin{eqnarray*}
&&{{\mathcal W}^{(1)}_s} =
\sum_{k=0}^{s} q^{s-k} V(k\ol{\La}_1),
\\
&&{\mathcal W}^{(2)}_s = 
\sum_{\begin{subarray}{c}
		j_1 + j_2 \le s\\ 
		j_1,j_2 \in \mathbb{Z}_{\ge 0}
	\end{subarray}}
\min(1+j_2,1+s-j_1-j_2 )
 q^{3s-2j_1-3j_2} \\ 
&& \qquad \qquad \times
\begin{bmatrix} j_1+1 \\ 1
\end{bmatrix}_{q}
V(j_1 \ol{\La}_1+j_2 \ol{\La}_2).
\end{eqnarray*}
The conjecture has been checked for 
$1 \le s \le 20$.
This is consistent with the $q=1$ result (\ref{eq:domino-dg}).

\section{\mathversion{bold} 
Explicit form of Proposition \ref{pr:dualspinon}}
\label{app:regII}

Let us write down (\ref{eq:mpro}) explicitly for each $\geh$.
In describing the data $\nu^\ast$ (\ref{eq:nuast}), we use the symbols:  
\begin{align*}
&\nu = \{\nu_k \}_{k\ge 1}, \quad s\nu = \{s\nu_k\}_{k\ge 1}\quad (s=2,3),\\
&\xi_2 = \{\nu_{\frac{k}{2}} \}_{k \ge 1}, \quad 
\xi_3 = \{\nu_{\frac{k}{3}} \}_{k \ge 1} 
\quad (\nu_{\rm non-integer} = 0),\\
&\eta_2 = \{ 2\nu_{2k}+\nu_{2k-1}+\nu_{2k+1} \}_{k \ge 1},\\
&\eta_3 = \{ 3\nu_{3k}+2(\nu_{3k-1}+\nu_{3k+1})+\nu_{3k-2}+\nu_{3k+2} \}_{k \ge 1}.
\end{align*}

To simplify the presentation 
we find it convenient to introduce the fermionic formula
${\mathcal M}^{(a)}_{q,\nu}(\geh\vert \mu)$ (see (\ref{eq:mcal})) 
formally for $\geh= ``B^{(1)}_2, D^{(1)}_3$ and $A^{(2)}_3$".
These algebras were out of the list 
in the beginning of Section \ref{subsec:affine}, but it is natural to set
\begin{equation*}
\begin{split}
&{\mathcal M}^{(a)}_{q,\nu}(B^{(1)}_2\vert \sum_{c=1}^2\mu_c\alpha_c)
= {\mathcal M}^{(3-a)}_{q,\nu}(C^{(1)}_2\vert \sum_{c=1}^2\mu_{3-c}\alpha_c)
\quad a=1,2,\\
&{\mathcal M}^{(a)}_{q,\nu}(D^{(1)}_3\vert \sum_{c=1}^3\mu_c\alpha_c)
= {\mathcal M}^{(a')}_{q,\nu}(A^{(1)}_3\vert \sum_{c=1}^3\mu_{c'}\alpha_c)
\quad 1 \le a \le 3,\; (1',2',3') = (2,1,3),\\
&{\mathcal M}^{(a)}_{q,\nu}(A^{(2)}_3\vert \sum_{c=1}^2\mu_c\alpha_c)
= {\mathcal M}^{(3-a)}_{q,\nu}(D^{(2)}_3\vert \sum_{c=1}^2\mu_{3-c}\alpha_c)
\quad a=1,2.
\end{split}
\end{equation*}
These relations are relevant to 
$(\geh,b) = (B^{(1)}_n,n-2), (D^{(1)}_n,n-3)$ and $(A^{(2)}_{2n-1},n-2)$
in the following.
In view of (\ref{eq:symmetry}), there is a symmetry 
\begin{equation}\label{eq:small2}
{\mathcal M}^{(a)}_{q,\nu}(\geh\vert \mu) = 
{\mathcal M}^{(\sigma(a))}_{q,\nu}(\geh\vert \sigma(\mu)),
\end{equation}
for the algebras 
$\geh = A^{(1)}_n, D^{(1)}_n$ and $E^{(2)}_6$.
Therefore the factor 
${\mathcal M}^{(2)}_{q, \nu}(A^{(1)}_{2}\vert \mu_1\alpha_1+\mu_2\alpha_2)$
for example in $(\geh,b) = (E^{(1)}_6, 3)$ can also be written 
as ${\mathcal M}^{(1)}_{q, \nu}(A^{(1)}_{2}\vert \mu_2\alpha_1+\mu_1\alpha_2)$.
There are several places in the following results where 
one has a freedom to make such replacements.
Finally we employ the convention 
${\mathcal M}^{(a)}_{q,\nu}(A^{(1)}_0\vert \mu) = 1$.

\noindent
$\geh = A^{(1)}_n:$
\begin{equation*}
\sum_\nu \frac{q^\psi}{\prod_i(q)_{\nu_i}}
{\mathcal M}^{(b-1)}_{q, \nu}(A^{(1)}_{b-1}\vert \sum_{j=1}^{b-1}\mu_j\alpha_j)
{\mathcal M}^{(1)}_{q, \nu}(A^{(1)}_{n-b}\vert  \sum_{j=1}^{n-b}\mu_{b+j}\alpha_j).
\end{equation*}
\noindent
$\geh = B^{(1)}_n:$
\begin{align*}
&b \le n-2,&\quad &\sum_\nu \frac{q^\psi}{\prod_i(q)_{\nu_i}}
{\mathcal M}^{(b-1)}_{q, \nu}(A^{(1)}_{b-1}\vert  \sum_{j=1}^{b-1}\mu_j\alpha_j)
{\mathcal M}^{(1)}_{q,\nu}(B^{(1)}_{n-b}\vert \sum_{j=1}^{n-b}\mu_{b+j}\alpha_j),\\
&b=n-1,&\quad &\sum_\nu \frac{q^\psi}{\prod_i(q)_{\nu_i}}
{\mathcal M}^{(n-2)}_{q, \nu}(A^{(1)}_{n-2}\vert \sum_{j=1}^{n-2}\mu_j\alpha_j)
{\mathcal M}^{(1)}_{q, \xi_2}(A^{(1)}_{1}\vert  \mu_n\alpha_1),\\
&b=n,&\quad &\sum_\nu \frac{q^\psi}{\prod_i(q)_{\nu_i}}
{\mathcal M}^{(n-1)}_{q, \eta_2}
(A^{(1)}_{n-1}\vert  \sum_{j=1}^{n-1}\mu_j\alpha_j).
\end{align*}
\noindent
$\geh = C^{(1)}_n:$
\begin{alignat*}{2}
&b \le n-2, &\quad &\sum_\nu \frac{q^\psi}{\prod_i(q)_{\nu_i}}
{\mathcal M}^{(b-1)}_{q, \nu}(A^{(1)}_{b-1}\vert  \sum_{j=1}^{b-1}\mu_j\alpha_j)
{\mathcal M}^{(1)}_{q, \nu}(C^{(1)}_{n-b}\vert  \sum_{j=1}^{n-b}\mu_{b+j}\alpha_j),\\
&b=n-1,&\quad &\sum_\nu \frac{q^\psi}{\prod_i(q)_{\nu_i}}
{\mathcal M}^{(n-2)}_{q, \nu}(A^{(1)}_{n-2}\vert  \sum_{j=1}^{n-2}\mu_j\alpha_j)
{\mathcal M}^{(1)}_{q, \eta_2}(A^{(1)}_{1}\vert  \mu_n\alpha_1),\\
&b=n,&\quad &\sum_\nu \frac{q^\psi}{\prod_i(q)_{\nu_i}}
{\mathcal M}^{(n-1)}_{q, \xi_2}
(A^{(1)}_{n-1}\vert  \sum_{j=1}^{n-1}\mu_j\alpha_j).
\end{alignat*}
\noindent
$\geh = D^{(1)}_n:$
\begin{alignat*}{2}
&b\le n-3,& &\sum_\nu \frac{q^\psi}{\prod_i(q)_{\nu_i}}
{\mathcal M}^{(b-1)}_{q, \nu}(A^{(1)}_{b-1}\vert  \sum_{j=1}^{b-1}\mu_j\alpha_j)
{\mathcal M}^{(1)}_{q, \nu}
(D^{(1)}_{n-b}\vert  \sum_{j=1}^{n-b}\mu_{b+j}\alpha_j),\\
&b=n-2,& &\sum_\nu \frac{q^\psi}{\prod_i(q)_{\nu_i}}
{\mathcal M}^{(n-3)}_{q, \nu}(A^{(1)}_{n-3}\vert  \sum_{j=1}^{n-3}\mu_j\alpha_j)
{\mathcal M}^{(1)}_{q, \nu}(A^{(1)}_{1}\vert \mu_{n-1}\alpha_1)
{\mathcal M}^{(1)}_{q, \nu}(A^{(1)}_{1}\vert \mu_{n}\alpha_1),\\
&b=n-1, \, n, & \; &\sum_\nu \frac{q^\psi}{\prod_i(q)_{\nu_i}}
{\mathcal M}^{(n-2)}_{q, \nu}(A^{(1)}_{n-1}\vert
 \sum_{j=1}^{n-2}\mu_j\alpha_j+\mu_{\ol{b}}\alpha_{n-1}),
\end{alignat*}
where in the  last line, $\ol{b}$ is determined by 
$\{b, \ol{b}\} = \{n-1,n\}$.

\noindent
$\geh = E^{(1)}_6:$
\begin{alignat*}{2}
&b=1,5,& \quad &\sum_\nu \frac{q^\psi}{\prod_i(q)_{\nu_i}}
{\mathcal M}^{(4)}_{q, \nu}(D^{(1)}_{5}\vert  
\sum_{j=1}^{4}\mu_{((b-3)j-3(b-5))/2}\alpha_j+\mu_6\alpha_5),\\
&b=2,4,& \quad &\sum_\nu \frac{q^\psi}{\prod_i(q)_{\nu_i}}
{\mathcal M}^{(1)}_{q, \nu}(A^{(1)}_{1}\vert  \mu_{2b-3}\alpha_1)
{\mathcal M}^{(3)}_{q, \nu}(A^{(1)}_{4}\vert  
\sum_{j=1}^3\mu_{(b-3)j-3(b-4)}\alpha_j+\mu_6\alpha_4),\\
&b=3,& \quad &\sum_\nu \frac{q^\psi}{\prod_i(q)_{\nu_i}}
{\mathcal M}^{(2)}_{q, \nu}(A^{(1)}_{2}\vert \mu_1\alpha_1+\mu_2\alpha_2)
{\mathcal M}^{(1)}_{q, \nu}(A^{(1)}_{2}\vert \mu_4\alpha_1+\mu_5\alpha_2)
{\mathcal M}^{(1)}_{q, \nu}(A^{(1)}_{1}\vert  \mu_6\alpha_1),\\
&b=6,& \quad &\sum_\nu \frac{q^\psi}{\prod_i(q)_{\nu_i}}
{\mathcal M}^{(3)}_{q, \nu}(A^{(1)}_{5}\vert  \sum_{j=1}^{5}\mu_j\alpha_j).
\end{alignat*}
\noindent
$\geh = E^{(1)}_7:$
\begin{align*}
&b=1,\quad \sum_\nu \frac{q^\psi}{\prod_i(q)_{\nu_i}}
{\mathcal M}^{(5)}_{q, \nu}(D^{(1)}_{6}\vert  
\sum_{j=1}^{5}\mu_{7-j}\alpha_j+\mu_7\alpha_6),\\
&b=2,\quad \sum_\nu \frac{q^\psi}{\prod_i(q)_{\nu_i}}
{\mathcal M}^{(1)}_{q, \nu}(A^{(1)}_{1}\vert  \mu_1\alpha_1)
{\mathcal M}^{(2)}_{q, \nu}(A^{(1)}_{5}\vert  
\mu_7\alpha_1+\sum_{j=2}^{5}\mu_{j+1}\alpha_j),\\
&b=3,\quad \sum_\nu \frac{q^\psi}{\prod_i(q)_{\nu_i}}
{\mathcal M}^{(2)}_{q, \nu}(A^{(1)}_{2}\vert \mu_1\alpha_1+\mu_2\alpha_2)
{\mathcal M}^{(1)}_{q, \nu}(A^{(1)}_{3}\vert  \sum_{j=1}^{3}\mu_{j+3}\alpha_j)
{\mathcal M}^{(1)}_{q, \nu}(A^{(1)}_{1}\vert  \mu_7\alpha_1),\\
&b=4,\quad \sum_\nu \frac{q^\psi}{\prod_i(q)_{\nu_i}}
{\mathcal M}^{(3)}_{q, \nu}(A^{(1)}_{4}\vert  
\sum_{j=1}^{3}\mu_j\alpha_j+\mu_7\alpha_4)
{\mathcal M}^{(1)}_{q, \nu}(A^{(1)}_{2}\vert \mu_5\alpha_1+\mu_6\alpha_2),\\
&b=5,\quad \sum_\nu \frac{q^\psi}{\prod_i(q)_{\nu_i}}
{\mathcal M}^{(4)}_{q, \nu}(D^{(1)}_{5}\vert  
\sum_{j=1}^{4}\mu_j\alpha_j+\mu_7\alpha_5)
{\mathcal M}^{(1)}_{q, \nu}(A^{(1)}_{1}\vert \mu_6\alpha_1),\\
&b=6,\quad \sum_\nu \frac{q^\psi}{\prod_i(q)_{\nu_i}}
{\mathcal M}^{(5)}_{q, \nu}(E^{(1)}_{6}\vert  
\sum_{j=1}^{5}\mu_j\alpha_j+\mu_7\alpha_6),\\
&b=7,\quad \sum_\nu \frac{q^\psi}{\prod_i(q)_{\nu_i}}
{\mathcal M}^{(3)}_{q, \nu}(A^{(1)}_{6}\vert  
\sum_{j=1}^{6}\mu_j\alpha_j).
\end{align*}

\noindent
$\geh = E^{(1)}_8:$
\begin{align*}
&b=1,\quad \sum_\nu \frac{q^\psi}{\prod_i(q)_{\nu_i}}
{\mathcal M}^{(6)}_{q, \nu}(E^{(1)}_{7}\vert  
\sum_{j=1}^{6}\mu_{8-j}\alpha_j+\mu_8\alpha_7),\\
&b=2,\quad \sum_\nu \frac{q^\psi}{\prod_i(q)_{\nu_i}}
{\mathcal M}^{(1)}_{q, \nu}(A^{(1)}_{1}\vert  \mu_1\alpha_1)
{\mathcal M}^{(1)}_{q, \nu}(E^{(1)}_{6}\vert  
\sum_{j=1}^{5}\mu_{j+2}\alpha_j+\mu_8\alpha_6),\\
&b=3,\quad \sum_\nu \frac{q^\psi}{\prod_i(q)_{\nu_i}}
{\mathcal M}^{(2)}_{q, \nu}(A^{(1)}_{2}\vert \mu_1\alpha_1+\mu_2\alpha_2)
{\mathcal M}^{(4)}_{q, \nu}(D^{(1)}_{5}\vert  
\sum_{j=1}^{4}\mu_{8-j}\alpha_j+\mu_8\alpha_5),\\
&b=4,\quad \sum_\nu \frac{q^\psi}{\prod_i(q)_{\nu_i}}
{\mathcal M}^{(3)}_{q, \nu}(A^{(1)}_{3}\vert \sum_{j=1}^{3}\mu_j\alpha_j)
{\mathcal M}^{(2)}_{q, \nu}(A^{(1)}_{4}\vert 
\mu_8\alpha_1+\sum_{j=2}^{4}\mu_{j+3}\alpha_j),\\
&b=5,\quad \sum_\nu \frac{q^\psi}{\prod_i(q)_{\nu_i}}
{\mathcal M}^{(4)}_{q, \nu}(A^{(1)}_{4}\vert \sum_{j=1}^{4}\mu_j\alpha_j)
{\mathcal M}^{(1)}_{q, \nu}(A^{(1)}_{1}\vert \mu_8\alpha_1)
{\mathcal M}^{(1)}_{q, \nu}(A^{(1)}_{2}\vert \mu_6\alpha_1+\mu_7\alpha_2),\\
&b=6,\quad \sum_\nu \frac{q^\psi}{\prod_i(q)_{\nu_i}}
{\mathcal M}^{(5)}_{q, \nu}(A^{(1)}_{6}\vert  
\sum_{j=1}^{5}\mu_j\alpha_j+\mu_8\alpha_6)
{\mathcal M}^{(1)}_{q, \nu}(A^{(1)}_{1}\vert \mu_7\alpha_1),\\
&b=7,\quad \sum_\nu \frac{q^\psi}{\prod_i(q)_{\nu_i}}
{\mathcal M}^{(6)}_{q, \nu}(D^{(1)}_{7}\vert  
\sum_{j=1}^{6}\mu_j\alpha_j+\mu_8\alpha_7),\\
&b=8,\quad \sum_\nu \frac{q^\psi}{\prod_i(q)_{\nu_i}}
{\mathcal M}^{(5)}_{q, \nu}(A^{(1)}_{7}\vert  
\sum_{j=1}^{7}\mu_j\alpha_j).
\end{align*}
\noindent
$\geh = F^{(1)}_4:$
\begin{align*}
&b=1,\quad \sum_\nu \frac{q^\psi}{\prod_i(q)_{\nu_i}}
{\mathcal M}^{(3)}_{q, \nu}(C^{(1)}_{3}\vert  \sum_{j=1}^{3}\mu_{5-j}\alpha_j),\\
&b=2,\quad \sum_\nu \frac{q^\psi}{\prod_i(q)_{\nu_i}}
{\mathcal M}^{(1)}_{q, \nu}(A^{(1)}_{1}\vert  \mu_1\alpha_1)
{\mathcal M}^{(1)}_{q, \xi_2}(A^{(1)}_{2}\vert  \mu_3\alpha_1+\mu_4\alpha_2),\\
&b=3,\quad \sum_\nu \frac{q^\psi}{\prod_i(q)_{\nu_i}}
{\mathcal M}^{(2)}_{q, \eta_2}(A^{(1)}_{2}\vert \mu_1\alpha_1+\mu_2\alpha_2)
{\mathcal M}^{(1)}_{q, \nu}(A^{(1)}_{1}\vert  \mu_4\alpha_1),\\
&b=4,\quad \sum_\nu \frac{q^\psi}{\prod_i(q)_{\nu_i}}
{\mathcal M}^{(3)}_{q, \nu}(B^{(1)}_{3}\vert  \sum_{j=1}^{3}\mu_j\alpha_j).
\end{align*}
\noindent
$\geh = G^{(1)}_2:$
\begin{align*}
&b=1,\quad \sum_\nu \frac{q^\psi}{\prod_i(q)_{\nu_i}}
{\mathcal M}^{(1)}_{q, \xi_3}(A^{(1)}_{1}\vert \mu_{2}\alpha_1),\\
&b=2,\quad \sum_\nu \frac{q^\psi}{\prod_i(q)_{\nu_i}}
{\mathcal M}^{(1)}_{q, \eta_3}(A^{(1)}_{1}\vert  \mu_1\alpha_1).
\end{align*}
\noindent
$\geh = A^{(2)}_{2n}:$
\begin{alignat*}{2}
&b \le n-1,&\quad &\sum_\nu \frac{q^{2\psi}}{\prod_i(q^2)_{\nu_i}}
{\mathcal M}^{(b-1)}_{q^2, \nu}
(A^{(1)}_{b-1}\vert  \sum_{j=1}^{b-1}\mu_j\alpha_j)
{\mathcal M}^{(1)}_{q,\nu}
(A^{(2)}_{2n-2b}\vert \sum_{j=1}^{n-b}\mu_{b+j}\alpha_j),\\
&b=n,&\quad &\sum_\nu \frac{q^\psi}{\prod_i(q^2)_{\nu_i}}
{\mathcal M}^{(n-1)}_{q^2, \nu}
(A^{(1)}_{n-1}\vert  \sum_{j=1}^{n-1}\mu_j\alpha_j).
\end{alignat*}
\noindent
$\geh = A^{(2)}_{2n-1}:$
\begin{alignat*}{2}
&b\le n-2,&\quad &\sum_\nu \frac{q^\psi}{\prod_i(q)_{\nu_i}}
{\mathcal M}^{(b-1)}_{q, \nu}(A^{(1)}_{b-1}\vert  \sum_{j=1}^{b-1}\mu_j\alpha_j)
{\mathcal M}^{(1)}_{q,\nu}
(A^{(2)}_{2n-2b-1}\vert \sum_{j=1}^{n-b}\mu_{b+j}\alpha_j),\\
&b=n-1,&\quad &\sum_\nu \frac{q^\psi}{\prod_i(q)_{\nu_i}}
{\mathcal M}^{(n-2)}_{q, \nu}(A^{(1)}_{n-2}\vert \sum_{j=1}^{n-2}\mu_j\alpha_j)
{\mathcal M}^{(1)}_{q^2, \nu}(A^{(1)}_{1}\vert  \mu_n\alpha_1),\\
&b=n,&\quad &\sum_\nu \frac{q^{2\psi}}{\prod_i(q^2)_{\nu_i}}
{\mathcal M}^{(n-1)}_{q, 2\nu}
(A^{(1)}_{n-1}\vert  \sum_{j=1}^{n-1}\mu_j\alpha_j).
\end{alignat*}
\noindent
$\geh = D^{(2)}_{n+1}:$
\begin{alignat*}{2}
&b \le n-2,&\quad &\sum_\nu \frac{q^{2\psi}}{\prod_i(q^2)_{\nu_i}}
{\mathcal M}^{(b-1)}_{q^2, \nu}(A^{(1)}_{b-1}\vert  \sum_{j=1}^{b-1}\mu_j\alpha_j)
{\mathcal M}^{(1)}_{q,\nu}
(D^{(2)}_{n-b+1}\vert \sum_{j=1}^{n-b}\mu_{b+j}\alpha_j),\\
&b=n-1,&\quad &\sum_\nu \frac{q^{2\psi}}{\prod_i(q^2)_{\nu_i}}
{\mathcal M}^{(n-2)}_{q^2, \nu}(A^{(1)}_{n-2}\vert \sum_{j=1}^{n-2}\mu_j\alpha_j)
{\mathcal M}^{(1)}_{q, 2\nu}(A^{(1)}_{1}\vert  \mu_n\alpha_1),\\
&b=n,&\quad &\sum_\nu \frac{q^{\psi}}{\prod_i(q)_{\nu_i}}
{\mathcal M}^{(n-1)}_{q^2, \nu}
(A^{(1)}_{n-1}\vert  \sum_{j=1}^{n-1}\mu_j\alpha_j).
\end{alignat*}
\noindent
$\geh = E^{(2)}_6:$
\begin{align*}
&b=1,\quad \sum_\nu \frac{q^\psi}{\prod_i(q)_{\nu_i}}
{\mathcal M}^{(3)}_{q, \nu}(D^{(2)}_{4}\vert  \sum_{j=1}^{3}\mu_{5-j}\alpha_j),\\
&b=2,\quad \sum_\nu \frac{q^\psi}{\prod_i(q)_{\nu_i}}
{\mathcal M}^{(1)}_{q, \nu}(A^{(1)}_{1}\vert  \mu_1\alpha_1)
{\mathcal M}^{(1)}_{q^2, \nu}(A^{(1)}_{2}\vert  \mu_3\alpha_1+\mu_4\alpha_2),\\
&b=3,\quad \sum_\nu \frac{q^{2\psi}}{\prod_i(q^2)_{\nu_i}}
{\mathcal M}^{(2)}_{q, 2\nu}(A^{(1)}_{2}\vert \mu_1\alpha_1+\mu_2\alpha_2)
{\mathcal M}^{(1)}_{q^2, \nu}(A^{(1)}_{1}\vert  \mu_4\alpha_1),\\
&b=4,\quad \sum_\nu \frac{q^{2\psi}}{\prod_i(q^2)_{\nu_i}}
{\mathcal M}^{(3)}_{q, \nu}(A^{(2)}_{5}\vert  \sum_{j=1}^{3}\mu_j\alpha_j).
\end{align*}
\noindent
$\geh = D^{(3)}_4:$
\begin{align*}
&b=1, \quad \sum_\nu \frac{q^\psi}{\prod_i(q)_{\nu_i}}
{\mathcal M}^{(1)}_{q^3, \nu}(A^{(1)}_{1}\vert \mu_{2}\alpha_1),\\
&b=2, \quad \sum_\nu \frac{q^{3\psi}}{\prod_i(q^3)_{\nu_i}}
{\mathcal M}^{(1)}_{q, 3\nu}(A^{(1)}_{1}\vert  \mu_1\alpha_1).
\end{align*}
The decomposition into subalgebras 
for nontwisted cases exhibits the same pattern as 
the Bethe ansatz calculation of central charges 
in ``regime II" RSOS models. 
See eq.(3.9) in \cite{Ku}.


\section{\mathversion{bold}
Solution to $B^{(1)}_{n}$ $Q$-system by $D_{n}$ characters}
\label{app:BD}

Here we present a solution of the unrestricted $Q$-system 
(\ref{eq:qsys}) for $\geh = B^{(1)}_{n}$
in terms of the character of a classical subalgebra 
$\geh'=D_{n} \subset B^{(1)}_{n}$. 
If $\geh$ is nontwisted and non-exceptional, i.e., 
if $\geh= X^{(1)}_n = A^{(1)}_n, B^{(1)}_n, C^{(1)}_n, D^{(1)}_n$, 
this is the unique case where the classical subalgebra $\geh' \subset \geh$ 
obtained by removing an end point of the Dynkin diagram differs {}from 
$X_n$.
A result for $({\mathfrak g},\geh')=(B_{n}^{(1)},B_{n})$ 
case is available in \cite{KR2,HKOTY1}.
Let us introduce the linear combination of characters:
\begin{eqnarray}
\chi^{(a)}_j =
\begin{cases}
\sum \mbox{ch}V(k_{1}\ol{\La}'_{1} + k_{2}\ol{\La}'_{ 2} 
+ \cdots +k_a\ol{\La}'_a)  & \text{if} \quad 1 \le a \le n-2,  \\ 
\sum \mbox{ch}V(k_{1}\ol{\La}'_{1} + k_{2}\ol{\La}'_{ 2} 
+ \cdots +k_{n-2}\ol{\La}'_{n-2} \\ 
 \hspace{90pt} +k_{n-1}\ol{\omega}_{n-1})
 & \text{if} \quad a=n-1,  \\ 
\sum \mbox{ch}V(k_{1}\ol{\La}'_{1} + k_{2}\ol{\La}'_{ 2} 
+ \cdots +k_{n-2}\ol{\La}'_{n-2} \\ 
 \hspace{40pt} 
+(k_{n-1}+\frac{1}{2}|k_{n}|)\ol{\omega}_{n-1}
+\frac{1}{2}k_{n}\ol{\omega}_{n}) 
& \text{if} \quad a=n.  
\end{cases} \label{eq:domino-bd}
\end{eqnarray}
Here $\mbox{ch}V(\lambda)$ denotes the irreducible 
$D_{n}$ character with highest weight $\lambda$.
$\ol{\La}'_a\; (1 \le a \le n)$ are the fundamental weights 
of $D_{n}$. 
We have also used the convention 
$\ol{\La}'_0 = 0$ and 
$\ol{\omega}_{n-1} = \ol{\La}'_{n-1}+\ol{\La}'_{n}$, 
$\ol{\omega}_{n} = -\ol{\La}'_{n-1}+\ol{\La}'_{n}$. 
The sums in (\ref{eq:domino-bd}) are taken as follows.
\begin{eqnarray*}
&& {\rm For} \quad a \le n-1: 
k_1, \ldots, k_{a}\in \Z_{\ge 0},\; 
k_1 \cdots + k_a \le j,\\
&&{\rm For} \quad a=n:
 k_1, \ldots, k_{n-1}\in \Z_{\ge 0},\;  k_n \in  \Z,\; 
2(k_1 + \cdots +k_{n-1})+ |k_n| \le j, \\ 
&& k_n - j \in 2 \Z .
\end{eqnarray*}
If one depicts the highest weights in 
(\ref{eq:domino-bd})${}_{a \le n-2}$ 
with Young diagrams as usual, they correspond to those
obtained {}from the $a \times j$ rectangle  by successively
removing $1 \times 1$ pieces. 
We note that 
$\chi^{(n-1)}_j = \mbox{ch} V(j\ol{\La}'_{n-1}+j\ol{\La}'_n) + \cdots$ and 
$\chi^{(n)}_j = 
\mbox{ch} V(j\ol{\La}'_{n-1}) + \mbox{ch} V(j\ol{\La}'_{n-1}) + \cdots$.

The $D_{n}$ character 
$\mbox{ch} V(\la)$  belongs to 
${\mathbb Z}[x^{\pm 1}_1, \ldots, x^{\pm 1}_n]$, where 
$x_a = e^{\ol{\La}'_a}$ is a complex variable.
One can show 
Theorem \ref{th:domino} (B) for 
$({\mathfrak g},\geh')=(B_{n}^{(1)},D_{n})$ for the 
$\chi^{(a)}_j$ given by (\ref{eq:domino-bd}).
However, (A) and (C) no longer hold.
For example instead of (C), one has
\begin{eqnarray}
\lim_{j \rightarrow \infty}
\left(\frac{\chi^{(a)}_j}{\chi^{(a)}_{j+1}}\right)= 
\begin{cases}
x^{-1}_{a}          & \text{if} \quad 1\le a \le n-2, \\ 
(x_{n-1}x_{n})^{-1} & \text{if} \quad a=n-1, \\ 
x^{-1}_{n-1}        & \text{if} \quad a=n \quad \text{and} \quad 
   \vert e^{\alpha'_{n-1}} \vert > \vert e^{\alpha'_{n}} \vert, \\ 
x^{-1}_{n}          & \text{if} \quad a=n \quad \text{and} \quad 
   \vert e^{\alpha'_{n-1}} \vert < \vert e^{\alpha'_{n}} \vert, 
\end{cases}
\label{eq:lim-bd}
\end{eqnarray}
in the domain $\vert e^{\alpha'_1} \vert, \ldots, 
\vert e^{\alpha'_n} \vert > 1$, where 
$\alpha'_a$ denotes a 
simple root of $D_{n}$.

The fact that the $D_{n}$ character 
$\chi^{(a)}_j$ satisfies the $Q$-system of type 
$\geh = B^{(1)}_n$ can be shown by
using Theorem 7.1 (B) for $B^{(1)}_{n}$ in \cite{HKOTY1}, 
 a branching rule for $O(2n+1)\downarrow O(2n)$  
 (see for example, \cite{W}) and 
 $O(n)\downarrow SO(n)$ 
 (see for example, p.471 of  \cite{KT}).


\section{\mathversion{bold} Examples of combinatorial $R$ matrices}
\label{app:example}

In this Appendix we provide explicit data of
the isomorphism $B^{1,s_1} \ot B^{1,s_2} \simeq B^{1,s_2} \ot B^{1,s_1}$
$(1 \le s_2 \le s_1 \le 2)$ and energy function for 
$\geh = B^{(1)}_3, C^{(1)}_2, D^{(1)}_4, A^{(2)}_3$, 
$A^{(2)}_4, D^{(2)}_3$ 
and $D^{(3)}_4$.
Suppose $b \ot c \simeq \tilde{c} \ot \tilde{b}$ 
and $H(b \ot c) = h$ under the 
isomorphism $B^{1,s_1}\ot B^{1,s_2} \simeq B^{1,s_2}\ot B^{1,s_1}$.
For each $B^{1,s_1}\ot B^{1,s_2}$, we list the data in a table 
whose row is labeled with $b$ and column with $c$.
At the position corresponding to the $b$-th row and 
$c$-th column, we put $\tilde{c} \cdot \tilde{b}_h$ if $s_1\neq s_2$.
If $s_1 = s_2$, the isomorphism is trivial in that 
$b = \tilde{c}$ and $c = \tilde{b}$, therefore we just put $h$.
Due to the limitation of space, tables 
for $B^{1,2} \ot B^{1,2}$ for $B^{(1)}_3, D^{(1)}_4$ and $D^{(3)}_4$
are not included.

Let us comment on the relevant earlier results.
For $\geh= A^{(1)}_n$, there is a simple algorithm to find these data
in \cite{NY1}.
For $\geh$ non-exceptional type other than $A^{(1)}_n$, 
an algorithm is available for 
general $s_1, s_2$ based on an insertion scheme \cite{HKOT1,HKOT2}.
If $\geh$ is non-exceptional and $s_1=s_2$, 
there is a piecewise linear formula
to find $H(b\ot c)$ for general $s_1$ in \cite{KKM}.
(Although for $C^{(1)}_n$ \cite{KKM} only treats  $s_1$ even case,
the results therein are valid also for odd $s_1$ under an 
appropriate adjustment.)

In \cite{KKM}, elements of crystals $B^{1,s}$ 
are labeled by the vectors with non-negative integer coordinates
$(x_1,\ldots,x_n,\ol{x}_n,\ldots,\ol{x}_1)$ 
($\geh = C^{(1)}_n, D^{(1)}_n, A^{(2)}_{2n}, A^{(2)}_{2n-1}$) or 
$(x_1,\ldots,x_n,x_0,\ol{x}_n,\ldots,\ol{x}_1)$ 
($\geh = B^{(1)}_n, D^{(1)}_{n+1}$).
Here we represent them as 
$1^{x_1}\ldots n^{x_n}(0^{x_0})\ol{n}^{\ol{x}_n}
\ldots \ol{1}^{\ol{x}_1}$
if not all of the coordinates are zero.
In case they are all zero, such a crystal element is denoted by $\phi$.
All the tables given below are  typeset by automatically transforming the 
output {}from a computer.

{
\begin{table}[tbp]
	\caption{$B^{1,1} \otimes B^{1,1}$ for $B^{(1)}_3$}
	\label{tab:B13_11}
    \newcommand{\lw}[1]{\smash{\lower1.6ex\hbox{#1}}}
	\begin{center}
	\begin{tabular}{c|ccccccc}	\hline
	\lw{$b$} & \multicolumn{7}{c}{$c$} \\

& $0$ & $1$ & $2$ & $3$ & $\ol{3}$ & $\ol{2}$ & $\ol{1}$ \\  \hline$0$ & 1 
& 2 & 2 & 2 & 1 & 1 & 1 \\ $1$ & 1 & 2 & 1 & 1 & 1 & 1 & 0 \\ $2$ & 1 & 2 & 2 
& 1 & 1 & 1 & 1 \\ $3$ & 1 & 2 & 2 & 2 & 1 & 1 & 1 \\ $\ol{3}$ & 2 & 2 & 2 & 
2 & 2 & 1 & 1 \\ $\ol{2}$ & 2 & 2 & 2 & 2 & 2 & 2 & 1 \\ $\ol{1}$ & 2 & 2 & 2 
& 2 & 2 & 2 & 2 \\  \hline \end{tabular}\end{center}\end{table}

\begin{table}[tbp]
	\caption{$B^{1,1} \otimes B^{1,1}$ for $D^{(1)}_4$}
	\label{tab:D14_11}
    \newcommand{\lw}[1]{\smash{\lower1.6ex\hbox{#1}}}
	\begin{center}
	{
	\begin{tabular}{c|cccccccc}	\hline
	\lw{$b$} & \multicolumn{8}{c}{$c$} \\
	
& $1$ & $2$ & $3$ & $4$ & $\ol{4}$ & $\ol{3}$ & $\ol{2}$ & $\ol{1}$ \\  
\hline$1$ & 2 & 1 & 1 & 1 & 1 & 1 & 1 & 0 \\ $2$ & 2 & 2 & 1 & 1 & 1 & 1 & 1 
& 1 \\ $3$ & 2 & 2 & 2 & 1 & 1 & 1 & 1 & 1 \\ $4$ & 2 & 2 & 2 & 2 & 1 & 1 & 1 
& 1 \\ $\ol{4}$ & 2 & 2 & 2 & 1 & 2 & 1 & 1 & 1 \\ $\ol{3}$ & 2 & 2 & 2 & 2 & 
2 & 2 & 1 & 1 \\ $\ol{2}$ & 2 & 2 & 2 & 2 & 2 & 2 & 2 & 1 \\ $\ol{1}$ & 2 & 2 
& 2 & 2 & 2 & 2 & 2 & 2 \\  \hline \end{tabular}}\end{center}\end{table}

\begin{table}[tbp]
	\caption{$B^{1,2} \otimes B^{1,1}$ for $B^{(1)}_3$}
	\label{tab:B13_21}
    \newcommand{\lw}[1]{\smash{\lower1.6ex\hbox{#1}}}
	\begin{center}
    {
	\begin{tabular}{c|ccccccc}	\hline
	\lw{$b$} & \multicolumn{7}{c}{$c$} \\

& $0$ & $1$ & $2$ & $3$ & $\ol{3}$ & $\ol{2}$ & $\ol{1}$ \\  \hline$10$ & 
$0 \cdot 10_1$ & $0 \cdot 11_2$ & $1 \cdot 20_1$ & $1 \cdot 30_1$ & $0 \cdot 
1\ol{3}_1$ & $0 \cdot 1\ol{2}_1$ & $1 \cdot 0\ol{1}_0$ \\ $20$ & $0 \cdot 
20_1$ & $0 \cdot 12_2$ & $0 \cdot 22_2$ & $2 \cdot 30_1$ & $0 \cdot 
2\ol{3}_1$ & $0 \cdot 1\ol{1}_1$ & $0 \cdot 2\ol{1}_1$ \\ $30$ & $0 \cdot 
30_1$ & $0 \cdot 13_2$ & $0 \cdot 23_2$ & $0 \cdot 33_2$ & $0 \cdot 
2\ol{2}_1$ & $0 \cdot 3\ol{2}_1$ & $0 \cdot 3\ol{1}_1$ \\ $0\ol{3}$ & $0 
\cdot 0\ol{3}_1$ & $\ol{3} \cdot 10_2$ & $\ol{3} \cdot 20_2$ & $\ol{3} \cdot 
30_2$ & $0 \cdot \ol{3}\ol{3}_1$ & $\ol{3} \cdot 0\ol{2}_1$ & $\ol{3} \cdot 
0\ol{1}_1$ \\ $0\ol{2}$ & $0 \cdot 0\ol{2}_1$ & $\ol{2} \cdot 10_2$ & $\ol{2} 
\cdot 20_2$ & $\ol{2} \cdot 30_2$ & $0 \cdot \ol{3}\ol{2}_1$ & $0 \cdot 
\ol{2}\ol{2}_1$ & $\ol{2} \cdot 0\ol{1}_1$ \\ $0\ol{1}$ & $0 \cdot 0\ol{1}_1$ 
& $\ol{1} \cdot 10_2$ & $\ol{1} \cdot 20_2$ & $\ol{1} \cdot 30_2$ & $0 \cdot 
\ol{3}\ol{1}_1$ & $0 \cdot \ol{2}\ol{1}_1$ & $0 \cdot \ol{1}\ol{1}_1$ \\ $11$ 
& $1 \cdot 10_1$ & $1 \cdot 11_2$ & $1 \cdot 12_1$ & $1 \cdot 13_1$ & $1 
\cdot 1\ol{3}_1$ & $1 \cdot 1\ol{2}_1$ & $1 \cdot 1\ol{1}_0$ \\ $12$ & $2 
\cdot 10_1$ & $2 \cdot 11_2$ & $1 \cdot 22_1$ & $2 \cdot 13_1$ & $2 \cdot 
1\ol{3}_1$ & $2 \cdot 1\ol{2}_1$ & $1 \cdot 2\ol{1}_0$ \\ $13$ & $3 \cdot 
10_1$ & $3 \cdot 11_2$ & $1 \cdot 23_1$ & $1 \cdot 33_1$ & $3 \cdot 
1\ol{3}_1$ & $3 \cdot 1\ol{2}_1$ & $1 \cdot 3\ol{1}_0$ \\ $1\ol{3}$ & $1 
\cdot 0\ol{3}_1$ & $\ol{3} \cdot 11_2$ & $1 \cdot 2\ol{3}_1$ & $1 \cdot 
3\ol{3}_1$ & $1 \cdot \ol{3}\ol{3}_1$ & $\ol{3} \cdot 1\ol{2}_1$ & $1 \cdot 
\ol{3}\ol{1}_0$ \\ $1\ol{2}$ & $1 \cdot 0\ol{2}_1$ & $\ol{2} \cdot 11_2$ & $1 
\cdot 2\ol{2}_1$ & $1 \cdot 3\ol{2}_1$ & $1 \cdot \ol{3}\ol{2}_1$ & $1 \cdot 
\ol{2}\ol{2}_1$ & $1 \cdot \ol{2}\ol{1}_0$ \\ $1\ol{1}$ & $2 \cdot 0\ol{2}_1$ 
& $\ol{1} \cdot 11_2$ & $2 \cdot 2\ol{2}_1$ & $2 \cdot 3\ol{2}_1$ & $2 \cdot 
\ol{3}\ol{2}_1$ & $2 \cdot \ol{2}\ol{2}_1$ & $1 \cdot \ol{1}\ol{1}_0$ \\ $22$ 
& $2 \cdot 20_1$ & $2 \cdot 12_2$ & $2 \cdot 22_2$ & $2 \cdot 23_1$ & $2 
\cdot 2\ol{3}_1$ & $2 \cdot 1\ol{1}_1$ & $2 \cdot 2\ol{1}_1$ \\ $23$ & $3 
\cdot 20_1$ & $3 \cdot 12_2$ & $3 \cdot 22_2$ & $2 \cdot 33_1$ & $3 \cdot 
2\ol{3}_1$ & $3 \cdot 1\ol{1}_1$ & $3 \cdot 2\ol{1}_1$ \\ $2\ol{3}$ & $2 
\cdot 0\ol{3}_1$ & $\ol{3} \cdot 12_2$ & $\ol{3} \cdot 22_2$ & $2 \cdot 
3\ol{3}_1$ & $2 \cdot \ol{3}\ol{3}_1$ & $\ol{3} \cdot 1\ol{1}_1$ & $\ol{3} 
\cdot 2\ol{1}_1$ \\ $2\ol{2}$ & $3 \cdot 0\ol{3}_1$ & $\ol{2} \cdot 12_2$ & 
$\ol{2} \cdot 22_2$ & $3 \cdot 3\ol{3}_1$ & $3 \cdot \ol{3}\ol{3}_1$ & 
$\ol{2} \cdot 1\ol{1}_1$ & $\ol{2} \cdot 2\ol{1}_1$ \\ $2\ol{1}$ & $2 \cdot 
0\ol{1}_1$ & $\ol{1} \cdot 12_2$ & $\ol{1} \cdot 22_2$ & $2 \cdot 3\ol{1}_1$ 
& $2 \cdot \ol{3}\ol{1}_1$ & $2 \cdot \ol{2}\ol{1}_1$ & $2 \cdot 
\ol{1}\ol{1}_1$ \\ $33$ & $3 \cdot 30_1$ & $3 \cdot 13_2$ & $3 \cdot 23_2$ & 
$3 \cdot 33_2$ & $3 \cdot 2\ol{2}_1$ & $3 \cdot 3\ol{2}_1$ & $3 \cdot 
3\ol{1}_1$ \\ $3\ol{3}$ & $0 \cdot 3\ol{3}_1$ & $\ol{3} \cdot 13_2$ & $\ol{3} 
\cdot 23_2$ & $\ol{3} \cdot 33_2$ & $\ol{3} \cdot 2\ol{2}_1$ & $\ol{3} \cdot 
3\ol{2}_1$ & $\ol{3} \cdot 3\ol{1}_1$ \\ $3\ol{2}$ & $3 \cdot 0\ol{2}_1$ & 
$\ol{2} \cdot 13_2$ & $\ol{2} \cdot 23_2$ & $\ol{2} \cdot 33_2$ & $3 \cdot 
\ol{3}\ol{2}_1$ & $3 \cdot \ol{2}\ol{2}_1$ & $\ol{2} \cdot 3\ol{1}_1$ \\ 
$3\ol{1}$ & $3 \cdot 0\ol{1}_1$ & $\ol{1} \cdot 13_2$ & $\ol{1} \cdot 23_2$ & 
$\ol{1} \cdot 33_2$ & $3 \cdot \ol{3}\ol{1}_1$ & $3 \cdot \ol{2}\ol{1}_1$ & 
$3 \cdot \ol{1}\ol{1}_1$ \\ $\ol{3}\ol{3}$ & $\ol{3} \cdot 0\ol{3}_2$ & 
$\ol{3} \cdot 1\ol{3}_2$ & $\ol{3} \cdot 2\ol{3}_2$ & $\ol{3} \cdot 
3\ol{3}_2$ & $\ol{3} \cdot \ol{3}\ol{3}_2$ & $\ol{3} \cdot \ol{3}\ol{2}_1$ & 
$\ol{3} \cdot \ol{3}\ol{1}_1$ \\ $\ol{3}\ol{2}$ & $\ol{2} \cdot 0\ol{3}_2$ & 
$\ol{2} \cdot 1\ol{3}_2$ & $\ol{2} \cdot 2\ol{3}_2$ & $\ol{2} \cdot 
3\ol{3}_2$ & $\ol{2} \cdot \ol{3}\ol{3}_2$ & $\ol{3} \cdot \ol{2}\ol{2}_1$ & 
$\ol{2} \cdot \ol{3}\ol{1}_1$ \\ $\ol{3}\ol{1}$ & $\ol{1} \cdot 0\ol{3}_2$ & 
$\ol{1} \cdot 1\ol{3}_2$ & $\ol{1} \cdot 2\ol{3}_2$ & $\ol{1} \cdot 
3\ol{3}_2$ & $\ol{1} \cdot \ol{3}\ol{3}_2$ & $\ol{3} \cdot \ol{2}\ol{1}_1$ & 
$\ol{3} \cdot \ol{1}\ol{1}_1$ \\ $\ol{2}\ol{2}$ & $\ol{2} \cdot 0\ol{2}_2$ & 
$\ol{2} \cdot 1\ol{2}_2$ & $\ol{2} \cdot 2\ol{2}_2$ & $\ol{2} \cdot 
3\ol{2}_2$ & $\ol{2} \cdot \ol{3}\ol{2}_2$ & $\ol{2} \cdot \ol{2}\ol{2}_2$ & 
$\ol{2} \cdot \ol{2}\ol{1}_1$ \\ $\ol{2}\ol{1}$ & $\ol{1} \cdot 0\ol{2}_2$ & 
$\ol{1} \cdot 1\ol{2}_2$ & $\ol{1} \cdot 2\ol{2}_2$ & $\ol{1} \cdot 
3\ol{2}_2$ & $\ol{1} \cdot \ol{3}\ol{2}_2$ & $\ol{1} \cdot \ol{2}\ol{2}_2$ & 
$\ol{2} \cdot \ol{1}\ol{1}_1$ \\ $\ol{1}\ol{1}$ & $\ol{1} \cdot 0\ol{1}_2$ & 
$\ol{1} \cdot 1\ol{1}_2$ & $\ol{1} \cdot 2\ol{1}_2$ & $\ol{1} \cdot 
3\ol{1}_2$ & $\ol{1} \cdot \ol{3}\ol{1}_2$ & $\ol{1} \cdot \ol{2}\ol{1}_2$ & 
$\ol{1} \cdot \ol{1}\ol{1}_2$ \\  \hline \end{tabular}}\end{center}\end{table}

\begin{table}[tbp]
	\caption{$B^{1,1} \otimes B^{1,1}$ for $C^{(1)}_2$}
	\label{tab:C12_11}
    \newcommand{\lw}[1]{\smash{\lower1.6ex\hbox{#1}}}
	\begin{center}
	{
	\begin{tabular}{c|cccc}	\hline
	\lw{$b$} & \multicolumn{4}{c}{$c$} \\

& $1$ & $2$ & $\ol{2}$ & $\ol{1}$ \\  \hline$1$ & 1 & 0 & 0 & 0 \\ $2$ & 1 
& 1 & 0 & 0 \\ $\ol{2}$ & 1 & 1 & 1 & 0 \\ $\ol{1}$ & 1 & 1 & 1 & 1 \\  
\hline \end{tabular}}\end{center}\end{table}

\begin{table}[tbp]
	\caption{$B^{1,2} \otimes B^{1,1}$ for $C^{(1)}_2$}
	\label{tab:C12_21}
    \newcommand{\lw}[1]{\smash{\lower1.6ex\hbox{#1}}}
	\begin{center}
	{
	\begin{tabular}{c|cccc}	\hline
	\lw{$b$} & \multicolumn{4}{c}{$c$} \\
	
& $1$ & $2$ & $\ol{2}$ & $\ol{1}$ \\  \hline$\phi $ & $1 \cdot 1\ol{1}_0$ & 
$1 \cdot 2\ol{1}_0$ & $1 \cdot \ol{2}\ol{1}_0$ & $1 \cdot \ol{1}\ol{1}_0$ \\ 
$11$ & $1 \cdot 11_1$ & $1 \cdot 12_0$ & $1 \cdot 1\ol{2}_0$ & $1 \cdot \phi 
_0$ \\ $12$ & $2 \cdot 11_1$ & $1 \cdot 22_0$ & $2 \cdot 1\ol{2}_0$ & $2 
\cdot \phi _0$ \\ $1\ol{2}$ & $\ol{2} \cdot 11_1$ & $1 \cdot 2\ol{2}_0$ & $1 
\cdot \ol{2}\ol{2}_0$ & $\ol{2} \cdot \phi _0$ \\ $1\ol{1}$ & $\ol{1} \cdot 
11_1$ & $2 \cdot 2\ol{2}_0$ & $2 \cdot \ol{2}\ol{2}_0$ & $\ol{1} \cdot \phi 
_0$ \\ $22$ & $2 \cdot 12_1$ & $2 \cdot 22_1$ & $2 \cdot 1\ol{1}_0$ & $2 
\cdot 2\ol{1}_0$ \\ $2\ol{2}$ & $\ol{2} \cdot 12_1$ & $\ol{2} \cdot 22_1$ & 
$\ol{2} \cdot 1\ol{1}_0$ & $\ol{2} \cdot 2\ol{1}_0$ \\ $2\ol{1}$ & $\ol{1} 
\cdot 12_1$ & $\ol{1} \cdot 22_1$ & $2 \cdot \ol{2}\ol{1}_0$ & $2 \cdot 
\ol{1}\ol{1}_0$ \\ $\ol{2}\ol{2}$ & $\ol{2} \cdot 1\ol{2}_1$ & $\ol{2} \cdot 
2\ol{2}_1$ & $\ol{2} \cdot \ol{2}\ol{2}_1$ & $\ol{2} \cdot \ol{2}\ol{1}_0$ \\ 
$\ol{2}\ol{1}$ & $\ol{1} \cdot 1\ol{2}_1$ & $\ol{1} \cdot 2\ol{2}_1$ & 
$\ol{1} \cdot \ol{2}\ol{2}_1$ & $\ol{2} \cdot \ol{1}\ol{1}_0$ \\ 
$\ol{1}\ol{1}$ & $\ol{1} \cdot 1\ol{1}_1$ & $\ol{1} \cdot 2\ol{1}_1$ & 
$\ol{1} \cdot \ol{2}\ol{1}_1$ & $\ol{1} \cdot \ol{1}\ol{1}_1$ \\  \hline 
\end{tabular}}\end{center}\end{table}

\begin{table}[tbp]
	\caption{$B^{1,2} \otimes B^{1,2}$ for $C^{(1)}_2$}
	\label{tab:C12_22}
    \newcommand{\lw}[1]{\smash{\lower1.6ex\hbox{#1}}}
	\begin{center}
	{
	\begin{tabular}{c|ccccccccccc}	\hline
	\lw{$b$} & \multicolumn{11}{c}{$c$} \\

 & $\phi $ & $11$ & $12$ & $1\ol{2}$ & $1\ol{1}$ & $22$ & $2\ol{2}$ & 
$2\ol{1}$ & $\ol{2}\ol{2}$ & $\ol{2}\ol{1}$ & $\ol{1}\ol{1}$ \\  \hline$\phi 
$ & 0 & 1 & 1 & 1 & 1 & 1 & 1 & 1 & 1 & 1 & 1 \\ $11$ & 1 & 2 & 1 & 1 & 0 & 0 
& 0 & 0 & 0 & 0 & 0 \\ $12$ & 1 & 2 & 1 & 1 & 0 & 1 & 0 & 0 & 0 & 0 & 0 \\ 
$1\ol{2}$ & 1 & 2 & 1 & 1 & 0 & 1 & 1 & 0 & 1 & 0 & 0 \\ $1\ol{1}$ & 1 & 2 & 
1 & 1 & 0 & 1 & 1 & 0 & 1 & 0 & 0 \\ $22$ & 1 & 2 & 2 & 1 & 1 & 2 & 0 & 1 & 0 
& 0 & 0 \\ $2\ol{2}$ & 1 & 2 & 2 & 1 & 1 & 2 & 0 & 1 & 0 & 0 & 0 \\ $2\ol{1}$ 
& 1 & 2 & 2 & 1 & 1 & 2 & 1 & 1 & 1 & 1 & 1 \\ $\ol{2}\ol{2}$ & 1 & 2 & 2 & 2 
& 1 & 2 & 2 & 1 & 2 & 1 & 0 \\ $\ol{2}\ol{1}$ & 1 & 2 & 2 & 2 & 1 & 2 & 2 & 1 
& 2 & 1 & 1 \\ $\ol{1}\ol{1}$ & 1 & 2 & 2 & 2 & 2 & 2 & 2 & 2 & 2 & 2 & 2 \\  
\hline \end{tabular}}\end{center}\end{table}

\begin{table}[tbp]
	\caption{$B^{1,2} \otimes B^{1,1}$ for $D^{(1)}_4$}
	\label{tab:D14_21}
    \newcommand{\lw}[1]{\smash{\lower1.6ex\hbox{#1}}}
	\begin{center}
    {
	\begin{tabular}{c|cccccccc}	\hline
	\lw{$b$} & \multicolumn{8}{c}{$c$} \\

& $1$ & $2$ & $3$ & $4$ & $\ol{4}$ & $\ol{3}$ & $\ol{2}$ & $\ol{1}$ \\  
\hline$11$ & $1 \cdot 11_2$ & $1 \cdot 12_1$ & $1 \cdot 13_1$ & $1 \cdot 
14_1$ & $1 \cdot 1\ol{4}_1$ & $1 \cdot 1\ol{3}_1$ & $1 \cdot 1\ol{2}_1$ & $1 
\cdot 1\ol{1}_0$ \\ $12$ & $2 \cdot 11_2$ & $1 \cdot 22_1$ & $2 \cdot 13_1$ & 
$2 \cdot 14_1$ & $2 \cdot 1\ol{4}_1$ & $2 \cdot 1\ol{3}_1$ & $2 \cdot 
1\ol{2}_1$ & $1 \cdot 2\ol{1}_0$ \\ $13$ & $3 \cdot 11_2$ & $1 \cdot 23_1$ & 
$1 \cdot 33_1$ & $3 \cdot 14_1$ & $3 \cdot 1\ol{4}_1$ & $3 \cdot 1\ol{3}_1$ & 
$3 \cdot 1\ol{2}_1$ & $1 \cdot 3\ol{1}_0$ \\ $14$ & $4 \cdot 11_2$ & $1 \cdot 
24_1$ & $1 \cdot 34_1$ & $1 \cdot 44_1$ & $4 \cdot 1\ol{4}_1$ & $4 \cdot 
1\ol{3}_1$ & $4 \cdot 1\ol{2}_1$ & $1 \cdot 4\ol{1}_0$ \\ $1\ol{4}$ & $\ol{4} 
\cdot 11_2$ & $1 \cdot 2\ol{4}_1$ & $1 \cdot 3\ol{4}_1$ & $\ol{4} \cdot 14_1$ 
& $1 \cdot \ol{4}\ol{4}_1$ & $\ol{4} \cdot 1\ol{3}_1$ & $\ol{4} \cdot 
1\ol{2}_1$ & $1 \cdot \ol{4}\ol{1}_0$ \\ $1\ol{3}$ & $\ol{3} \cdot 11_2$ & $1 
\cdot 2\ol{3}_1$ & $1 \cdot 3\ol{3}_1$ & $1 \cdot 4\ol{3}_1$ & $1 \cdot 
\ol{4}\ol{3}_1$ & $1 \cdot \ol{3}\ol{3}_1$ & $\ol{3} \cdot 1\ol{2}_1$ & $1 
\cdot \ol{3}\ol{1}_0$ \\ $1\ol{2}$ & $\ol{2} \cdot 11_2$ & $1 \cdot 
2\ol{2}_1$ & $1 \cdot 3\ol{2}_1$ & $1 \cdot 4\ol{2}_1$ & $1 \cdot 
\ol{4}\ol{2}_1$ & $1 \cdot \ol{3}\ol{2}_1$ & $1 \cdot \ol{2}\ol{2}_1$ & $1 
\cdot \ol{2}\ol{1}_0$ \\ $1\ol{1}$ & $\ol{1} \cdot 11_2$ & $2 \cdot 
2\ol{2}_1$ & $2 \cdot 3\ol{2}_1$ & $2 \cdot 4\ol{2}_1$ & $2 \cdot 
\ol{4}\ol{2}_1$ & $2 \cdot \ol{3}\ol{2}_1$ & $2 \cdot \ol{2}\ol{2}_1$ & $1 
\cdot \ol{1}\ol{1}_0$ \\ $22$ & $2 \cdot 12_2$ & $2 \cdot 22_2$ & $2 \cdot 
23_1$ & $2 \cdot 24_1$ & $2 \cdot 2\ol{4}_1$ & $2 \cdot 2\ol{3}_1$ & $2 \cdot 
1\ol{1}_1$ & $2 \cdot 2\ol{1}_1$ \\ $23$ & $3 \cdot 12_2$ & $3 \cdot 22_2$ & 
$2 \cdot 33_1$ & $3 \cdot 24_1$ & $3 \cdot 2\ol{4}_1$ & $3 \cdot 2\ol{3}_1$ & 
$3 \cdot 1\ol{1}_1$ & $3 \cdot 2\ol{1}_1$ \\ $24$ & $4 \cdot 12_2$ & $4 \cdot 
22_2$ & $2 \cdot 34_1$ & $2 \cdot 44_1$ & $4 \cdot 2\ol{4}_1$ & $4 \cdot 
2\ol{3}_1$ & $4 \cdot 1\ol{1}_1$ & $4 \cdot 2\ol{1}_1$ \\ $2\ol{4}$ & $\ol{4} 
\cdot 12_2$ & $\ol{4} \cdot 22_2$ & $2 \cdot 3\ol{4}_1$ & $\ol{4} \cdot 24_1$ 
& $2 \cdot \ol{4}\ol{4}_1$ & $\ol{4} \cdot 2\ol{3}_1$ & $\ol{4} \cdot 
1\ol{1}_1$ & $\ol{4} \cdot 2\ol{1}_1$ \\ $2\ol{3}$ & $\ol{3} \cdot 12_2$ & 
$\ol{3} \cdot 22_2$ & $2 \cdot 3\ol{3}_1$ & $2 \cdot 4\ol{3}_1$ & $2 \cdot 
\ol{4}\ol{3}_1$ & $2 \cdot \ol{3}\ol{3}_1$ & $\ol{3} \cdot 1\ol{1}_1$ & 
$\ol{3} \cdot 2\ol{1}_1$ \\ $2\ol{2}$ & $\ol{2} \cdot 12_2$ & $\ol{2} \cdot 
22_2$ & $3 \cdot 3\ol{3}_1$ & $3 \cdot 4\ol{3}_1$ & $3 \cdot \ol{4}\ol{3}_1$ 
& $3 \cdot \ol{3}\ol{3}_1$ & $\ol{2} \cdot 1\ol{1}_1$ & $\ol{2} \cdot 
2\ol{1}_1$ \\ $2\ol{1}$ & $\ol{1} \cdot 12_2$ & $\ol{1} \cdot 22_2$ & $2 
\cdot 3\ol{1}_1$ & $2 \cdot 4\ol{1}_1$ & $2 \cdot \ol{4}\ol{1}_1$ & $2 \cdot 
\ol{3}\ol{1}_1$ & $2 \cdot \ol{2}\ol{1}_1$ & $2 \cdot \ol{1}\ol{1}_1$ \\ $33$ 
& $3 \cdot 13_2$ & $3 \cdot 23_2$ & $3 \cdot 33_2$ & $3 \cdot 34_1$ & $3 
\cdot 3\ol{4}_1$ & $3 \cdot 2\ol{2}_1$ & $3 \cdot 3\ol{2}_1$ & $3 \cdot 
3\ol{1}_1$ \\ $34$ & $4 \cdot 13_2$ & $4 \cdot 23_2$ & $4 \cdot 33_2$ & $3 
\cdot 44_1$ & $4 \cdot 3\ol{4}_1$ & $4 \cdot 2\ol{2}_1$ & $4 \cdot 3\ol{2}_1$ 
& $4 \cdot 3\ol{1}_1$ \\ $3\ol{4}$ & $\ol{4} \cdot 13_2$ & $\ol{4} \cdot 
23_2$ & $\ol{4} \cdot 33_2$ & $\ol{4} \cdot 34_1$ & $3 \cdot \ol{4}\ol{4}_1$ 
& $\ol{4} \cdot 2\ol{2}_1$ & $\ol{4} \cdot 3\ol{2}_1$ & $\ol{4} \cdot 
3\ol{1}_1$ \\ $3\ol{3}$ & $\ol{3} \cdot 13_2$ & $\ol{3} \cdot 23_2$ & $\ol{3} 
\cdot 33_2$ & $\ol{4} \cdot 44_1$ & $4 \cdot \ol{4}\ol{4}_1$ & $\ol{3} \cdot 
2\ol{2}_1$ & $\ol{3} \cdot 3\ol{2}_1$ & $\ol{3} \cdot 3\ol{1}_1$ \\ $3\ol{2}$ 
& $\ol{2} \cdot 13_2$ & $\ol{2} \cdot 23_2$ & $\ol{2} \cdot 33_2$ & $3 \cdot 
4\ol{2}_1$ & $3 \cdot \ol{4}\ol{2}_1$ & $3 \cdot \ol{3}\ol{2}_1$ & $3 \cdot 
\ol{2}\ol{2}_1$ & $\ol{2} \cdot 3\ol{1}_1$ \\ $3\ol{1}$ & $\ol{1} \cdot 13_2$ 
& $\ol{1} \cdot 23_2$ & $\ol{1} \cdot 33_2$ & $3 \cdot 4\ol{1}_1$ & $3 \cdot 
\ol{4}\ol{1}_1$ & $3 \cdot \ol{3}\ol{1}_1$ & $3 \cdot \ol{2}\ol{1}_1$ & $3 
\cdot \ol{1}\ol{1}_1$ \\ $44$ & $4 \cdot 14_2$ & $4 \cdot 24_2$ & $4 \cdot 
34_2$ & $4 \cdot 44_2$ & $4 \cdot 3\ol{3}_1$ & $4 \cdot 4\ol{3}_1$ & $4 \cdot 
4\ol{2}_1$ & $4 \cdot 4\ol{1}_1$ \\ $4\ol{3}$ & $\ol{3} \cdot 14_2$ & $\ol{3} 
\cdot 24_2$ & $\ol{3} \cdot 34_2$ & $\ol{3} \cdot 44_2$ & $4 \cdot 
\ol{4}\ol{3}_1$ & $4 \cdot \ol{3}\ol{3}_1$ & $\ol{3} \cdot 4\ol{2}_1$ & 
$\ol{3} \cdot 4\ol{1}_1$ \\ $4\ol{2}$ & $\ol{2} \cdot 14_2$ & $\ol{2} \cdot 
24_2$ & $\ol{2} \cdot 34_2$ & $\ol{2} \cdot 44_2$ & $4 \cdot \ol{4}\ol{2}_1$ 
& $4 \cdot \ol{3}\ol{2}_1$ & $4 \cdot \ol{2}\ol{2}_1$ & $\ol{2} \cdot 
4\ol{1}_1$ \\ $4\ol{1}$ & $\ol{1} \cdot 14_2$ & $\ol{1} \cdot 24_2$ & $\ol{1} 
\cdot 34_2$ & $\ol{1} \cdot 44_2$ & $4 \cdot \ol{4}\ol{1}_1$ & $4 \cdot 
\ol{3}\ol{1}_1$ & $4 \cdot \ol{2}\ol{1}_1$ & $4 \cdot \ol{1}\ol{1}_1$ \\ 
$\ol{4}\ol{4}$ & $\ol{4} \cdot 1\ol{4}_2$ & $\ol{4} \cdot 2\ol{4}_2$ & 
$\ol{4} \cdot 3\ol{4}_2$ & $\ol{4} \cdot 3\ol{3}_1$ & $\ol{4} \cdot 
\ol{4}\ol{4}_2$ & $\ol{4} \cdot \ol{4}\ol{3}_1$ & $\ol{4} \cdot 
\ol{4}\ol{2}_1$ & $\ol{4} \cdot \ol{4}\ol{1}_1$ \\ $\ol{4}\ol{3}$ & $\ol{3} 
\cdot 1\ol{4}_2$ & $\ol{3} \cdot 2\ol{4}_2$ & $\ol{3} \cdot 3\ol{4}_2$ & 
$\ol{4} \cdot 4\ol{3}_1$ & $\ol{3} \cdot \ol{4}\ol{4}_2$ & $\ol{4} \cdot 
\ol{3}\ol{3}_1$ & $\ol{3} \cdot \ol{4}\ol{2}_1$ & $\ol{3} \cdot 
\ol{4}\ol{1}_1$ \\ $\ol{4}\ol{2}$ & $\ol{2} \cdot 1\ol{4}_2$ & $\ol{2} \cdot 
2\ol{4}_2$ & $\ol{2} \cdot 3\ol{4}_2$ & $\ol{4} \cdot 4\ol{2}_1$ & $\ol{2} 
\cdot \ol{4}\ol{4}_2$ & $\ol{4} \cdot \ol{3}\ol{2}_1$ & $\ol{4} \cdot 
\ol{2}\ol{2}_1$ & $\ol{2} \cdot \ol{4}\ol{1}_1$ \\ $\ol{4}\ol{1}$ & $\ol{1} 
\cdot 1\ol{4}_2$ & $\ol{1} \cdot 2\ol{4}_2$ & $\ol{1} \cdot 3\ol{4}_2$ & 
$\ol{4} \cdot 4\ol{1}_1$ & $\ol{1} \cdot \ol{4}\ol{4}_2$ & $\ol{4} \cdot 
\ol{3}\ol{1}_1$ & $\ol{4} \cdot \ol{2}\ol{1}_1$ & $\ol{4} \cdot 
\ol{1}\ol{1}_1$ \\ $\ol{3}\ol{3}$ & $\ol{3} \cdot 1\ol{3}_2$ & $\ol{3} \cdot 
2\ol{3}_2$ & $\ol{3} \cdot 3\ol{3}_2$ & $\ol{3} \cdot 4\ol{3}_2$ & $\ol{3} 
\cdot \ol{4}\ol{3}_2$ & $\ol{3} \cdot \ol{3}\ol{3}_2$ & $\ol{3} \cdot 
\ol{3}\ol{2}_1$ & $\ol{3} \cdot \ol{3}\ol{1}_1$ \\ $\ol{3}\ol{2}$ & $\ol{2} 
\cdot 1\ol{3}_2$ & $\ol{2} \cdot 2\ol{3}_2$ & $\ol{2} \cdot 3\ol{3}_2$ & 
$\ol{2} \cdot 4\ol{3}_2$ & $\ol{2} \cdot \ol{4}\ol{3}_2$ & $\ol{2} \cdot 
\ol{3}\ol{3}_2$ & $\ol{3} \cdot \ol{2}\ol{2}_1$ & $\ol{2} \cdot 
\ol{3}\ol{1}_1$ \\ $\ol{3}\ol{1}$ & $\ol{1} \cdot 1\ol{3}_2$ & $\ol{1} \cdot 
2\ol{3}_2$ & $\ol{1} \cdot 3\ol{3}_2$ & $\ol{1} \cdot 4\ol{3}_2$ & $\ol{1} 
\cdot \ol{4}\ol{3}_2$ & $\ol{1} \cdot \ol{3}\ol{3}_2$ & $\ol{3} \cdot 
\ol{2}\ol{1}_1$ & $\ol{3} \cdot \ol{1}\ol{1}_1$ \\ $\ol{2}\ol{2}$ & $\ol{2} 
\cdot 1\ol{2}_2$ & $\ol{2} \cdot 2\ol{2}_2$ & $\ol{2} \cdot 3\ol{2}_2$ & 
$\ol{2} \cdot 4\ol{2}_2$ & $\ol{2} \cdot \ol{4}\ol{2}_2$ & $\ol{2} \cdot 
\ol{3}\ol{2}_2$ & $\ol{2} \cdot \ol{2}\ol{2}_2$ & $\ol{2} \cdot 
\ol{2}\ol{1}_1$ \\ $\ol{2}\ol{1}$ & $\ol{1} \cdot 1\ol{2}_2$ & $\ol{1} \cdot 
2\ol{2}_2$ & $\ol{1} \cdot 3\ol{2}_2$ & $\ol{1} \cdot 4\ol{2}_2$ & $\ol{1} 
\cdot \ol{4}\ol{2}_2$ & $\ol{1} \cdot \ol{3}\ol{2}_2$ & $\ol{1} \cdot 
\ol{2}\ol{2}_2$ & $\ol{2} \cdot \ol{1}\ol{1}_1$ \\ $\ol{1}\ol{1}$ & $\ol{1} 
\cdot 1\ol{1}_2$ & $\ol{1} \cdot 2\ol{1}_2$ & $\ol{1} \cdot 3\ol{1}_2$ & 
$\ol{1} \cdot 4\ol{1}_2$ & $\ol{1} \cdot \ol{4}\ol{1}_2$ & $\ol{1} \cdot 
\ol{3}\ol{1}_2$ & $\ol{1} \cdot \ol{2}\ol{1}_2$ & $\ol{1} \cdot 
\ol{1}\ol{1}_2$ \\  \hline \end{tabular}}\end{center}\end{table}

\begin{table}[tbp]
	\caption{$B^{1,1} \otimes B^{1,1}$ for $A^{(2)}_3$}
	\label{tab:A23_11}
    \newcommand{\lw}[1]{\smash{\lower1.6ex\hbox{#1}}}
	\begin{center}
	\begin{tabular}{c|cccc}	\hline
	\lw{$b$} & \multicolumn{4}{c}{$c$} \\

& $1$ & $2$ & $\ol{2}$ & $\ol{1}$ \\  \hline$1$ & 2 & 1 & 1 & 0 \\ $2$ & 2 
& 2 & 1 & 1 \\ $\ol{2}$ & 2 & 2 & 2 & 1 \\ $\ol{1}$ & 2 & 2 & 2 & 2 \\  
\hline \end{tabular}\end{center}\end{table}

\begin{table}[tbp]
	\caption{$B^{1,2} \otimes B^{1,1}$ for $A^{(2)}_3$}
	\label{tab:A23_21}
    \newcommand{\lw}[1]{\smash{\lower1.6ex\hbox{#1}}}
	\begin{center}
	\begin{tabular}{c|cccc}	\hline
	\lw{$b$} & \multicolumn{4}{c}{$c$} \\

& $1$ & $2$ & $\ol{2}$ & $\ol{1}$ \\  \hline$11$ & $1 \cdot 11_2$ & $1 
\cdot 12_1$ & $1 \cdot 1\ol{2}_1$ & $1 \cdot 1\ol{1}_0$ \\ $12$ & $2 \cdot 
11_2$ & $1 \cdot 22_1$ & $2 \cdot 1\ol{2}_1$ & $1 \cdot 2\ol{1}_0$ \\ 
$1\ol{2}$ & $\ol{2} \cdot 11_2$ & $1 \cdot 2\ol{2}_1$ & $1 \cdot 
\ol{2}\ol{2}_1$ & $1 \cdot \ol{2}\ol{1}_0$ \\ $1\ol{1}$ & $\ol{1} \cdot 11_2$ 
& $2 \cdot 2\ol{2}_1$ & $2 \cdot \ol{2}\ol{2}_1$ & $1 \cdot \ol{1}\ol{1}_0$ 
\\ $22$ & $2 \cdot 12_2$ & $2 \cdot 22_2$ & $2 \cdot 1\ol{1}_1$ & $2 \cdot 
2\ol{1}_1$ \\ $2\ol{2}$ & $\ol{2} \cdot 12_2$ & $\ol{2} \cdot 22_2$ & $\ol{2} 
\cdot 1\ol{1}_1$ & $\ol{2} \cdot 2\ol{1}_1$ \\ $2\ol{1}$ & $\ol{1} \cdot 
12_2$ & $\ol{1} \cdot 22_2$ & $2 \cdot \ol{2}\ol{1}_1$ & $2 \cdot 
\ol{1}\ol{1}_1$ \\ $\ol{2}\ol{2}$ & $\ol{2} \cdot 1\ol{2}_2$ & $\ol{2} \cdot 
2\ol{2}_2$ & $\ol{2} \cdot \ol{2}\ol{2}_2$ & $\ol{2} \cdot \ol{2}\ol{1}_1$ \\ 
$\ol{2}\ol{1}$ & $\ol{1} \cdot 1\ol{2}_2$ & $\ol{1} \cdot 2\ol{2}_2$ & 
$\ol{1} \cdot \ol{2}\ol{2}_2$ & $\ol{2} \cdot \ol{1}\ol{1}_1$ \\ 
$\ol{1}\ol{1}$ & $\ol{1} \cdot 1\ol{1}_2$ & $\ol{1} \cdot 2\ol{1}_2$ & 
$\ol{1} \cdot \ol{2}\ol{1}_2$ & $\ol{1} \cdot \ol{1}\ol{1}_2$ \\  \hline 
\end{tabular}\end{center}\end{table}

\begin{table}[tbp]
	\caption{$B^{1,2} \otimes B^{1,2}$ for $A^{(2)}_3$}
	\label{tab:A23_22}
    \newcommand{\lw}[1]{\smash{\lower1.6ex\hbox{#1}}}
	\begin{center}
	\begin{tabular}{c|cccccccccc}	\hline
	\lw{$b$} & \multicolumn{10}{c}{$c$} \\

& $11$ & $12$ & $1\ol{2}$ & $1\ol{1}$ & $22$ & $2\ol{2}$ & $2\ol{1}$ & 
$\ol{2}\ol{2}$ & $\ol{2}\ol{1}$ & $\ol{1}\ol{1}$ \\  \hline$11$ & 4 & 3 & 3 & 
2 & 2 & 2 & 1 & 2 & 1 & 0 \\ $12$ & 4 & 3 & 3 & 2 & 3 & 2 & 2 & 2 & 1 & 1 \\ 
$1\ol{2}$ & 4 & 3 & 3 & 2 & 3 & 3 & 2 & 3 & 2 & 1 \\ $1\ol{1}$ & 4 & 3 & 3 & 
2 & 3 & 3 & 2 & 3 & 2 & 2 \\ $22$ & 4 & 4 & 3 & 3 & 4 & 2 & 3 & 2 & 2 & 2 \\ 
$2\ol{2}$ & 4 & 4 & 3 & 3 & 4 & 2 & 3 & 2 & 2 & 2 \\ $2\ol{1}$ & 4 & 4 & 3 & 
3 & 4 & 3 & 3 & 3 & 3 & 3 \\ $\ol{2}\ol{2}$ & 4 & 4 & 4 & 3 & 4 & 4 & 3 & 4 & 
3 & 2 \\ $\ol{2}\ol{1}$ & 4 & 4 & 4 & 3 & 4 & 4 & 3 & 4 & 3 & 3 \\ 
$\ol{1}\ol{1}$ & 4 & 4 & 4 & 4 & 4 & 4 & 4 & 4 & 4 & 4 \\  \hline 
\end{tabular}\end{center}\end{table}

\begin{table}[tbp]
	\caption{$B^{1,1} \otimes B^{1,1}$ for $A^{(2)}_4$}
	\label{tab:A24_11}
    \newcommand{\lw}[1]{\smash{\lower1.6ex\hbox{#1}}}
	\begin{center}
	\begin{tabular}{c|ccccc}	\hline
	\lw{$b$} & \multicolumn{5}{c}{$c$} \\

& $\phi $ & $1$ & $2$ & $\ol{2}$ & $\ol{1}$ \\  \hline$\phi $ & 0 
& 1 & 1 & 1 & 1 \\ $1$ & 1 & 2 & 0 & 0 & 0 \\ $2$ & 1 & 2 & 2 & 0 & 0 \\ 
$\ol{2}$ & 1 & 2 & 2 & 2 & 0 \\ $\ol{1}$ & 1 & 2 & 2 & 2 & 2 \\  \hline 
\end{tabular}\end{center}\end{table}

\begin{table}[tbp]
	\caption{$B^{1,2} \otimes B^{1,1}$ for $A^{(2)}_4$}
	\label{tab:A24_21}
    \newcommand{\lw}[1]{\smash{\lower1.6ex\hbox{#1}}}
	\begin{center}
	\begin{tabular}{c|ccccc}	\hline
	\lw{$b$} & \multicolumn{5}{c}{$c$} \\
	
& $\phi $ & $1$ & $2$ & $\ol{2}$ & $\ol{1}$ \\  \hline$\phi $ & $\phi  
\cdot \phi _0$ & $1 \cdot 1\ol{1}_0$ & $1 \cdot 2\ol{1}_0$ & $1 \cdot 
\ol{2}\ol{1}_0$ & $1 \cdot \ol{1}\ol{1}_0$ \\ $1$ & $\phi  \cdot 1_0$ & $\phi 
 \cdot 11_1$ & $1 \cdot 2_0$ & $1 \cdot \ol{2}_0$ & $1 \cdot \ol{1}_0$ \\ $2$ 
& $\phi  \cdot 2_0$ & $\phi  \cdot 12_1$ & $\phi  \cdot 22_1$ & $2 \cdot 
\ol{2}_0$ & $2 \cdot \ol{1}_0$ \\ $\ol{2}$ & $\phi  \cdot \ol{2}_0$ & $\phi  
\cdot 1\ol{2}_1$ & $\phi  \cdot 2\ol{2}_1$ & $\phi  \cdot \ol{2}\ol{2}_1$ & 
$\ol{2} \cdot \ol{1}_0$ \\ $\ol{1}$ & $\phi  \cdot \ol{1}_0$ & $\phi  \cdot 
1\ol{1}_1$ & $\phi  \cdot 2\ol{1}_1$ & $\phi  \cdot \ol{2}\ol{1}_1$ & $\phi  
\cdot \ol{1}\ol{1}_1$ \\ $11$ & $1 \cdot 1_1$ & $1 \cdot 11_2$ & $1 \cdot 
12_0$ & $1 \cdot 1\ol{2}_0$ & $1 \cdot \phi _0$ \\ $12$ & $2 \cdot 1_1$ & $2 
\cdot 11_2$ & $1 \cdot 22_0$ & $2 \cdot 1\ol{2}_0$ & $2 \cdot \phi _0$ \\ 
$1\ol{2}$ & $\ol{2} \cdot 1_1$ & $\ol{2} \cdot 11_2$ & $1 \cdot 2\ol{2}_0$ & 
$1 \cdot \ol{2}\ol{2}_0$ & $\ol{2} \cdot \phi _0$ \\ $1\ol{1}$ & $\ol{1} 
\cdot 1_1$ & $\ol{1} \cdot 11_2$ & $2 \cdot 2\ol{2}_0$ & $2 \cdot 
\ol{2}\ol{2}_0$ & $\ol{1} \cdot \phi _0$ \\ $22$ & $2 \cdot 2_1$ & $2 \cdot 
12_2$ & $2 \cdot 22_2$ & $2 \cdot 1\ol{1}_0$ & $2 \cdot 2\ol{1}_0$ \\ 
$2\ol{2}$ & $\ol{2} \cdot 2_1$ & $\ol{2} \cdot 12_2$ & $\ol{2} \cdot 22_2$ & 
$\ol{2} \cdot 1\ol{1}_0$ & $\ol{2} \cdot 2\ol{1}_0$ \\ $2\ol{1}$ & $\ol{1} 
\cdot 2_1$ & $\ol{1} \cdot 12_2$ & $\ol{1} \cdot 22_2$ & $2 \cdot 
\ol{2}\ol{1}_0$ & $2 \cdot \ol{1}\ol{1}_0$ \\ $\ol{2}\ol{2}$ & $\ol{2} \cdot 
\ol{2}_1$ & $\ol{2} \cdot 1\ol{2}_2$ & $\ol{2} \cdot 2\ol{2}_2$ & $\ol{2} 
\cdot \ol{2}\ol{2}_2$ & $\ol{2} \cdot \ol{2}\ol{1}_0$ \\ $\ol{2}\ol{1}$ & 
$\ol{1} \cdot \ol{2}_1$ & $\ol{1} \cdot 1\ol{2}_2$ & $\ol{1} \cdot 2\ol{2}_2$ 
& $\ol{1} \cdot \ol{2}\ol{2}_2$ & $\ol{2} \cdot \ol{1}\ol{1}_0$ \\ 
$\ol{1}\ol{1}$ & $\ol{1} \cdot \ol{1}_1$ & $\ol{1} \cdot 1\ol{1}_2$ & $\ol{1} 
\cdot 2\ol{1}_2$ & $\ol{1} \cdot \ol{2}\ol{1}_2$ & $\ol{1} \cdot 
\ol{1}\ol{1}_2$ \\  \hline \end{tabular}\end{center}\end{table}

\begin{table}[tbp]
	\caption{$B^{1,2} \otimes B^{1,2}$ for $A^{(2)}_4$}
	\label{tab:A24_22}
    \newcommand{\lw}[1]{\smash{\lower1.6ex\hbox{#1}}}
	\begin{center}
	\begin{tabular}{c|ccccccccccccccc}	\hline
	\lw{$b$} & \multicolumn{15}{c}{$c$} \\
	
& $\phi $ & $1$ & $2$ & $\ol{2}$ & $\ol{1}$ & $11$ & $12$ & $1\ol{2}$ & 
$1\ol{1}$ & $22$ & $2\ol{2}$ & $2\ol{1}$ & $\ol{2}\ol{2}$ & $\ol{2}\ol{1}$ & 
$\ol{1}\ol{1}$ \\  \hline$\phi $ & 0 & 1 & 1 & 1 & 1 & 2 & 2 & 2 & 2 & 2 & 2 
& 2 & 2 & 2 & 2 \\ $1$ & 1 & 2 & 0 & 0 & 0 & 3 & 1 & 1 & 1 & 1 & 1 & 1 & 1 & 
1 & 1 \\ $2$ & 1 & 2 & 2 & 0 & 0 & 3 & 3 & 1 & 1 & 3 & 1 & 1 & 1 & 1 & 1 \\ 
$\ol{2}$ & 1 & 2 & 2 & 2 & 0 & 3 & 3 & 3 & 1 & 3 & 3 & 1 & 3 & 1 & 1 \\ 
$\ol{1}$ & 1 & 2 & 2 & 2 & 2 & 3 & 3 & 3 & 3 & 3 & 3 & 3 & 3 & 3 & 3 \\ $11$ 
& 2 & 3 & 1 & 1 & 1 & 4 & 2 & 2 & 0 & 0 & 0 & 0 & 0 & 0 & 0 \\ $12$ & 2 & 3 & 
1 & 1 & 1 & 4 & 2 & 2 & 0 & 2 & 0 & 0 & 0 & 0 & 0 \\ $1\ol{2}$ & 2 & 3 & 1 & 
1 & 1 & 4 & 2 & 2 & 0 & 2 & 2 & 0 & 2 & 0 & 0 \\ $1\ol{1}$ & 2 & 3 & 1 & 1 & 
1 & 4 & 2 & 2 & 0 & 2 & 2 & 0 & 2 & 0 & 0 \\ $22$ & 2 & 3 & 3 & 1 & 1 & 4 & 4 
& 2 & 2 & 4 & 0 & 2 & 0 & 0 & 0 \\ $2\ol{2}$ & 2 & 3 & 3 & 1 & 1 & 4 & 4 & 2 
& 2 & 4 & 0 & 2 & 0 & 0 & 0 \\ $2\ol{1}$ & 2 & 3 & 3 & 1 & 1 & 4 & 4 & 2 & 2 
& 4 & 2 & 2 & 2 & 2 & 2 \\ $\ol{2}\ol{2}$ & 2 & 3 & 3 & 3 & 1 & 4 & 4 & 4 & 2 
& 4 & 4 & 2 & 4 & 2 & 0 \\ $\ol{2}\ol{1}$ & 2 & 3 & 3 & 3 & 1 & 4 & 4 & 4 & 2 
& 4 & 4 & 2 & 4 & 2 & 2 \\ $\ol{1}\ol{1}$ & 2 & 3 & 3 & 3 & 3 & 4 & 4 & 4 & 4 
& 4 & 4 & 4 & 4 & 4 & 4 \\  \hline \end{tabular}\end{center}\end{table}

\begin{table}[tbp]
	\caption{$B^{1,1} \otimes B^{1,1}$ for $D^{(2)}_3$}
	\label{tab:D23_11}
    \newcommand{\lw}[1]{\smash{\lower1.6ex\hbox{#1}}}
	\begin{center}
	\begin{tabular}{c|cccccc}	\hline
	\lw{$b$} & \multicolumn{6}{c}{$c$} \\
	
& $0$ & $1$ & $2$ & $\ol{2}$ & $\ol{1}$ & $\phi $ \\  \hline$0$ & 0 & 2 & 2 
& 0 & 0 & 1 \\ $1$ & 0 & 2 & 0 & 0 & 0 & 1 \\ $2$ & 0 & 2 & 2 & 0 & 0 & 1 \\ 
$\ol{2}$ & 2 & 2 & 2 & 2 & 0 & 1 \\ $\ol{1}$ & 2 & 2 & 2 & 2 & 2 & 1 \\ $\phi 
$ & 1 & 1 & 1 & 1 & 1 & 0 \\  \hline \end{tabular}\end{center}\end{table}

\begin{table}[tbp]
	\caption{$B^{1,2} \otimes B^{1,1}$ for $D^{(2)}_3$}
	\label{tab:D23_21}
    \newcommand{\lw}[1]{\smash{\lower1.6ex\hbox{#1}}}
	\begin{center}
	\begin{tabular}{c|cccccc}	\hline
	\lw{$b$} & \multicolumn{6}{c}{$c$} \\
	
& $0$ & $1$ & $2$ & $\ol{2}$ & $\ol{1}$ & $\phi $ \\  \hline$10$ & $0 \cdot 
10_0$ & $0 \cdot 11_2$ & $1 \cdot 20_0$ & $0 \cdot 1\ol{2}_0$ & $0 \cdot \phi 
_0$ & $0 \cdot 1_1$ \\ $20$ & $0 \cdot 20_0$ & $0 \cdot 12_2$ & $0 \cdot 
22_2$ & $0 \cdot 1\ol{1}_0$ & $0 \cdot 2\ol{1}_0$ & $0 \cdot 2_1$ \\ 
$0\ol{2}$ & $0 \cdot 0\ol{2}_0$ & $\ol{2} \cdot 10_2$ & $\ol{2} \cdot 20_2$ & 
$0 \cdot \ol{2}\ol{2}_0$ & $\ol{2} \cdot 0\ol{1}_0$ & $\ol{2} \cdot 0_1$ \\ 
$0\ol{1}$ & $0 \cdot 0\ol{1}_0$ & $\ol{1} \cdot 10_2$ & $\ol{1} \cdot 20_2$ & 
$0 \cdot \ol{2}\ol{1}_0$ & $0 \cdot \ol{1}\ol{1}_0$ & $\ol{1} \cdot 0_1$ \\ 
$11$ & $1 \cdot 10_0$ & $1 \cdot 11_2$ & $1 \cdot 12_0$ & $1 \cdot 1\ol{2}_0$ 
& $1 \cdot \phi _0$ & $1 \cdot 1_1$ \\ $12$ & $2 \cdot 10_0$ & $2 \cdot 11_2$ 
& $1 \cdot 22_0$ & $2 \cdot 1\ol{2}_0$ & $2 \cdot \phi _0$ & $2 \cdot 1_1$ \\ 
$1\ol{2}$ & $1 \cdot 0\ol{2}_0$ & $\ol{2} \cdot 11_2$ & $1 \cdot 2\ol{2}_0$ & 
$1 \cdot \ol{2}\ol{2}_0$ & $\ol{2} \cdot \phi _0$ & $\ol{2} \cdot 1_1$ \\ 
$1\ol{1}$ & $2 \cdot 0\ol{2}_0$ & $\ol{1} \cdot 11_2$ & $2 \cdot 2\ol{2}_0$ & 
$2 \cdot \ol{2}\ol{2}_0$ & $\ol{1} \cdot \phi _0$ & $\ol{1} \cdot 1_1$ \\ 
$22$ & $2 \cdot 20_0$ & $2 \cdot 12_2$ & $2 \cdot 22_2$ & $2 \cdot 1\ol{1}_0$ 
& $2 \cdot 2\ol{1}_0$ & $2 \cdot 2_1$ \\ $2\ol{2}$ & $0 \cdot 2\ol{2}_0$ & 
$\ol{2} \cdot 12_2$ & $\ol{2} \cdot 22_2$ & $\ol{2} \cdot 1\ol{1}_0$ & 
$\ol{2} \cdot 2\ol{1}_0$ & $\ol{2} \cdot 2_1$ \\ $2\ol{1}$ & $2 \cdot 
0\ol{1}_0$ & $\ol{1} \cdot 12_2$ & $\ol{1} \cdot 22_2$ & $2 \cdot 
\ol{2}\ol{1}_0$ & $2 \cdot \ol{1}\ol{1}_0$ & $\ol{1} \cdot 2_1$ \\ 
$\ol{2}\ol{2}$ & $\ol{2} \cdot 0\ol{2}_2$ & $\ol{2} \cdot 1\ol{2}_2$ & 
$\ol{2} \cdot 2\ol{2}_2$ & $\ol{2} \cdot \ol{2}\ol{2}_2$ & $\ol{2} \cdot 
\ol{2}\ol{1}_0$ & $\ol{2} \cdot \ol{2}_1$ \\ $\ol{2}\ol{1}$ & $\ol{1} \cdot 
0\ol{2}_2$ & $\ol{1} \cdot 1\ol{2}_2$ & $\ol{1} \cdot 2\ol{2}_2$ & $\ol{1} 
\cdot \ol{2}\ol{2}_2$ & $\ol{2} \cdot \ol{1}\ol{1}_0$ & $\ol{1} \cdot 
\ol{2}_1$ \\ $\ol{1}\ol{1}$ & $\ol{1} \cdot 0\ol{1}_2$ & $\ol{1} \cdot 
1\ol{1}_2$ & $\ol{1} \cdot 2\ol{1}_2$ & $\ol{1} \cdot \ol{2}\ol{1}_2$ & 
$\ol{1} \cdot \ol{1}\ol{1}_2$ & $\ol{1} \cdot \ol{1}_1$ \\ $0$ & $0 \cdot 
0_0$ & $\phi  \cdot 10_1$ & $\phi  \cdot 20_1$ & $0 \cdot \ol{2}_0$ & $0 
\cdot \ol{1}_0$ & $\phi  \cdot 0_0$ \\ $1$ & $1 \cdot 0_0$ & $\phi  \cdot 
11_1$ & $1 \cdot 2_0$ & $1 \cdot \ol{2}_0$ & $1 \cdot \ol{1}_0$ & $\phi  
\cdot 1_0$ \\ $2$ & $2 \cdot 0_0$ & $\phi  \cdot 12_1$ & $\phi  \cdot 22_1$ & 
$2 \cdot \ol{2}_0$ & $2 \cdot \ol{1}_0$ & $\phi  \cdot 2_0$ \\ $\ol{2}$ & 
$\phi  \cdot 0\ol{2}_1$ & $\phi  \cdot 1\ol{2}_1$ & $\phi  \cdot 2\ol{2}_1$ & 
$\phi  \cdot \ol{2}\ol{2}_1$ & $\ol{2} \cdot \ol{1}_0$ & $\phi  \cdot 
\ol{2}_0$ \\ $\ol{1}$ & $\phi  \cdot 0\ol{1}_1$ & $\phi  \cdot 1\ol{1}_1$ & 
$\phi  \cdot 2\ol{1}_1$ & $\phi  \cdot \ol{2}\ol{1}_1$ & $\phi  \cdot 
\ol{1}\ol{1}_1$ & $\phi  \cdot \ol{1}_0$ \\ $\phi $ & $1 \cdot 0\ol{1}_0$ & 
$1 \cdot 1\ol{1}_0$ & $1 \cdot 2\ol{1}_0$ & $1 \cdot \ol{2}\ol{1}_0$ & $1 
\cdot \ol{1}\ol{1}_0$ & $\phi  \cdot \phi _0$ \\  \hline 
\end{tabular}\end{center}\end{table}

\begin{table}[tbp]
	\caption{$B^{1,2} \otimes B^{1,2}$ for $D^{(2)}_3$}
	\label{tab:D23_22}
    \newcommand{\lw}[1]{\smash{\lower1.6ex\hbox{#1}}}
	\begin{center}
    {
	\begin{tabular}{c|cccccccccccccccccccc}	\hline
	\lw{$b$} & \multicolumn{20}{c}{$c$} \\
	
& $10$ & $20$ & $0\ol{2}$ & $0\ol{1}$ & $11$ & $12$ & $1\ol{2}$ & $1\ol{1}$ 
& $22$ & $2\ol{2}$ & $2\ol{1}$ & $\ol{2}\ol{2}$ & $\ol{2}\ol{1}$ & 
$\ol{1}\ol{1}$ & $0$ & $1$ & $2$ & $\ol{2}$ & $\ol{1}$ & $\phi $ \\  
\hline$10$ & 2 & 0 & 0 & 0 & 4 & 2 & 2 & 0 & 2 & 0 & 0 & 0 & 0 & 0 & 1 & 3 & 
1 & 1 & 1 & 2 \\ $20$ & 2 & 2 & 0 & 0 & 4 & 4 & 2 & 2 & 4 & 0 & 2 & 0 & 0 & 0 
& 1 & 3 & 3 & 1 & 1 & 2 \\ $0\ol{2}$ & 2 & 2 & 2 & 0 & 4 & 4 & 2 & 2 & 4 & 2 
& 2 & 2 & 0 & 0 & 1 & 3 & 3 & 1 & 1 & 2 \\ $0\ol{1}$ & 2 & 2 & 2 & 2 & 4 & 4 
& 2 & 2 & 4 & 2 & 2 & 2 & 2 & 2 & 1 & 3 & 3 & 1 & 1 & 2 \\ $11$ & 2 & 0 & 0 & 
0 & 4 & 2 & 2 & 0 & 0 & 0 & 0 & 0 & 0 & 0 & 1 & 3 & 1 & 1 & 1 & 2 \\ $12$ & 2 
& 0 & 0 & 0 & 4 & 2 & 2 & 0 & 2 & 0 & 0 & 0 & 0 & 0 & 1 & 3 & 1 & 1 & 1 & 2 
\\ $1\ol{2}$ & 2 & 2 & 2 & 0 & 4 & 2 & 2 & 0 & 2 & 2 & 0 & 2 & 0 & 0 & 1 & 3 
& 1 & 1 & 1 & 2 \\ $1\ol{1}$ & 2 & 2 & 2 & 0 & 4 & 2 & 2 & 0 & 2 & 2 & 0 & 2 
& 0 & 0 & 1 & 3 & 1 & 1 & 1 & 2 \\ $22$ & 2 & 2 & 0 & 0 & 4 & 4 & 2 & 2 & 4 & 
0 & 2 & 0 & 0 & 0 & 1 & 3 & 3 & 1 & 1 & 2 \\ $2\ol{2}$ & 2 & 2 & 0 & 0 & 4 & 
4 & 2 & 2 & 4 & 0 & 2 & 0 & 0 & 0 & 1 & 3 & 3 & 1 & 1 & 2 \\ $2\ol{1}$ & 2 & 
2 & 2 & 2 & 4 & 4 & 2 & 2 & 4 & 2 & 2 & 2 & 2 & 2 & 1 & 3 & 3 & 1 & 1 & 2 \\ 
$\ol{2}\ol{2}$ & 4 & 4 & 4 & 2 & 4 & 4 & 4 & 2 & 4 & 4 & 2 & 4 & 2 & 0 & 3 & 
3 & 3 & 3 & 1 & 2 \\ $\ol{2}\ol{1}$ & 4 & 4 & 4 & 2 & 4 & 4 & 4 & 2 & 4 & 4 & 
2 & 4 & 2 & 2 & 3 & 3 & 3 & 3 & 1 & 2 \\ $\ol{1}\ol{1}$ & 4 & 4 & 4 & 4 & 4 & 
4 & 4 & 4 & 4 & 4 & 4 & 4 & 4 & 4 & 3 & 3 & 3 & 3 & 3 & 2 \\ $0$ & 1 & 1 & 1 
& 1 & 3 & 3 & 1 & 1 & 3 & 1 & 1 & 1 & 1 & 1 & 0 & 2 & 2 & 0 & 0 & 1 \\ $1$ & 
1 & 1 & 1 & 1 & 3 & 1 & 1 & 1 & 1 & 1 & 1 & 1 & 1 & 1 & 0 & 2 & 0 & 0 & 0 & 1 
\\ $2$ & 1 & 1 & 1 & 1 & 3 & 3 & 1 & 1 & 3 & 1 & 1 & 1 & 1 & 1 & 0 & 2 & 2 & 
0 & 0 & 1 \\ $\ol{2}$ & 3 & 3 & 3 & 1 & 3 & 3 & 3 & 1 & 3 & 3 & 1 & 3 & 1 & 1 
& 2 & 2 & 2 & 2 & 0 & 1 \\ $\ol{1}$ & 3 & 3 & 3 & 3 & 3 & 3 & 3 & 3 & 3 & 3 & 
3 & 3 & 3 & 3 & 2 & 2 & 2 & 2 & 2 & 1 \\ $\phi $ & 2 & 2 & 2 & 2 & 2 & 2 & 2 
& 2 & 2 & 2 & 2 & 2 & 2 & 2 & 1 & 1 & 1 & 1 & 1 & 0 \\  \hline 
\end{tabular}}\end{center}\end{table}

\begin{table}[tbp]
	\caption{$B^{1,1} \otimes B^{1,1}$ for $D^{(3)}_4$}
	\label{tab:D34_11}
    \newcommand{\lw}[1]{\smash{\lower1.6ex\hbox{#1}}}
	\begin{center}
	\begin{tabular}{c|cccccccc}	\hline
	\lw{$b$} & \multicolumn{8}{c}{$c$} \\

& $\phi $ & $\ol{1}$ & $\ol{2}$ & $\ol{3}$ & $0$ & $3$ & $2$ & $1$ \\  
\hline$\phi $ & 0 & 1 & 1 & 1 & 1 & 1 & 1 & 1 \\ $\ol{1}$ & 1 & 2 & 2 & 2 & 2 
& 2 & 2 & 2 \\ $\ol{2}$ & 1 & 1 & 2 & 2 & 2 & 2 & 2 & 2 \\ $\ol{3}$ & 1 & 1 & 
1 & 2 & 2 & 2 & 2 & 2 \\ $0$ & 1 & 0 & 1 & 1 & 1 & 2 & 2 & 2 \\ $3$ & 1 & 0 & 
1 & 1 & 1 & 2 & 2 & 2 \\ $2$ & 1 & 0 & 0 & 1 & 1 & 1 & 2 & 2 \\ $1$ & 1 & 0 & 
0 & 0 & 0 & 1 & 1 & 2 \\  \hline \end{tabular}\end{center}\end{table}

\begin{table}[tbp]
	\caption{$B^{1,2} \otimes B^{1,1}$ for $D^{(3)}_4$}
	\label{tab:D34_21}
    \newcommand{\lw}[1]{\smash{\lower1.6ex\hbox{#1}}}
	\begin{center}
	\begin{tabular}{c|cccccccc}	\hline
	\lw{$b$} & \multicolumn{8}{c}{$c$} \\

& $\phi $ & $\ol{1}$ & $\ol{2}$ & $\ol{3}$ & $0$ & $3$ & $2$ & $1$ \\  
\hline$\ol{1}\ol{1}$ & $\ol{1} \cdot \ol{1}_1$ & $\ol{1} \cdot 
\ol{1}\ol{1}_2$ & $\ol{1} \cdot \ol{2}\ol{1}_2$ & $\ol{1} \cdot 
\ol{3}\ol{1}_2$ & $\ol{1} \cdot 0\ol{1}_2$ & $\ol{1} \cdot 3\ol{1}_2$ & 
$\ol{1} \cdot 2\ol{1}_2$ & $\ol{1} \cdot 1\ol{1}_2$ \\ $\ol{3}\ol{3}$ & 
$\ol{3} \cdot \ol{3}_1$ & $\ol{3} \cdot \ol{3}\ol{1}_1$ & $\ol{3} \cdot 
0\ol{1}_1$ & $\ol{3} \cdot \ol{3}\ol{3}_2$ & $\ol{3} \cdot 0\ol{3}_2$ & 
$\ol{3} \cdot 3\ol{3}_2$ & $\ol{3} \cdot 2\ol{3}_2$ & $\ol{3} \cdot 
1\ol{3}_2$ \\ $20$ & $0 \cdot 2_1$ & $0 \cdot \ol{3}_0$ & $0 \cdot 0_0$ & $0 
\cdot 2\ol{3}_1$ & $0 \cdot 1\ol{3}_1$ & $0 \cdot 10_1$ & $0 \cdot 22_2$ & $0 
\cdot 12_2$ \\ $\ol{2}\ol{2}$ & $\ol{2} \cdot \ol{2}_1$ & $\ol{2} \cdot 
\ol{2}\ol{1}_1$ & $\ol{2} \cdot \ol{2}\ol{2}_2$ & $\ol{2} \cdot 
\ol{3}\ol{2}_2$ & $\ol{2} \cdot 0\ol{2}_2$ & $\ol{2} \cdot 3\ol{2}_2$ & 
$\ol{2} \cdot 2\ol{2}_2$ & $\ol{2} \cdot 1\ol{2}_2$ \\ $\ol{3}\ol{2}$ & 
$\ol{2} \cdot \ol{3}_1$ & $\ol{2} \cdot \ol{3}\ol{1}_1$ & $\ol{2} \cdot 
0\ol{1}_1$ & $\ol{2} \cdot \ol{3}\ol{3}_2$ & $\ol{2} \cdot 0\ol{3}_2$ & 
$\ol{2} \cdot 3\ol{3}_2$ & $\ol{2} \cdot 2\ol{3}_2$ & $\ol{2} \cdot 
1\ol{3}_2$ \\ $0\ol{2}$ & $\ol{2} \cdot 0_1$ & $\ol{2} \cdot \ol{1}_0$ & 
$\ol{2} \cdot 3\ol{1}_1$ & $\ol{2} \cdot 2\ol{1}_1$ & $\ol{2} \cdot 
1\ol{1}_1$ & $\ol{2} \cdot 30_2$ & $\ol{2} \cdot 20_2$ & $\ol{2} \cdot 10_2$ 
\\ $0\ol{3}$ & $\ol{3} \cdot 0_1$ & $\ol{3} \cdot \ol{1}_0$ & $\ol{3} \cdot 
3\ol{1}_1$ & $\ol{3} \cdot 2\ol{1}_1$ & $\ol{3} \cdot 1\ol{1}_1$ & $\ol{3} 
\cdot 30_2$ & $\ol{3} \cdot 20_2$ & $\ol{3} \cdot 10_2$ \\ $3\ol{3}$ & 
$\ol{3} \cdot 3_1$ & $\ol{3} \cdot \ol{2}_0$ & $\ol{3} \cdot 3\ol{2}_1$ & 
$\ol{3} \cdot 2\ol{2}_1$ & $\ol{3} \cdot 1\ol{2}_1$ & $\ol{3} \cdot 33_2$ & 
$\ol{3} \cdot 23_2$ & $\ol{3} \cdot 13_2$ \\ $30$ & $0 \cdot 3_1$ & $0 \cdot 
\ol{2}_0$ & $0 \cdot 3\ol{2}_1$ & $0 \cdot 2\ol{2}_1$ & $0 \cdot 1\ol{2}_1$ & 
$0 \cdot 33_2$ & $0 \cdot 23_2$ & $0 \cdot 13_2$ \\ $33$ & $3 \cdot 3_1$ & $3 
\cdot \ol{2}_0$ & $3 \cdot 3\ol{2}_1$ & $3 \cdot 2\ol{2}_1$ & $3 \cdot 
1\ol{2}_1$ & $3 \cdot 33_2$ & $3 \cdot 23_2$ & $3 \cdot 13_2$ \\ $23$ & $3 
\cdot 2_1$ & $3 \cdot \ol{3}_0$ & $3 \cdot 0_0$ & $3 \cdot 2\ol{3}_1$ & $3 
\cdot 1\ol{3}_1$ & $3 \cdot 10_1$ & $3 \cdot 22_2$ & $3 \cdot 12_2$ \\ $22$ & 
$2 \cdot 2_1$ & $2 \cdot \ol{3}_0$ & $2 \cdot 0_0$ & $2 \cdot 2\ol{3}_1$ & $2 
\cdot 1\ol{3}_1$ & $2 \cdot 10_1$ & $2 \cdot 22_2$ & $2 \cdot 12_2$ \\ 
$3\ol{2}$ & $\ol{2} \cdot 3_1$ & $3 \cdot \ol{2}\ol{1}_0$ & $3 \cdot 
\ol{2}\ol{2}_1$ & $3 \cdot \ol{3}\ol{2}_1$ & $3 \cdot 0\ol{2}_1$ & $\ol{2} 
\cdot 33_2$ & $\ol{2} \cdot 23_2$ & $\ol{2} \cdot 13_2$ \\ $2\ol{2}$ & 
$\ol{2} \cdot 2_1$ & $3 \cdot \ol{3}\ol{1}_0$ & $3 \cdot 0\ol{1}_0$ & $3 
\cdot \ol{3}\ol{3}_1$ & $3 \cdot 0\ol{3}_1$ & $3 \cdot 3\ol{3}_1$ & $\ol{2} 
\cdot 22_2$ & $\ol{2} \cdot 12_2$ \\ $2\ol{3}$ & $\ol{3} \cdot 2_1$ & $2 
\cdot \ol{3}\ol{1}_0$ & $2 \cdot 0\ol{1}_0$ & $2 \cdot \ol{3}\ol{3}_1$ & $2 
\cdot 0\ol{3}_1$ & $2 \cdot 3\ol{3}_1$ & $\ol{3} \cdot 22_2$ & $\ol{3} \cdot 
12_2$ \\ $\ol{2}\ol{1}$ & $\ol{1} \cdot \ol{2}_1$ & $\ol{2} \cdot 
\ol{1}\ol{1}_1$ & $\ol{1} \cdot \ol{2}\ol{2}_2$ & $\ol{1} \cdot 
\ol{3}\ol{2}_2$ & $\ol{1} \cdot 0\ol{2}_2$ & $\ol{1} \cdot 3\ol{2}_2$ & 
$\ol{1} \cdot 2\ol{2}_2$ & $\ol{1} \cdot 1\ol{2}_2$ \\ $\ol{3}\ol{1}$ & 
$\ol{1} \cdot \ol{3}_1$ & $\ol{3} \cdot \ol{1}\ol{1}_1$ & $\ol{3} \cdot 
\ol{2}\ol{1}_1$ & $\ol{1} \cdot \ol{3}\ol{3}_2$ & $\ol{1} \cdot 0\ol{3}_2$ & 
$\ol{1} \cdot 3\ol{3}_2$ & $\ol{1} \cdot 2\ol{3}_2$ & $\ol{1} \cdot 
1\ol{3}_2$ \\ $0\ol{1}$ & $\ol{1} \cdot 0_1$ & $0 \cdot \ol{1}\ol{1}_0$ & 
$\ol{3} \cdot \ol{2}\ol{2}_1$ & $\ol{3} \cdot \ol{3}\ol{2}_1$ & $\ol{3} \cdot 
0\ol{2}_1$ & $\ol{1} \cdot 30_2$ & $\ol{1} \cdot 20_2$ & $\ol{1} \cdot 10_2$ 
\\ $3\ol{1}$ & $\ol{1} \cdot 3_1$ & $0 \cdot \ol{2}\ol{1}_0$ & $0 \cdot 
\ol{2}\ol{2}_1$ & $0 \cdot \ol{3}\ol{2}_1$ & $0 \cdot 0\ol{2}_1$ & $\ol{1} 
\cdot 33_2$ & $\ol{1} \cdot 23_2$ & $\ol{1} \cdot 13_2$ \\ $2\ol{1}$ & 
$\ol{1} \cdot 2_1$ & $0 \cdot \ol{3}\ol{1}_0$ & $0 \cdot 0\ol{1}_0$ & $0 
\cdot \ol{3}\ol{3}_1$ & $0 \cdot 0\ol{3}_1$ & $0 \cdot 3\ol{3}_1$ & $\ol{1} 
\cdot 22_2$ & $\ol{1} \cdot 12_2$ \\ $1\ol{1}$ & $\ol{1} \cdot 1_1$ & $\ol{1} 
\cdot \phi _0$ & $0 \cdot 3\ol{1}_0$ & $0 \cdot 2\ol{1}_0$ & $0 \cdot 
1\ol{1}_0$ & $0 \cdot 30_1$ & $0 \cdot 20_1$ & $\ol{1} \cdot 11_2$ \\ 
$1\ol{2}$ & $\ol{2} \cdot 1_1$ & $\ol{2} \cdot \phi _0$ & $3 \cdot 3\ol{1}_0$ 
& $3 \cdot 2\ol{1}_0$ & $3 \cdot 1\ol{1}_0$ & $3 \cdot 30_1$ & $3 \cdot 20_1$ 
& $\ol{2} \cdot 11_2$ \\ $1\ol{3}$ & $\ol{3} \cdot 1_1$ & $\ol{3} \cdot \phi 
_0$ & $2 \cdot 3\ol{1}_0$ & $2 \cdot 2\ol{1}_0$ & $2 \cdot 1\ol{1}_0$ & $2 
\cdot 30_1$ & $2 \cdot 20_1$ & $\ol{3} \cdot 11_2$ \\ $10$ & $0 \cdot 1_1$ & 
$0 \cdot \phi _0$ & $2 \cdot 3\ol{2}_0$ & $2 \cdot 2\ol{2}_0$ & $2 \cdot 
1\ol{2}_0$ & $2 \cdot 33_1$ & $2 \cdot 23_1$ & $0 \cdot 11_2$ \\ $13$ & $3 
\cdot 1_1$ & $3 \cdot \phi _0$ & $1 \cdot 3\ol{2}_0$ & $1 \cdot 2\ol{2}_0$ & 
$1 \cdot 1\ol{2}_0$ & $1 \cdot 33_1$ & $1 \cdot 23_1$ & $3 \cdot 11_2$ \\ 
$12$ & $2 \cdot 1_1$ & $2 \cdot \phi _0$ & $2 \cdot 3_0$ & $1 \cdot 
2\ol{3}_0$ & $1 \cdot 1\ol{3}_0$ & $2 \cdot 13_1$ & $1 \cdot 22_1$ & $2 \cdot 
11_2$ \\ $11$ & $1 \cdot 1_1$ & $1 \cdot \phi _0$ & $1 \cdot 3_0$ & $1 \cdot 
2_0$ & $1 \cdot 10_0$ & $1 \cdot 13_1$ & $1 \cdot 12_1$ & $1 \cdot 11_2$ \\ 
$\phi $ & $\phi  \cdot \phi _0$ & $1 \cdot \ol{1}\ol{1}_0$ & $1 \cdot 
\ol{2}\ol{1}_0$ & $1 \cdot \ol{3}\ol{1}_0$ & $1 \cdot 0\ol{1}_0$ & $1 \cdot 
3\ol{1}_0$ & $1 \cdot 2\ol{1}_0$ & $1 \cdot 1\ol{1}_0$ \\ $\ol{1}$ & $\phi  
\cdot \ol{1}_0$ & $\phi  \cdot \ol{1}\ol{1}_1$ & $\phi  \cdot \ol{2}\ol{1}_1$ 
& $\phi  \cdot \ol{3}\ol{1}_1$ & $\phi  \cdot 0\ol{1}_1$ & $\phi  \cdot 
3\ol{1}_1$ & $\phi  \cdot 2\ol{1}_1$ & $\phi  \cdot 1\ol{1}_1$ \\ $\ol{2}$ & 
$\phi  \cdot \ol{2}_0$ & $3 \cdot \ol{1}\ol{1}_0$ & $\phi  \cdot 
\ol{2}\ol{2}_1$ & $\phi  \cdot \ol{3}\ol{2}_1$ & $\phi  \cdot 0\ol{2}_1$ & 
$\phi  \cdot 3\ol{2}_1$ & $\phi  \cdot 2\ol{2}_1$ & $\phi  \cdot 1\ol{2}_1$ 
\\ $\ol{3}$ & $\phi  \cdot \ol{3}_0$ & $2 \cdot \ol{1}\ol{1}_0$ & $2 \cdot 
\ol{2}\ol{1}_0$ & $\phi  \cdot \ol{3}\ol{3}_1$ & $\phi  \cdot 0\ol{3}_1$ & 
$\phi  \cdot 3\ol{3}_1$ & $\phi  \cdot 2\ol{3}_1$ & $\phi  \cdot 1\ol{3}_1$ 
\\ $0$ & $\phi  \cdot 0_0$ & $0 \cdot \ol{1}_0$ & $2 \cdot \ol{2}\ol{2}_0$ & 
$2 \cdot \ol{3}\ol{2}_0$ & $2 \cdot 0\ol{2}_0$ & $\phi  \cdot 30_1$ & $\phi  
\cdot 20_1$ & $\phi  \cdot 10_1$ \\ $3$ & $\phi  \cdot 3_0$ & $3 \cdot 
\ol{1}_0$ & $1 \cdot \ol{2}\ol{2}_0$ & $1 \cdot \ol{3}\ol{2}_0$ & $1 \cdot 
0\ol{2}_0$ & $\phi  \cdot 33_1$ & $\phi  \cdot 23_1$ & $\phi  \cdot 13_1$ \\ 
$2$ & $\phi  \cdot 2_0$ & $2 \cdot \ol{1}_0$ & $2 \cdot \ol{2}_0$ & $1 \cdot 
\ol{3}\ol{3}_0$ & $1 \cdot 0\ol{3}_0$ & $1 \cdot 3\ol{3}_0$ & $\phi  \cdot 
22_1$ & $\phi  \cdot 12_1$ \\ $1$ & $\phi  \cdot 1_0$ & $1 \cdot \ol{1}_0$ & 
$1 \cdot \ol{2}_0$ & $1 \cdot \ol{3}_0$ & $1 \cdot 0_0$ & $1 \cdot 30_0$ & $1 
\cdot 20_0$ & $\phi  \cdot 11_1$ \\  \hline 
\end{tabular}\end{center}\end{table}

}

\end{document}